%% file: fluid-shell-preprint.tex
\documentclass[runningheads,reqno]{amsart}
\usepackage{graphicx}
\usepackage{mathptmx}

\usepackage{amsmath}
\usepackage{amssymb}
\usepackage{amsfonts}
\usepackage{amsxtra}
\usepackage{color}

\usepackage{srcltx}
\numberwithin{equation}{section}

\usepackage[rose,operator,hyperref]{paper_diening}

\newtheorem{Theorem}[equation]{Theorem}
\newtheorem{Lemma}[equation]{Lemma}
\newtheorem{Proposition}[equation]{Proposition}
\newtheorem{Corollary}[equation]{Corollary}
 \newtheorem{Definition}[equation]{Definition}
\newtheorem{Remark}[equation]{Remark}
\newcommand{\R}{\mathcal{R}}
\newcommand{\F}{\mathcal{F}}
\newcommand{\M}{\mathcal{M}}
\newcommand{\T}{\mathcal{T}}
\newcommand{\bnue}{\boldsymbol{\nu}_{\eta(t)}}
\newcommand{\iot}{\int_{\Omega_{\eta(t)}}}

\newcommand{\iorte}{\int_{\Omega_{\mathcal{R}_\epsilon\eta(t)}}}
\newcommand{\ios}{\int_{\Omega_{\eta(s)}}}

\newcommand{\iortd}{\int_{\Omega_{\mathcal{R}\delta(t)}}}
\newcommand{\iortdn}{\int_{\Omega_{\mathcal{R}\delta_n(t)}}}
\newcommand{\idot}{\int_{\pa\Omega_{\eta(t)}\setminus\Gamma}}
\newcommand{\im}{\int_M}

\newcommand{\onorm}[1]{\norm{#1}_{L^2(\Omega_{\eta(t)})}}
\newcommand{\mnorm}[1]{\norm{#1}_{L^2(M)}}
\newcommand{\oet}{\Omega_{\eta(t)}}

\newcommand{\beweis}{\noindent{\bf Proof:}\ }

\DeclareMathOperator{\dv}{div}

\DeclareMathOperator{\tr}{tr_\eta}
\DeclareMathOperator{\trt}{tr_{\tilde\eta}}
\DeclareMathOperator{\trent}{tr_{\tilde\eta_n}}
\DeclareMathOperator{\trnormal}{tr^n_\eta}
\DeclareMathOperator{\trnormaln}{tr^n_{\eta_0}}
\DeclareMathOperator{\trnormale}{tr^n_{\R_\epsilon\eta_0}}
\DeclareMathOperator{\trnormaleo}{tr^n_{\R\eta_0}}
\DeclareMathOperator{\trnormald}{tr^n_{\R\delta}}
\DeclareMathOperator{\tren}{tr_{\eta_n}}

\DeclareMathOperator{\trre}{tr_{\mathcal{R}_\epsilon\eta}}

\DeclareMathOperator{\trrd}{tr_{\mathcal{R}\delta}}
\DeclareMathOperator{\trrdt}{tr_{\mathcal{R}\delta(t)}}
\DeclareMathOperator{\trd}{tr_{\delta}}
\DeclareMathOperator{\trrdn}{tr_{\mathcal{R}\delta_n}}

\DeclareMathOperator{\grad}{grad}
\DeclareMathOperator{\id}{id}

\DeclareMathOperator{\spann}{span}
\DeclareMathOperator{\inn}{int}
\DeclareMathOperator{\supp}{supp}

\begin{document}

\title[Interaction of a Newtonian fluid with a Koiter shell]{Global
  weak solutions for an incompressible, Newtonian fluid interacting
  with a linearly elastic Koiter shell}


\author{Daniel Lengeler \and Michael R\r u\v zi\v cka }


\address{Daniel Lengeler,
  Fakult{\"a}t f{\"u}r Mathematik\\
  Universit{\"a}t Regensburg\\
  Universit{\"a}tsstr.~31, 93040 Regensburg, Germany \\
} 
\email{daniel.lengeler@mathematik.uni-regensburg.de}

\address{Michael R\r u\v zi\v cka \\
  Mathematisches Institut\\
  Universit{\"a}t Freiburg \\
  Eckerstr.~1, 79104 Freiburg, Germany\\
}
\email{rose@mathematik.uni-freiburg.de} 


\maketitle

\begin{abstract}
  In this paper we analyze the interaction of an incompressible 
  Newtonian fluid with a linearly elastic Koiter shell whose motion is
  restricted to transverse displacements. The middle surface of the
  shell constitutes the mathematical boundary of the three-dimensional
  fluid domain. We show that weak solutions exist as long as the
  magnitude of the displacement stays below some (possibly large)
  bound which is determined by the geometry of the undeformed shell.
\end{abstract}

\section{Introduction}
\label{intro}
Fluid-solid interaction problems involving moving interfaces have been
studied intensively during the last two decades. The interaction with
elastic solids has proven to be particularly difficult, due to
apparent regularity incompatibilities between the parabolic fluid
phase and the hyperbolic or dispersive solid phase. Therefore, the
first results in this direction involved damped or regularized
elasticity laws; see for instance \cite{b38,b37,b27,b61} and the references therein.  In
particular, in \cite{b27} the global-in-time existence of weak
solutions for the interaction of an incompressible, Newtonian fluid
with a damped Kirchhoff-Love plate is shown. The first results
involving the classical elasticity laws without any regularisation are
\cite{b40,b41} where the local-in-time existence of a unique solution
for an incompressible, Newtonian fluid interacting with linear and
quasilinear three-dimensional elasticity, respectively, is shown.
Here, the data is assumed to be very smooth. Similar results are shown
in \cite{b60,b42} where the interaction of an incompressible,
Newtonian fluid with two-dimensional nonlinearly elastic structures,
so-called shells, constituting the mathematical boundary of the
three-dimensional fluid domain, is analyzed. However, in both papers
inertia had to be neglected in the three-dimensional case, leading to
parabolic-elliptic systems. In \cite{b27} the added damping term is
mainly needed in the proof of compactness of sequences of approximate
weak solutions. The subsequent work \cite{b28} introduces a technique
which allows to get rid of this damping term. In both papers the
compactness argument relies on the strategy to test the system with
difference quotients in time. Unfortunately, it seems that this
argument rests heavily upon the simple flat geometry of the elastic
plate. In particular, it seems to be very difficult to apply it to the
case of a general shell which is the problem we are interested in. For
this reason, in the present paper we introduce a completely new
approach to this compactness problem which is closely related to the
classical Aubin-Lions theorem.

The present paper summarizes results from Daniel Lengeler's Ph.D. thesis \cite{b62}. It is organized as follows. In Subsection 1.1 we introduce Koiter's energy for elastic shells, in Subsection 1.2 we introduce the coupled fluid-shell system, and in Subsection 1.3 we derive formal a-priori estimates for this system. In Section 2 we derive some results concerning domains with non-Lipschitz boundaries. Then, in Section 3 we state the main result of the paper. The rest of the section is devoted to the proof of this result. In Subsection 3.1 we give the proof of compactness of sequences of weak solutions. Subsequently, in Subsection 3.2 we analyse a decoupled variant of our original system, while in Subsection 3.3 we apply a fixed-point argument to this decoupled system. In Subsection 3.4 we conclude the proof by letting the regularisation parameter, which we introduced earlier, tend to zero. Finally, some further results and technical computations can be found in the appendix.

\subsection{Koiter's energy}
Throughout the paper, let $\Omega\subset\setR^3$ be a bounded,
non-empty domain of class $C^4$ with outer unit normal $\bnu$. We
denote by $g$ and $h$ the first and the second fundamental form of
$\pa\Omega$, induced by the ambient Euclidean space, and by
$dA$ the surface measure of $\pa\Omega$. Furthermore, let
$\Gamma\subset\pa\Omega$ be a union of domains of class $C^{1,1}$
having non-trivial intersection with all connected components of
$\pa\Omega$. We set $M:=\pa\Omega\setminus\Gamma$; note that $M$
is compact. Let $\pa\Omega$ represent the middle surface of an
elastic shell of thickness $2\,\epsilon_0>0$ in its rest state, where
$\epsilon_0$ is taken to be small compared to the reciprocal of the
principal curvatures. Furthermore, we assume that the elastic shell
consists of a homogeneous, isotropic material whose linear elastic
behavior may be characterized by the Lam{\'e} constants $\lambda$ and
$\mu$. We shall describe deformations of the middle surface and hence
of the shell by a vector field $\boldsymbol{\eta}:\ M
\rightarrow\setR^3$ vanishing at the boundary $\pa M$. Hence, we
assume the part $\Gamma$ of the middle surface to be fixed. Letting
$g(\boldsymbol{\eta})$ and $h(\boldsymbol{\eta})$ denote the pullback
of the first and second fundamental form, respectively, of the
deformed middle surface, the elastic energy of the deformation may be
modeled by \emph{Koiter's energy for a nonlinearly elastic shell}
\[
K_N(\boldsymbol{\eta})=\frac{1}{2}\int_M \epsilon_0\,\langle C,
\Sigma(\boldsymbol{\eta})\otimes\Sigma(\boldsymbol{\eta})\rangle +
\frac{\epsilon_0^3}{3}\,\langle C,
\Xi(\boldsymbol{\eta})\otimes\Xi(\boldsymbol{\eta})\rangle\ dA;
\]
see \cite{b54}, \cite{b20}, \cite{b45}, \cite{b23}, and the references
therein.  Here
\begin{equation*}
  \begin{aligned}
    C_{\alpha\beta\gamma\delta}:=\frac{4\lambda\mu}{\lambda+2\mu} 
    \,g_{\alpha\beta}\,g_{\gamma\delta}
    + 2\mu\,(g_{\alpha\gamma}\,g_{\beta\delta} +
    g_{\alpha\delta}\,g_{\beta\gamma})
 \end{aligned}
\end{equation*}
is the elasticity tensor of the shell, and
$\Sigma(\boldsymbol{\eta}):=1/2\, (g(\boldsymbol{\eta})-g)$ and
$\Xi(\boldsymbol{\eta}):=h(\boldsymbol{\eta})-h$ denote the
differences of the first and the second fundamental forms,
respectively. In \cite{b20} this energy is derived from
three-dimensional elasticity under the additional assumptions of small
strains and plane stresses parallel to the middle surface. The part of
the energy scaling with $\epsilon_0$ is the membrane energy, the part
scaling with $\epsilon_0^3$ is the bending energy. For a rigorous
justification of this energy we refer to \cite{b45,b21,b22}. Motivated
by the approach in \cite{b27,b28}, we restrict the deformation of the
middle surface to displacements along the unit normal field $\bnu$ of
$M$. Hence, the displacement $\boldsymbol{\eta}=:\eta\,\bnu$ may be
described by a scalar function $\eta$.  Furthermore, we linearize the
dependence of the strain tensors $\Sigma(\eta)$ and $\Xi(\eta)$ on
$\eta$ at $\eta=0$, resulting in the tensor fields\footnote{We denote
  the covariant derivative of tensor fields on Riemannian manifolds
  (including subsets of $\setR^n$) by $\nabla$. Furthermore, we write
  $\Delta$ for the associated Laplace operator.}
\begin{equation*}
 \begin{aligned}
  \sigma(\eta)=-h\, \eta,\qquad  \xi(\eta)=\nabla^2 \eta - k\, \eta,
 \end{aligned}
\end{equation*}
where $k_{\alpha\beta}:=h_\alpha^\sigma\, h_{\sigma\beta}$; see
Theorem 4.2-1 and Theorem 4.2-2 in \cite{b23}. This way, we obtain
\emph{Koiter's energy for a linearly elastic shell and transverse
  displacements}
\begin{equation*}
 \begin{aligned}
  K(\eta)=K(\eta,\eta)=\frac{1}{2}\int_M \epsilon_0\,\langle
C,
\sigma(\eta)\otimes\sigma(\eta)\rangle + \frac{\epsilon_0^3}{3}\,\langle
C,
\xi(\eta)\otimes\xi(\eta)\rangle\ dA. 
 \end{aligned}
\end{equation*}
$K$ is a quadratic form in $\eta$ which is coercive on $H^2_0(M)$,
i.e. there exists a constant $c_0$ such that
\begin{equation}\label{eqn:kkoerziv}
 \begin{aligned}
K(\eta)\ge c_0\,\norm{\eta}_{H^2_0(M)}^2;
  \end{aligned}
\end{equation}
see the proof of Theorem 4.4-2 in \cite{b23}. Using integration by
parts and taking into account some facts from Riemannian geometry one
can show that the $L^2$ gradient of this energy has the form
\[\grad_{L^2}K(\eta)=\epsilon_0^3\,\frac{8\mu(\lambda+\mu)}{
  3(\lambda+2\mu)}\,\Delta^2\eta+B\eta,\] where $B$ is a second order
differential operator which vanishes on flat parts of $M$, i.e. where
$h=0$. The details can be found in \cite{b62}. Thus, we obtain a generalization of the
linear Kirchhoff-Love plate equation for transverse displacements; cf.
for instance \cite{b44}.  By Hamilton's principle, the displacement
$\eta$ of the shell must be a stationary point of the action
functional
\[
\mathcal{A}(\eta)=\int_I \epsilon_0\rho_S\int_M
\frac{(\pa_t\eta(t,\cdot))^2}{2}\ dA-K(\eta(t,\cdot))\ dt,
\]
where $I:=(0,T)$, $T>0$. Here we assume
that the mass density of $M$ may be described by a constant
$\epsilon_0\rho_S$. Hence, the integrand with respect to time is the
difference of the kinetic and the potential energy of the shell. The
corresponding Euler-Lagrange equation is
\begin{equation*}
  \begin{aligned}
    \epsilon_0\rho_S\,\pa^2_t\eta + \grad_{L^2}K(\eta)=0\text{ in }
    I\times M.
  \end{aligned}
\end{equation*}
Note that this equation is dispersive but not hyperbolic since the
principal part factorizes into two Schr{\"o}dinger operators, i.e.
\[
\pa^2_t + \pa^4_x=(i\pa_t+\pa^2_x)(-i\pa_t+\pa^2_x).
\]
Hence, we have an infinite speed of propagation of disturbances. For
the sake of a simple notation we assume $\epsilon_0\rho_S=1$.

\subsection{Statement of the problem}\label{statement}

We denote by $\Omega_{\eta(t)}$, $t\in I$, the deformed domain
(cf.~\eqref{eq:def-om-eta}) and by \[\Omega_\eta^I:=\bigcup_{t\in I}\,
\{t\}\times \Omega_{\eta(t)}\] the deformed spacetime cylinder. Now,
let us suppose that the variable domain $\Omega_\eta$ is filled by a
viscous, incompressible Newtonian fluid whose velocity field $\bu$ and
pressure field $\pi$ is described by the \emph{Navier-Stokes
  equations}, i.e.
\begin{equation}\label{eqn:fluid}
 \begin{aligned}
   \rho_F\big(\pa_t \bu + (\bu\cdot\nabla)\bu\big) - \dv \big(2\sigma D\bu
   - \pi\id\big) &= \ff &&\mbox{ in }
   \Omega_{\eta}^I, \\
   \dv \bu &= 0 &&\mbox{ in }  \Omega_{\eta}^I,\\
   \bu(\,\cdot\,,\,\cdot\, + \eta\,\bnu) &= \pa_t\eta\,\bnu &&\mbox{ on } I\times M,\\
   \bu &= 0 &&\mbox{ on } I\times\Gamma.
\end{aligned}
\end{equation}
Here, $\id$ denotes the $3\times 3$ unit matrix and $\ff$ is a given
force density. Furthermore, we set $(\bu\cdot\nabla)\bu:=u^i\,\pa_i\bu.$
For the sake of a simple exposition we set the constant fluid density
$\rho_F$ and the constant dynamical viscosity $\sigma$ equal to $1$. Then,
$2D\bu - \pi\id$ is the Cauchy stress tensor and $2D\bu$ its viscous
part.  Due to the divergence-free constraint we have $\dv
2D\bu=\Delta\bu$.  \eqref{eqn:fluid}$_{3,4}$ is the no-slip condition
in the case of a moving boundary, i.e. the velocity of the fluid at
the boundary equals the velocity of the boundary. The force exerted by
the fluid on the boundary is given by the evaluation of the stress
tensor at the deformed boundary in the direction of the inner normal
$-\bnu_{\eta(t)}$, i.e. by
\begin{equation}\label{eqn:kraft}
 \begin{aligned}
-2D\bu(t,\cdot)\,\bnu_{\eta(t)} + \pi(t,\cdot)\,\bnu_{\eta(t)}.  
 \end{aligned}
\end{equation}
Thus, the equation for the displacement of the shell takes the form
\begin{equation}\label{eqn:shell}
\begin{aligned}
  \pa^2_{t} \eta + \grad_{L^2}K(\eta)  &= g + {\bf F}\cdot\bnu
  &&\text{ in } I\times M, \\ 
  \eta=0,\ \nabla \eta &= 0 &&\text{ on } I\times\pa M,
\end{aligned}
\end{equation}
where $g$ is an additional force density and
\[{\bf F}(t,\cdot) := \big(-2D\bu(t,\cdot)\,\bnu_{\eta(t)} +
\pi(t,\cdot)\,\bnu_{\eta(t)}\big)\circ\Phi_{\eta(t)}\ |\det
d\Phi_{\eta(t)}|\] with
$\Phi_{\eta(t)}:\pa\Omega\rightarrow\pa\Omega_{\eta(t)},\ q\mapsto
\eta(t,q)\,\bnu(q)$. Finally, we specify initial values
\begin{equation}\label{eqn:data}
 \begin{aligned}
   \eta(0,\cdot)=\eta_0,\ \pa_t\eta(0,\cdot)=\eta_1 \ \text{ in }M
   \quad\text{ and }\quad \bu(0,\cdot)=\bu_0 \ \text{ in
   }\Omega_{\eta_0}.
 \end{aligned}
\end{equation}
In the following we will analyse the system \eqref{eqn:fluid},
\eqref{eqn:shell}, \eqref{eqn:data}. 

\subsection{Formal a-priori estimates}\label{subsec:apriori}
Let us now formally derive energy estimates for this
parabolic-dispersive system. To this end, we multiply
\eqref{eqn:fluid}$_1$ by $\bu$, integrate the resulting identity over
$\Omega_{\eta(t)}$, and obtain after integrating by parts the stress
tensor\footnote{For the sake of a better readability we suppress the
  dependence of the unknown on the independent variables, e.g~we write
    $\bu=\bu(t,\cdot)$.}
\begin{align}
   &\iot \pa_t\bu\cdot\bu\ dx + \iot (\bu\cdot\nabla)\bu\cdot\bu\ dx
   \label{eqn:mult}
   \\
   &\quad= -\iot |\nabla\bu|^2\ dx+ \iot \ff\cdot\bu\ dx +\idot (2D\bu\,
   \bnue - \pi\,\bnue)\cdot \bu\ dA_{\eta(t)}.\notag
\end{align}
Here, $dA_{\eta(t)}$ denotes the surface measure of the deformed boundary $\pa\Omega_{\eta(t)}$. Concerning the viscous stress tensor we used Korn's identity
\[2\int_{\Omega_{\eta(t)}}D\bu:D\bu\
dx=2\int_{\Omega_{\eta(t)}}D\bu:\nabla\bu\
dx=\int_{\Omega_{\eta(t)}}\nabla\bu:\nabla\bu\ dx.\] For a (formal) proof of
the second equality see Remark \ref{bem:korn}. Taking into account that
\begin{equation}\label{eqn:wirbelid}
 \begin{aligned}
   \iot (\bu\cdot\nabla)\bu\cdot\bu\ dx=-\iot
   (\bu\cdot\nabla)\bu\cdot\bu\ dx + \int_{\pa\Omega_{\eta(t)}}
   \bu\cdot\bnu_{\eta(t)}|\bu|^2\ dA_{\eta(t)}
 \end{aligned}
\end{equation}
we may apply Reynold's transport theorem \ref{theorem:reynolds} to the first
two integrals in \eqref{eqn:mult} to obtain
\begin{equation}\label{eqn:energiefluid}
 \begin{aligned}
  \frac{1}{2}\frac{d}{dt} \iot |\bu|^2\ dx = &-\iot
|\nabla\bu|^2\ dx+ \iot \ff\cdot\bu\ dx\\
&+\idot
(2D\bu\, \bnue - \pi\,\bnue)\cdot \bu\ dA_{\eta(t)}.
 \end{aligned}
\end{equation}
Multiplying \eqref{eqn:shell}$_1$ by $\pa_t\eta$, integrating the
resulting identity over $M$, integrating by parts, and using the fact
that $(\grad_{L^2}K(\eta),\pa_t\eta)_{L^2}=2K(\eta,\pa_t\eta)$ we
obtain
\begin{equation}\label{eqn:energieshell}
 \begin{aligned}
   \frac{1}{2}\frac{d}{dt} \im |\pa_t\eta|^2\ dA + \frac{d}{dt}K(\eta)
   = \im g\, \pa_t\eta\ dA + \im {\bf F}\cdot\bnu\, \pa_t\eta\ dA.
 \end{aligned}
\end{equation}
Adding \eqref{eqn:energiefluid} and \eqref{eqn:energieshell}, taking
into account the definition of ${\bf F}$, \eqref{eqn:fluid}$_3$, and
applying a change of variables to the boundary integral we obtain the
energy identity
\begin{equation}\label{eqn:energiesatz}
 \begin{aligned}
   \frac{1}{2}\frac{d}{dt} \iot |\bu|^2\ dx& + \frac{1}{2}\frac{d}{dt}
   \im |\pa_t\eta|^2\ dA +
   \frac{d}{dt}K(\eta)\\
   & = - \iot |\nabla\bu|^2\ dx + \iot \ff\cdot\bu\ dx + \im g\,
   \pa_t\eta\ dA.
 \end{aligned}
\end{equation}
In view of \eqref{eqn:kkoerziv} an application of Gronwall's lemma
gives
\begin{align}
  &\onorm{\bu(t,\cdot)}^2 +
  \int_0^t\norm{\nabla\bu(s,\cdot)}^2_{L^2(\Omega_{\eta(s)})}\ ds +
  \mnorm{\pa_t\eta(t,\cdot)}^2 + \norm{\eta(t,\cdot)}_{H^2(M)}^2.\notag \\
  &\hspace{1cm}\le c\, e^t\,\Big(\norm{\bu_0}_{L^2(\Omega_{\eta_0})}^2 +
  \mnorm{\eta_1}^2 +
  \norm{\eta_0}_{H^2(M)}^2\label{ab:apriori} \\
&\hspace{4cm} + \int_0^t
  \norm{\ff(s,\cdot)}_{L^2(\Omega_{\eta(s)})}^2 + \mnorm{g(s,\cdot)}^2 ds\Big).\notag
 \end{align}
Hence we have
\begin{equation*}
 \begin{aligned}
 \norm{\eta}_{W^{1,\infty}(I,L^2(M))\cap L^\infty(I,H_0^2(M))} +
\norm{\bu}_{L^\infty(I,L^2(\Omega_{\eta(t)}))\cap L^2(I,H^1(\oet))}\le c(T,\text{data}).
 \end{aligned}
\end{equation*}
We shall construct weak solutions in this regularity class. In view of
the embedding $H^2(\pa\Omega)\embedding C^{0,\theta}(\pa\Omega)$ for
$\theta<1$ this implies that the boundary of our variable domain will
be the graph of a H{\"o}lder continuous function which, in general, is not
Lipschitz continuous. Due to this fact we need to take a closer look
at domains with non-Lipschitz boundaries in the next section.

\section{Variable domains}
\label{sec:1}

We denote by $S_{\alpha}$, $\alpha>0$, the open set of points in
$\setR^3$ whose distance from $\pa\Omega$ is less than $\alpha$. It's a
well known fact from elementary differential geometry, see for instance \cite{b1}, that there
exists a maximal $\kappa>0$ such that the mapping
\begin{equation*}
 \begin{aligned}
\Lambda: \partial\Omega\times (-\kappa,\kappa)\rightarrow
S_{\kappa},\
(q,s)\mapsto q + s\,\bnu(q)  
 \end{aligned}
\end{equation*}
is a $C^3$ diffeomorphism. For the inverse $\Lambda^{-1}$ we shall
write $x\mapsto(q(x),s(x))$; cf. Figure \ref{Bild1}. Note that $\kappa$ is not necessarily small; if $\Omega$ is the ball of radius $R$, then $\kappa=R$. Note furthermore that
$\Lambda$ has singularities on the boundary of its domain.
\begin{figure}[h]
\centering
\input{Bild1b.pstex_t}
\caption{}
\label{Bild1}
\end{figure}
Let $B_\alpha:=\Omega\cup S_{\alpha}$ for $0<\alpha<\kappa$. The
mapping $\Lambda(\,\cdot\, ,\alpha):\pa\Omega\rightarrow\pa B_\alpha$
is a $C^3$ diffeomorphism as well.  Hence, $B_\alpha$ is a bounded
domain with $C^3$ boundary. For a continuous function
$\eta:\pa\Omega\rightarrow (-\kappa,\kappa)$ we set
\begin{equation}
 \begin{aligned}
  \Omega_\eta:=\Omega\setminus S_\kappa\ \cup \{x\in S_\kappa\ |\
  s(x)<\eta(q(x))\};
 \end{aligned}\label{eq:def-om-eta}
\end{equation}
cf. Figure \ref{Bild1}. $\Omega_\eta$ is an open set. For $\eta\in
C^k(\pa\Omega)$, $k\in\{1,2,3\}$ we denote by $\bnu_\eta$ and
$dA_\eta$ the outer unit normal and the surface measure of
$\pa\Omega_\eta$, respectively. Now, let us construct a homeomorphism
from $\overline\Omega$ to $\overline{\Omega_\eta}$, the so called
Hanzawa transform. The details of this construction may be found in
\cite{b62}. Let the real-valued function $\beta\in C^\infty(\setR)$
be $0$ in a neighborhood of $-1$ and $1$ in a neighborhood of $0$. For
a continuous function $\eta:\pa\Omega\rightarrow (-\kappa,\kappa)$ we
define the mapping $\Psi_\eta:
\overline\Omega\rightarrow\overline{\Omega_\eta}$ in
$S_\kappa\cap\overline\Omega$ by
\begin{equation}\label{eqn:defpsi}
  \begin{aligned}
    &x\mapsto x + \bnu(q(x))\,
    \eta(q(x))\,\beta(s(x)/\kappa)\\
    &\hspace{0.6cm}=q(x)+\bnu(q(x))\big(s(x)+\eta(q(x))\,\beta(s(x)/\kappa)\big).
  \end{aligned}
\end{equation}
In $\Omega\setminus S_\kappa$  let $\Psi_\eta$ be the identity. For
this mapping to have an inverse (which is continuously differentiable
provided that $\eta$ is sufficiently regular) we need to assume that
$|\beta'(s)|<\kappa/|\eta(q)|$ for all $s\in [-1,0]$ and all
$q\in\pa\Omega$. This is possible since
$\norm{\eta}_{L^\infty(\pa\Omega)}<\kappa$. Then, it's
not hard to see that $\Psi_\eta$ is a homeomorphism and, provided that
$\eta\in C^k(\pa\Omega)$ with $k\in\{1,2,3\}$, even a $C^k$
diffeomorphism. Similarly, the homeomorphism
\[\Phi_\eta:=\Psi_\eta|_{\pa\Omega}:\pa\Omega\rightarrow\pa\Omega_\eta,\
q\mapsto q+\eta(q)\,\bnu(q)\] with inverse $x\mapsto q(x)$ is a
$C^k$ diffeomorphism provided that $\eta\in C^k(\pa\Omega)$,
$k\in\{1,2,3\}$. Furthermore, one can show that the components of the
Jacobians of $\Psi_\eta$, $\Phi_\eta$, and their inverses have the
form
\begin{equation}\label{eqn:funk}
  \begin{aligned}
    b_0+\sum_i b_i\ (\pa_i\eta)\circ q
  \end{aligned}
\end{equation}
for bounded, continuous functions $b_0,b_i$ whose supports are contained in $S_\kappa$. Let
\begin{equation}\label{eqn:taueta}
  \begin{aligned}
    \tau(\eta):=\left\{\begin{array}{cl}
        (1-\norm{\eta}_{L^\infty(\pa\Omega)}/\kappa)^{-1} &
        \quad \text{, if }\norm{\eta}_{L^\infty(\pa\Omega)}<\kappa\\
        \infty & \quad \text{, else}\end{array}\right..
 \end{aligned}
\end{equation}
For $\tau(\eta)\rightarrow\infty$ some of the functions $b_0$, $b_i$
blow up and the distance between their supports and the boundary of
$S_\kappa$ tends to zero. In particular, the mapping $q$ in
\eqref{eqn:funk} is evaluated near its singularities for large
$\tau(\eta)$. We conclude that there are two sources of singularities
of $\Psi_\eta$ and $\Phi_\eta$: a singularity of $\eta$ and the
maximum of $|\eta|$ getting close to $\kappa$. This is why the
continuity constants of the linear mappings between function spaces
constructed below will always depend on $\tau(\eta)$. Furthermore,
note that $\Psi_\eta$ depends on the cutoff function $\beta$ which in
turn may be chosen to depend only on $\tau(\eta)$. Whenever we will be
dealing with sequences $(\eta_n)$ with $\sup_n\tau(\eta_n)<\infty$ we
will tacitly assume $\beta$ to be independent of $n$.

\begin{Remark}\label{bem:tdelta}
For $\eta\in C^2(\pa\Omega)$ with $\norm{\eta}_{L^\infty(\pa\Omega)}<\kappa$ and
$\bphi:\Omega\rightarrow\setR^3$ we denote by $\T_\eta\bphi$ the \emph{pushforward} of $(\det
d\Psi_\eta)^{-1} \bphi$ under $\Psi_\eta$, i.e. 
\[\T_\eta\bphi:=\big(d\Psi_\eta\, (\det
d\Psi_\eta)^{-1}\bphi\big)\circ\Psi^{-1}_\eta.\]
The mapping $\T_\eta$ with the inverse
$\T_\eta^{-1}\bphi:=\big(d\Psi_\eta^{-1}\, (\det
d\Psi_\eta)\,\bphi\big)\circ\Psi_\eta$
obviously defines isomorphisms between the Lebesgue and Sobolev spaces on $\Omega$ respectively
$\Omega_\eta$ as long as the order of differentiability is not larger than $1$. Moreover, the mapping preserves vanishing
boundary values. 

The diffeomorphism $\Psi_\eta$ is an isometry from 
$\overline\Omega$ with the Riemannian metric $h:=(d\Psi_v)^T d\Psi_\eta$ to $\overline{
\Omega_\eta}$ with the Euclidean metric. From the naturality of the Levi-Civita connection with respect to isometries as well as the
identities $\sqrt{\det h_{\beta\gamma}}=\det d\Psi_t$ we infer that
\[(\dv \T_\eta\bphi)\circ\Psi_\delta=\dv_h((\det d\Psi_\eta)^{-1}\bphi)=(\det
d\Psi_\eta)^{-1}\dv\bphi.\]
Thus, $\T_\delta$ preserves the divergence-free property and hence defines isomorphisms between the corresponding function spaces
on $\Omega$ and $\Omega_\eta$, respectively.
\end{Remark}

A bi-Lipschitz mapping of domains induces isomorphisms of the
corresponding $L^p$ and $W^{1,p}$ spaces. For $\eta\in H^2(\pa\Omega)$
the mapping $\Psi_\eta$ is barely not bi-Lipschitz, due to the embedding
$H^2(\pa\Omega)\embedding C^{0,\theta}(\pa\Omega)$ for $\theta<1$.
Hence a small loss, made quantitative in the next lemma, will occur.
\begin{Lemma}\label{lemma:psi}
  Let $1<p\le\infty$ and $\eta\in H^2(\pa\Omega)$ with
  $\norm{\eta}_{L^\infty(\pa\Omega)}<\kappa$.  Then the linear mapping
  $v\mapsto v\circ\Psi_\eta$ is continuous from $L^p(\Omega_\eta)$ to
  $L^r(\Omega)$ and from $W^{1,p}(\Omega_\eta)$ to $W^{1,r}(\Omega)$
  for all $1\le r<p$. The analogous claim with $\Psi_\eta$ replaced by
  $\Psi_\eta^{-1}$ is true. The continuity constants depend only on
  $\Omega$, $p$, $\norm{\eta}_{H^2(\pa\Omega)}$, $\tau(\eta)$, and
  $r$; they stay bounded as long as $\norm{\eta}_{H^2(\pa\Omega)}$ and
  $\tau(\eta)$ stay bounded.
\end{Lemma}
\beweis Without loss of generality we may assume $p<\infty$. Let us
approximate $\eta$ by functions $(\eta_n)\subset C^2(\pa\Omega)$ in
$H^2(\pa\Omega)\embedding C(\pa\Omega)$. Due to \eqref{eqn:funk} and
\[
H^2(\pa\Omega)\embedding
W^{1,s}(\pa\Omega)\quad  \text{ for all }1\le s<\infty
\]
the components of the Jacobian of $\Psi_{\eta_n}^{-1}$, and hence its
Jacobian determinant, are bound\-ed in $L^s(\Omega_{\eta_n})$ for each
$s<\infty$, the bound depending on $\tau(\eta_n)$. Hence, by a
change of variables and H{\"o}lder's inequality we obtain for $v\in
C_0^\infty(\setR^3)$, $r<p$, and $1/r=1/p+1/s$ that\footnote{We denote the differential of a mapping $f$ between manifolds (including subsets of $\setR^n$) by $df$. In the Euclidean case we usually identify $df$ with the Jacobian matrix.}
\begin{equation*}
  \norm{v\circ\Psi_{\eta_n}}_{L^r(\Omega)}=\norm{v\,
    (\det d\Psi_{\eta_n}^{-1})^{1/r}}_{L^r(\Omega_{\eta_n})}\le \norm{(\det
    d\Psi_{\eta_n}^{-1})^{1/r}}_{L^s(\Omega_{\eta_n})}\,\norm{v}_{L^p(\Omega_{
      \eta_n})}.
\end{equation*}
From the convergence of $(\eta_n)$ in $C(\pa\Omega)\cap W^{1,s}(\pa\Omega)$ we infer
\begin{equation*}
 \begin{aligned}
\norm{v\circ\Psi_{\eta}}_{L^r(\Omega)}&\le \norm{(\det
d\Psi_{\eta}^{-1})^{1/r}}_{L^s(\Omega_{\eta})}\,\norm{v}_{L^p(\Omega_{
\eta})}\\ &\le c(\Omega,p,\norm{\eta}_{H^2(\pa\Omega)},\tau(\eta),r)\,\norm{v}_{L^p(\Omega_{
\eta})}.
 \end{aligned}
\end{equation*}
By the denseness of smooth functions in $L^p(\Omega_\eta)$ we deduce
the boundedness of the mapping in Lebesgue spaces.

By the chain rule and \eqref{eqn:funk} we show for $r<p$ that
\[\norm{\nabla(v\circ\Psi_{\eta_n})}_{L^r(\Omega)}\le
c(\Omega,\norm{\eta_n}_{H^2(\pa\Omega)},\tau(\eta_n),r)\,\norm{(\nabla
  v)\circ\Psi_{\eta_n}}_{L^{(r+p)/2}(\Omega)}.\] Furthermore, from the
uniform convergence of $((\nabla v)\circ\Psi_{\eta_n})$ and the
convergence of the components of the Jacobian of $\Psi_{\eta_n}$
in $L^s(\Omega)$ for $s<\infty$ we conclude the convergence of
$(\nabla(v\circ\Psi_{\eta_n}))$ to $\nabla(v\circ\Psi_{\eta})$ in
$L^r(\Omega)$. Hence, using the result shown above and the denseness
of smooth functions in $W^{1,p}(\Omega_\eta)$, see Proposition
\ref{theorem:dicht}, we obtain the boundedness of the mapping in Sobolev spaces.

The proof of the analogous claim with $\Psi_\eta$ replaced by $\Psi_\eta^{-1}$ is very similar.  
\qed

\begin{Remark}\label{bem:konv}
  Let the sequence $(\eta_n)$ converge to $\eta$ weakly in
  $H^2(\pa\Omega)$, in particular uniformly, due to $(0<\theta<1)$
\begin{align}
  \label{eq:1}
  H^2(\pa\Omega)\embedding C^{0,\theta}(\pa\Omega)\compactembedding
  C(\pa\Omega).
\end{align}
 Furthermore, we assume $\sup_n\tau(\eta_n)<\infty$.  Extending a
  function $v\in L^p(\Omega_\eta)$ by $0$ to $\setR^3$ the sequence
  $(v\circ\Psi_{\eta_n})$ converges to $v\circ\Psi_{\eta}$ in
  $L^r(\Omega)$, $r<p$. This follows, using Lemma \ref{lemma:psi} and
  approximating $v$ by $C_0^\infty(\setR^3)$ functions $\tilde v$,
  from the estimate (in $L^r(\Omega)$ norms)
\begin{equation*}
 \begin{aligned}
  \norm{v\circ\Psi_\eta-v\circ\Psi_{\eta_n}}\le \norm{\tilde v\circ\Psi_\eta-\tilde
v\circ\Psi_{\eta_n}} + \norm{(v-\tilde v)\circ\Psi_{\eta}} +
\norm{(v-\tilde v)\circ\Psi_{\eta_n}}.
 \end{aligned}
\end{equation*}
Now, let $v\in L^p(\Omega)$. Extending the functions
$v\circ\Psi^{-1}_{\eta_n}$ and $v\circ\Psi^{-1}_{\eta}$ by $0$ to
$\setR^3$, we similarly show the convergence of
$(v\circ\Psi^{-1}_{\eta_n})$ to $v\circ\Psi^{-1}_{\eta}$ in
$L^r(\setR^3)$, $r<p$.
\end{Remark}

Let us construct a trace operator for displacements $\eta\in
H^2(\pa\Omega)$. In the following we denote by $\ \cdot\
|_{\pa\Omega}$ the usual trace operator for Lipschitz domains.
\begin{Corollary}\label{lemma:spur}
  Let $1<p\le\infty$ and $\eta\in H^2(\pa\Omega)$ with
  $\norm{\eta}_{L^\infty(\pa\Omega)}<\kappa$. Then the linear mapping
  $\tr: v\mapsto (v\circ\Psi_\eta)|_{\pa\Omega}$ is well defined and
  continuous from $W^{1,p}(\Omega_\eta)$ to $W^{1-1/r,r}(\pa\Omega)$
  for all $1<r<p$. The continuity constant depends only on $\Omega$,
  $\norm{\eta}_{H^2(\pa\Omega)}$, $\tau(\eta)$, and $r$; it stays
  bounded as long as $\norm{\eta}_{H^2(\pa\Omega)}$ and $\tau(\eta)$
  stay bounded.
\end{Corollary}
\beweis The claim is a direct consequence of Lemma \ref{lemma:psi} and the continuity properties of the usual trace operator.
\qed
\medskip

From Lemma \ref{lemma:psi} and the Sobolev embeddings for regular domains we deduce Sobolev embeddings for our special domains.
\begin{Corollary}\label{lemma:sobolev}
  Let $1<p<3$ and $\eta\in H^2(\pa\Omega)$ with
  $\norm{\eta}_{L^\infty(\pa\Omega)}<\kappa$. Then
  \[
  W^{1,p}(\Omega_\eta)\compactembedding L^s(\Omega_\eta)
  \] 
  for $1\le s < p^*=3p/(3-p)$. The embedding constant depends only on
  $\Omega$, \linebreak $\norm{\eta}_{H^2(\pa\Omega)}$, $\tau(\eta)$,
  $p$, and $s$; it stays bounded as long as
  $\norm{\eta}_{H^2(\pa\Omega)}$ and $\tau(\eta)$ stay bounded.
\end{Corollary}

We denote by $H$ the mean curvature (with respect to the outer normal)
and by $G$ the Gauss curvature of $\pa\Omega$.

\begin{Proposition}\label{lemma:partInt}
  Let $1<p\le\infty$ and $\eta\in H^2(\pa\Omega)$ with
  $\norm{\eta}_{L^\infty(\pa\Omega)}<\kappa$.  Then, for $\bphi\in
  W^{1,p}(\Omega_\eta)$ with $\tr \bphi=b\,\bnu$, $b$ a scalar
  function, and $\psi\in C^1(\overline{\Omega_\eta})$ we have
  \[
  \int_{\Omega_\eta} \bphi\cdot\nabla\psi\ dx=-\int_{\Omega_\eta}
  \dv\bphi\ \psi\ dx+ \int_{\pa\Omega}b\, (1-2H\eta+G\,\eta^2)\, \tr
  \psi\ dA.
  \]
\end{Proposition}
\beweis We may assume $p<\infty$ without loss of generality. Let us
approximate $\bphi$ by functions $(\bphi_k)\subset
C^\infty_0(\setR^3)$ in $W^{1,p}(\Omega_\eta)$ and $\eta$ by
$(\eta_n)\subset C^2(\pa\Omega)$ in $H^2(\pa\Omega)$. Then,
integration by parts gives
\[
\int_{\Omega_{\eta_n}} \bphi_k\cdot\nabla\psi\
dx=-\int_{\Omega_{\eta_n}} \dv\bphi_k\ \psi\ dx +
\int_{\pa\Omega_{\eta_n}} \bphi_k\cdot\bnu_{\eta_n}\, \psi\
dA_{\eta_n}.
\]
A change of variables in the boundary integral yields
\[\int_{\pa\Omega}\tren\bphi_k\cdot(\bnu_{\eta_n}\circ \Phi_{\eta_n})\
\tren\psi\ |\det d\Phi_{\eta_n}|\ dA.\] 
In a neighborhood of every point in $\pa\Omega$ we may construct
orthonormal, tangential $C^3$ vector fields $\be_1,\be_2$ such that
$\be_1\times\be_2=\bnu$ by the Gram-Schmidt algorithm. We set $\bv_i^n:=d\Phi_{\eta_n}\be_i$. The
Jacobian determinant of $\Phi_{\eta_n}$ is equal to the area of the
parallelogram spanned by the vectors $\bv_i^n$, i.e. $|\bv_1^n\times
\bv_2^n|$. Furthermore, the normal $\bnu_{\eta_n}\circ\Phi_{\eta_n}$
is equal to $(\bv_1^n\times \bv_2^n)/|\bv_1^n\times \bv_2^n|$. In particular,
the vector fields $\bv^n:=\bv_1^n\times \bv_2^n$ are independent of
the choice of the fields $\be_i$. Hence, this local definition, in
fact, yields globally, i.e. on $\pa\Omega$, defined vector fields
$\bv^n$. Thus, the boundary integral may be written in the form
\[\int_{\pa\Omega}\tren\bphi_k\cdot\bv^n\,
\tren\psi\ dA.\] For $q\in\pa\Omega$ and $c$ a curve in $\pa\Omega$
with $c(0)=q$ and $\frac{d}{dt}\big|_{t=0}\ c(t) = \be_i(q)$ we have
\begin{equation}\label{eqn:v}
 \begin{aligned}
\bv_i^n(q)&=\frac{d}{dt}\Big|_{t=0}\Phi_{\eta_n}(c(t))=\frac{d}{dt}\big|_{t=0}
(c(t)+\eta_n(c(t))\,\bnu(c(t)))\\
&=\be_i(q)+d\eta_n\be_i(q)\, \bnu(q)+\eta_n(q)\,\frac{d}{dt}\Big|_{t=0}
\bnu(c(t))\\
&=\be_i(q)+d\eta_n\be_i(q)\, \bnu(q)-\eta_n(q)\, h_i^j(q)\, \be_j(q).
 \end{aligned}
\end{equation}
Here, $h^j_i$ denote the components of the Weingarten map with respect
to the orthonormal basis $\be_1,\be_2$. Since $(\eta_n)$ converges in
$H^2(\pa\Omega)$, the sequence $(\bv_i^n)$ converges in $L^r$ on its
domain of definition for all $r<\infty$. Hence, $(\bv^n)$ converges to
a field $\bv$ in $L^r(\pa\Omega)$ for all $r<\infty$. If we let first
$n$ and then $k$ tend to infinity, we obtain
\[\int_{\Omega_\eta} \bphi\cdot\nabla\psi\ dx=-\int_{\Omega_\eta} \dv\bphi\
\psi\ dx + \int_{\pa\Omega} b\,\bnu\cdot\bv\, \tr\psi\ dA.\] 
Now we merely have to compute $\bnu\cdot\bv$. From \eqref{eqn:v} we
see that
\[\bnu\cdot(\bv_1^n\times\bv_2^n)=1-(h_1^1+h^2_2)\,\eta_n+(h_1^1h_2^2-
h_1^2h_2^1)\,\eta_n^2=1-2H\eta_n+G\,\eta_n^2.\] Taking the limit
completes the proof.
\qed\medskip

We shall make repeated use of the field $\bv\in L^r(\pa\Omega)$,
$r<\infty$, constructed in the proof above. We denote it by $\bv_\eta$ in
order to emphasize its dependence on $\eta$.
\begin{Remark}\label{bem:trafo}
  The proof of Proposition \ref{lemma:partInt} shows that for $\psi\in
  L^1(\pa\Omega)$ and sufficiently smooth $\eta$ we have
  \[
  \int_{\pa\Omega_{\eta}}
  \big((\psi\,\bnu)\circ\Phi_{\eta}^{-1}\big)\cdot\bnu_\eta\
  dA_\eta=\int_{\pa\Omega} \psi\,\bnu\cdot\bv_\eta\
  dA=\int_{\pa\Omega} \psi\,(1-2H\eta+G\,\eta^2)\ dA.
  \]
\end{Remark}
\smallskip

\begin{Remark}\label{bem:groessernull}
  Let us  show that the function 
  \begin{align}
    \label{eq:2}
    \gamma(\eta):=1-2H\eta+G\,\eta^2
  \end{align}
  is positive as long as $|\eta|<\kappa$. For $G\not=0$ the roots of
  this polynomial are
  \[
  \frac{H}{G}\pm\frac{\sqrt{H^2-G}}{|G|}=\frac{1}{2}
  \bigg(\frac{h_1+h_2}{h_1h_2}\pm\frac{|h_1-h_2|}{ |h_1h_2|}\bigg),
  \]
  where $h_1,h_2$ are the principal curvatures. Thus, their
  absolute values are
  \[
  \frac{|h_1+h_2\pm(h_1-h_2)|}{2|h_1h_2|}=|h_1|^{-1},|h_2|^{-1}.
  \]
  Hence, for $G\not=0$, the claim is proved provided that we can show that
  \[\kappa\le\min(|h_1(q)|^{-1},|h_2(q)|^{-1})\] for all
  $q\in\pa\Omega$. But, with this estimate at hand, the proof is trivial for $G=0$. In order to
  show the estimate we consider the diffeomorphism
  \[
  \Lambda_\alpha:=\Lambda(\,\cdot\, ,\alpha):\pa\Omega\rightarrow
  \Lambda_\alpha(\pa\Omega)\subset\setR^3
  \]
  with $\alpha\in(-\kappa,\kappa)$. For an eigenbasis $\be_1,\be_2$ of
  the Weingarten map in $q\in \pa\Omega$ we have
  $d\Lambda_\alpha\,\be_i=(1-\alpha\, h_i(q))\,\be_i$. Thus, we deduce
  the estimate on $\kappa$ from the fact that the differential of
  $\Lambda_\alpha$ gets singular for
  $|\alpha|\nearrow\min(|h_1(q)|^{-1},|h_2(q)|^{-1})$.
\end{Remark}

Consider the space
\[E^p(\Omega_\eta):=\{\bphi\in
L^p(\Omega_\eta)\ |\ \dv \bphi\in
L^p(\Omega_\eta)\}\]
for $1\le p\le\infty$, endowed with the canonical norm. 
\begin{Proposition}\label{lemma:nspur}
  Let $1<p<\infty$ and $\eta\in H^2(\pa\Omega)$ with
  $\norm{\eta}_{L^\infty(\pa\Omega)}<\kappa$.  Then there exists a
  continuous, linear operator
\[\trnormal: E^p(\Omega_\eta)\rightarrow (W^{1,p'}(\pa\Omega))'\]
such that for $\bphi\in E^p(\Omega_\eta)$ and $\psi\in
C^1(\overline{\Omega_\eta})$
\[\int_{\Omega_\eta} \bphi\cdot\nabla\psi\ dx=-\int_{\Omega_\eta} \dv\bphi\
\psi\ dx + \langle
\trnormal\bphi,\tr\psi\rangle_{W^{1,p'}(\pa\Omega)}.\] The continuity
constant depends only on $\Omega$, $\tau(\eta)$, and $p$; it stays
bounded as long as $\tau(\eta)$ stays bounded.
\end{Proposition}
\beweis It suffices to consider $\bphi\in C^1(\overline{\Omega_\eta})$
since these functions are dense in $E^p(\Omega_\eta)$; see Proposition
\ref{theorem:dicht}. Analogously to the proof of Proposition
\ref{lemma:partInt} we may show that
\begin{equation}\label{eqn:approxnspur}
 \begin{aligned}
\int_{\Omega_\eta} \bphi\cdot\nabla\psi\ dx+\int_{\Omega_\eta} \dv\bphi\
 \psi\ dx = \int_{\pa\Omega} \tr\bphi\cdot\bv_\eta\,\tr\psi\ dA
 \end{aligned}
\end{equation}
for $\psi\in C^1(\overline\Omega_\eta)$. The left-hand side is
dominated by $2\norm{\bphi}_{E^p
  (\Omega_\eta)}\norm{\psi}_{W^{1,p'}(\Omega_\eta)}$.  Obviously, the
identity holds for all $\psi\in W^{1,p'}(\Omega_\eta)$.

Let us now show that the mapping $b\mapsto b\circ q$ defines a
bounded, linear extension operator from $W^{1,p'}(\pa\Omega)$ to
$W^{1,p'}(S_\alpha\cap\Omega_\eta)$ for fixed, but arbitrary $\alpha$
with $\norm{\eta}_{L^\infty(\pa\Omega)}<\alpha<\kappa$. By a change of
variables we obtain
\[
\int_{S_\alpha\cap\Omega_\eta}|b\circ q|^{p'}\
dx=\int_{\pa\Omega}|b|^{p'}\int_{-\alpha}^\eta|\det d\Lambda|\ dsdA.
\]
Hence, the $L^{p'}$ norm of $b\circ q$ is bounded by a constant times
the $L^{p'}$ norm of $b$, the constant depending on $\Omega$ and
$\alpha$. Approximating $b$ by sufficiently smooth functions we can
show the chain rule
\[
\nabla(b\circ q)=(\nabla b)\circ q\ dq.
\]
Thus, the $L^{p'}$ norm of $\nabla(b\circ q)$ is bounded by a constant
times the $L^{p'}$ norm of $\nabla b$, the constant depending on
$\Omega$, $\alpha$, and $p'$. By another approximation argument we
obtain the identity $\tr (b\circ q)=b$.

Multiplying the constructed extension by the cut-off function
$x\mapsto\beta(|s(x)|)$, where $\beta\in C^\infty(\setR)$, $\beta=1$
in a neighborhood of the interval $[0,\norm{\eta}_{L^\infty(\pa\Omega)}]$ and
$\beta=0$ in a neighborhood of $\alpha$, we finally obtain
a bounded, linear extension operator from $W^{1,p'}(\pa\Omega)$ to
$W^{1,p'}(\Omega_\eta)$. Thus, the right hand side of
\eqref{eqn:approxnspur} defines the trace operator we were looking
for.  \qed\smallskip

\begin{Proposition}\label{lemma:divdicht}
  Let $\eta\in H^2(\pa\Omega)$ with
  $\norm{\eta}_{L^\infty(\pa\Omega)}<\kappa$. Then the subspace of
  functions whose supports are contained in $\Omega_\eta$ is dense in the
  (canonically normed) space
  \[
  H(\Omega_\eta):=\{\bphi\in L^2(\Omega_\eta)\ |\ \dv\bphi=0,\
  \trnormal\bphi =0 \}.
  \] 
\end{Proposition}
\beweis Let $\bphi\in H(\Omega_\eta)$ with
\begin{equation*}
 \begin{aligned}
  \int_{\Omega_\eta} \bphi\cdot\bpsi\ dx =0
 \end{aligned}
\end{equation*}
for all $\bpsi\in L^2(\Omega_\eta)$ with $\dv\bpsi=0$ and $\supp
\bpsi\subset\Omega_\eta$. By DeRham's theorem there exists
a function $p\in L^2_{\loc}(\Omega_\eta)$ with $\bphi=\nabla p$.
Furthermore, from Proposition \ref{lemma:nspur} we deduce that
\begin{equation*}
 \begin{aligned}
  \int_{\Omega_\eta} \bphi\cdot\nabla \psi\ dx =0 
 \end{aligned}
\end{equation*}
for all $\psi\in C^1(\overline{\Omega_\eta})$. The proof of
Proposition \ref{theorem:dicht} shows that $\nabla p$ may be approximated
in $L^2(\Omega_\eta)$ by functions $\nabla\psi$, $\psi\in
C^1(\overline{\Omega_\eta})$. Thus, we have
\begin{equation*}
  \begin{aligned}
    \int_{\Omega_\eta} |\bphi|^2\ dx = \int_{\Omega_\eta}
    \bphi\cdot\nabla p\ dx = 0 .
  \end{aligned}
\end{equation*} 
\qed
\medskip

\pagebreak
Now, we shall construct an operator extending suitable boundary values
$b\,\bnu$ to divergence-free vector fields.
\begin{Proposition}\label{lemma:FortVonRand}
  Let $1<p<\infty$, $\eta\in H^2(\pa\Omega)$ with
  $\norm{\eta}_{L^\infty(\pa\Omega)}<\kappa$, and $\alpha$ such that
  $\norm{\eta}_{L^\infty(\pa\Omega)}<\alpha<\kappa$. Then there exists a
  bounded, linear extension operator
  \[
  \F_\eta: \Big\{b\in W^{1,p}(\pa\Omega)\ \big|\
  \int_{\pa\Omega}b\,\gamma(\eta)\ 
  dA=0\Big\}\rightarrow W^{1,p}_{\dv}(B_\alpha);
  \]
  in particular $\tr\F_\eta b=b\,\bnu.$ The continuity constant
  depends only on $\Omega$, \linebreak $\norm{\eta}_{H^2(\pa\Omega)}$,
  $\alpha$, and $p$; it stays bounded as long as
  $\norm{\eta}_{H^2(\pa\Omega)}$ and $\tau(\alpha)$ stay bounded.
\end{Proposition}
\beweis Let $b\in W^{1,p}(\pa\Omega)$ with
\begin{equation}\label{eqn:integralnull}
 \begin{aligned}
\int_{\pa\Omega}b\,\gamma(\eta)\ dA=0.  
 \end{aligned}
\end{equation}
For $x\in S_\alpha$ we define
\begin{equation}\label{eqn:fort}
 \begin{aligned}
(\F_\eta b)(x):=\exp\Big(\int_{\eta(q(x))}^{s(x)}
\beta(q(x)+\tau\,\bnu(q(x)))\ d\tau\Big)\,
(b\,\bnu)(q(x)),  
 \end{aligned}
\end{equation}
where $\beta:=-\dv(\bnu\circ q)$. Obviously $\beta$ is a $C^2$
function. By the chain rule we see that $\F_\eta b$ is weakly
differentiable and
\begin{align}
  &\pa_i\, \F_\eta b \label{eqn:fortnabla}
  \\
  &= \Big[\pa_i((b\,\bnu)\circ q) + (b\,\bnu)\circ q\
  \Big(\int_{\eta\circ q}^s \pa_i(\beta(q+\tau\,\bnu\circ q))\ d\tau +
  \beta(q+s\,\bnu\circ q)\ \pa_i s \notag \\
  &\hspace{40mm}\quad- \beta(q+(\eta\,\bnu)\circ q)\ \pa_i (\eta\circ
  q)\Big)\Big]\,e^{\int_{\eta\circ q}^s \beta(q+\tau\bnu\circ q)
    d\tau}.\notag 
\end{align}
Since $\pa_i((b\,\bnu)\circ q)\in L^p(S_\alpha)$, $(b\,\bnu)\circ q\in
L^{r}(S_\alpha)$ for some $r>p$, $\pa_i (\eta\circ q)\in
L^r(S_\alpha)$ for all $r<\infty$, and all other terms are bounded, we
see that $\F_\eta b$ is bounded in $W^{1,p}(S_\alpha)$. Furthermore,
we have\footnote{For a scalar field $f$ and a vector field $\bX$ we have $\dv (f\,\bX)=f\dv\bX + \nabla f\cdot\bX$.}
\begin{equation*}
  \begin{aligned}
    \dv\F_\eta b&=\dv (\bnu\circ q)\, e^{\int_{\eta\circ q}^s
      \beta(q+\tau\bnu\circ q) d\tau}\,b\circ q + (\bnu\circ q)\cdot
    \nabla (e^{\int_{\eta\circ q}^s
      \beta(q+\tau\bnu\circ q) d\tau}\,b\circ q)\\
    &= (-\beta\, e^{\int_{\eta\circ q}^s\beta(q+\tau\bnu\circ q)
      d\tau} + \pa_s e^{\int_{\eta\circ q}^s \beta(q+\tau\bnu\circ q)
      d\tau})\, b\circ q =0.
  \end{aligned}
\end{equation*}
For the second equality we used the defintion of $\beta$ and the fact
that for $x\in S_\alpha$
\[
dq\,\bnu(q(x))=\frac{d}{dt}\Big|_{t=0}q(x+t\,\bnu(q(x)))= 
\frac{d}{dt}\Big|_{t=0}q(x)=0
\]
and 
\[
ds\,\bnu(q(x))=\frac{d}{dt}\Big|_{t=0}s(x+t\,\bnu(q(x)))=
\frac{d}{dt}\Big|_{t=0}s(x)+t=1.
\]

Approximating $\eta$ and $b$ by $C^2$ functions $(\eta_n)$ in
$H^2(\pa\Omega)$ and $C^1$ functions $(b_n)$ in $W^{1,p}(\pa\Omega)$,
resp., the $C^1(\overline{S_\alpha})$ functions
\[
\bphi_n:=\exp\Big(\int_{\eta_n\circ q}^{s} \beta(q+\tau\,\bnu\circ q)\
d\tau\Big)\, (b_n\,\bnu)\circ q
\]
converge to $\F_\eta b$ in $W^{1,p}(S_\alpha)$. Furthermore, the traces
\begin{equation}\label{eqn:approxspur}
  \begin{aligned}
    \tr\bphi_n=\exp\Big(\int_{\eta_n}^{\eta} \beta(q+\tau\,\bnu\circ
    q)\ d\tau\Big)\,b_n\,\bnu
  \end{aligned}
\end{equation}
converge to $b\,\bnu$. Thus, $\tr\F_\eta b=b\,\bnu$.

Let us now extend $\F_\eta b$ to $\Omega\setminus\overline{S_\alpha}$.
By Proposition \ref{lemma:partInt} we have
\begin{equation*}
  \begin{aligned}
    \int_{\pa(\Omega\setminus\overline{S_\alpha})}(\F_\eta
    b)\cdot\bnu\ dA=-\int_{\pa\Omega}b\, \gamma(\eta)\
    dA+\int_{S_\alpha\cap\Omega_\eta} \dv\F_\eta b\ dx=0.
  \end{aligned}
\end{equation*}
Here, $\bnu$ denotes the inner unit normal of
$\pa(\Omega\setminus\overline{S_\alpha})$. Due to this identity and
since $\pa(\Omega\setminus\overline{S_\alpha})$ is a $C^3$ boundary we
may solve the Stokes system in $\Omega\setminus \overline{S_\alpha}$
with boundary values $(\F_\eta
b)|_{\pa(\Omega\setminus\overline{S_\alpha})}$; cf.~Theorem
\ref{theorem:stokes}. This defines the extension we were looking
for.\qed

\begin{Remark}\label{bem:nspur}
  If the identity \eqref{eqn:integralnull} holds for a function
  $b\in L^p(\pa\Omega)$, $1<p<\infty$, we may construct a
  divergence-free extension $\F_\eta b\in L^p(B_\alpha)$ exactly like
  in the proof above. This way we obtain a bounded, linear operator
  \[
  \F_\eta: \Big\{b\in L^{p}(\pa\Omega)\ \big|\ \int_{\pa\Omega}b\,
  \gamma(\eta)\ dA=0\Big\}\rightarrow \{\bphi\in L^p(B_\alpha)\ |\
  \dv\bphi=0\},
  \]
  whose continuity constant depends on the data as in Proposition
  \ref{lemma:FortVonRand}.

  When constructing the extension we need to apply the solution
  operator of the Stokes-system to the formal trace
  \[
  \bphi|_{\pa(\Omega\setminus\overline{S_\alpha})}=\exp\Big(\int_{\eta\circ
    q}^{-\alpha} \beta(q+\tau\,\bnu\circ q)\ d\tau\Big)\,
  (b\,\bnu)\circ q.
  \]
  Obviously, this trace is bounded in
  $L^p(\pa(\Omega\setminus\overline{S_\alpha}))$ and we have
  \begin{equation}\label{eqn:formalmitt}
    \begin{aligned}
      \int_{\pa(\Omega\setminus\overline{S_\alpha})}\bphi\cdot\bnu\
      dA=0.
    \end{aligned}
  \end{equation}
  In the proof that the resulting vector field in $B_\alpha$ is
  divergence-free we need to use the fact that for $\psi\in
  C_0^\infty(\Omega)$ we have
  \begin{equation}\label{eqn:horst}
    \begin{aligned}
      \int_{S_\alpha}\bphi\cdot\nabla\psi\
      dx=\int_{\pa(\Omega\setminus\overline{S_\alpha})}\bphi\cdot\bnu\,
      \psi\
      dA=-\int_{\Omega\setminus\overline{S_\alpha}}\bphi\cdot\nabla\psi\
      dx.
    \end{aligned}
  \end{equation}
  The identities \eqref{eqn:formalmitt} and \eqref{eqn:horst} can be
  easily proven by approximating $b$ by sufficiently smooth
  functions.  Again by an approximation argument, we may deduce from
  \eqref{eqn:approxspur} and \eqref{eqn:approxnspur} the identity
  \begin{equation}\label{eqn:nspur}
    \begin{aligned}
      \trnormal \bphi=b\, \gamma(\eta),
    \end{aligned}
  \end{equation}
  justifying the term \emph{extension operator}.
\end{Remark}
\smallskip

\begin{Proposition}\label{lemma:nullfort}
  Let $1<p<3$ and $\eta\in H^2(\pa\Omega)$ with
  $\norm{\eta}_{L^\infty(\pa\Omega)}<\kappa$. Then extension by $0$
  defines a bounded, linear operator from $W^{1,p}(\Omega_\eta)$ to
  $W^{1/4,p}(\setR^3)$. The continuity constant depends only on
  $\Omega$, $p$, $\norm{\eta}_{H^2(\pa\Omega)}$, and $\tau(\eta)$; it
  stays bounded as long as $\norm{\eta}_{H^2(\pa\Omega)}$ and
  $\tau(\eta)$ stay bounded.
\end{Proposition}
\beweis To begin with, let us show that extension by $0$ defines a
bounded, linear operator from $W^{1,r}(\Omega)$ to $W^{s,r}(\setR^3)$
with $1\le r<3$ and $s<1/3$. To this end, it suffices to estimate the
integral
\begin{equation*}
  \begin{aligned}
    &\int_{\setR^3}\int_{\setR^3} \frac{|v(x)-v(y)|^r}{|x-y|^{3+sr}}\
    dydx \\
    &\quad = \int_{\Omega}\int_{\Omega}
    \frac{|v(x)-v(y)|^r}{|x-y|^{3+sr}}\ dydx +2
    \int_{\Omega}\int_{\setR^3\setminus\Omega}
    \frac{|v(x)|^r}{|x-y|^{3+sr}}\ dydx\\ 
    &\quad=\int_{\Omega}\int_{\Omega} \frac{|v(x)-v(y)|^r}{|x-y|^{3+sr}}\
    dydx+ 2\int_{\Omega} |v(x)|^r \int_{\setR^3\setminus\Omega}
    \frac{1}{|x-y|^{3+sr}}\ dydx
 \end{aligned}
\end{equation*}
for $v\in W^{1,r}(\Omega)$. While the first term on the right-hand
side is dominated by $c\,\norm{v}_{W^{1,r}(\Omega)}$, we can estimate
the interior integral of the second term by
\[
\int_{|z|>d(x)}\frac{1}{|z|^{3+sr}}\ dz=\frac{c(s,r)}{d(x)^{sr}},
\]
where $d(x)$ denotes the distance from $x$ to $\pa\Omega$. Thus,
applying H{\"o}lder's inequality with the exponents $r^*/r$ and $(r^*/r)'$
to the second summand we see that it is dominated by
$c(s,r)\,\norm{v}_{r^*}^r\,\norm{d(\cdot)^{-s}}_{L^3(\Omega)}^r$. The
identity
\[
\int_{S_{\kappa/2}\cap\Omega}|d(x)|^{-3s}\
dx=\int_{\pa\Omega}\int_{-\kappa/2}^0|\det
d\Lambda|\, \alpha^{-3s} \ d\alpha dA,
\]
a consequence of a change of variables, and the inequality $3s<1$ show
that the factor $\norm{d(\cdot)^{-s}}_{L^3(\Omega)}^r$ is finite. This
proves the claim.

Concatenating the mapping from Lemma \ref{lemma:psi} with the
extension by $0$ thus yields a bounded, linear operator
\[
W^{1,p}(\Omega_\eta)\rightarrow W^{1,r}(\Omega)\rightarrow
W^{s,r}(\setR^3)
\]
for $r<p$ and $s<1/3$. Now, let $\delta\in C^4(\pa\Omega)$ with
$\eta<\delta<\kappa$. If we extend the mapping
$\Psi_\delta$ to $\overline{B_\alpha}$, $\alpha>0$ sufficiently small,
by using the definition \eqref{eqn:defpsi} for $x\in S_\kappa\cap
\overline{B_\alpha}$, we obtain a $C^3$ diffeomorphism
\[
\widetilde\Psi_\delta: \overline{B_\alpha}\rightarrow
\overline{\Omega_{\delta+\alpha}}.
\] 
Since the fractional Sobolev spaces are interpolation spaces the
linear mapping $v\mapsto v\circ\widetilde\Psi_\delta^{-1}$ is bounded
from $W^{s,r}(B_\alpha)$ to $W^{s,r}(\Omega_{\delta+\alpha})$. By the
same reason, the restriction of functions is bounded from
$W^{s,r}(\setR^3)$ to $W^{s,r}(B_\alpha)$. Furthermore, for suitable
$s<1/3$ and $r<p$ we have
\[
W^{s,r}(\Omega_{\delta+\alpha})\embedding
W^{1/4,p}(\Omega_{\delta+\alpha}).
\]
Concatenating these mappings we see that the extension by $0$ defines
a bounded, linear operator from $W^{1,p}(\Omega_\eta)$ to
$W^{1/4,p}(\Omega_{\delta+\alpha})$. Since $\Omega_\eta$ has a
positive distance from $\setR^3\setminus\Omega_{\delta+\alpha}$ we may
conclude the proof by the estimate
\begin{equation*}
  \begin{aligned}
    |v|_{1/4,p;\setR^3}^p&=|v|_{1/4,p;\Omega_{\delta+\alpha}}^p
    +\int_{\Omega_\eta}\int_{ \setR^3\setminus\Omega_{\delta+\alpha}}
    \frac{|v(x)|^p}{|x-y|^{3+p/4}}\ dydx\\
    &\le |v|_{1/4,p;\Omega_{\delta+\alpha}}^p +
    c\,\norm{v}_{L^p(\Omega_\eta)}^p.
  \end{aligned}
\end{equation*}
\qed
\medskip

The usual Bochner spaces are not the right objects to deal with
functions defined on time-dependent domains. For this reason we now
define an (obvious) substitute for these spaces. For $I:=(0,T)$,
$T>0$, and $\eta\in C(\bar I\times\pa\Omega)$ with
$\norm{\eta}_{L^\infty(I\times\pa\Omega)}<\kappa$ we set
$\Omega_\eta^I:=\bigcup_{t\in I}\, \{t\}\times \Omega_{\eta(t)}$.
Note that $\Omega_\eta^I$ is a domain in $\setR^4$. For $1\le
p,r\le\infty$ we set
\begin{equation*}
 \begin{aligned}
   L^p(I,L^r(\oet))&:=\{v\in L^1(\Omega_\eta^I)\ |\ v(t,\cdot)\in
   L^r(\Omega_{\eta(t)})\text{ for almost all $t$ and}\\
   &\hspace{.8cm}\norm{v(t,\cdot)}_{L^r(\Omega_{\eta(t)})}\in
   L^p(I)\},\\
   L^p(I,W^{1,r}(\Omega_{\eta(t)}))&:=\{v\in L^p(I,L^r(\oet))\ |\
   \nabla v\in
   L^p(I,L^r(\oet))\},\\
   L^p(I,W_{\dv}^{1,r}(\Omega_{\eta(t)}))&:=\{\bv\in
   L^p(I,W^{1,r}(\Omega_{\eta(t)}))\ |\
   \dv \bv=0\},\\
   W^{1,p}(I,W^{1,r}(\Omega_{\eta(t)}))&:=\{\bv\in
   L^p(I,W^{1,r}(\Omega_{\eta(t)}))\ |\ \pa_t\bv\in
   L^p(I,W^{1,r}(\Omega_{\eta(t)}))\}.
 \end{aligned}
\end{equation*}
Here $\nabla$ and $\dv$ are acting with respect to the space variables. Furthermore, we set
\[
\Psi_\eta: \bar I\times \overline\Omega\rightarrow
\overline{\Omega_\eta^I},\ (t,x)\mapsto (t,\Psi_{\eta(t)}(x))
\]
and
\[
\Phi_\eta: \bar I\times \pa\Omega\rightarrow\bigcup_{t\in \bar I}\,
\{t\}\times \pa\Omega_{\eta(t)},\ (t,x)\mapsto (t,\Phi_{\eta(t)}(x)).
\] 
If $\eta\in L^\infty(I,H^2(\pa\Omega))$ we obtain
``instationary'' versions of the claims made so far by applying these
at (almost) every $t\in I$. For instance, from Corollary \ref{lemma:sobolev} we
deduce that
\[
L^2(I,H^{1}(\Omega_{\eta(t)}))\embedding L^2(I,L^s(\Omega_{\eta(t)}))
\]
for $1\le s<2^*$. Note that the construction given above does not
provide a substitute for Bochner spaces of dual-space valued
functions.

Note that for all $1/2<\theta<1$ we have
\begin{align}
  \begin{aligned}
    &W^{1,\infty}(I,L^2(\pa\Omega))\cap
    L^\infty(I,H^2(\pa\Omega))\label{eqn:hoeldereinb}
    \\
    &\hspace{2cm} \embedding C^{0,1-\theta}(\bar
    I,H^{2\theta}(\pa\Omega))\embedding C^{0,1-\theta}(\bar I,
    C^{0,2\theta -1}(\pa\Omega)).
  \end{aligned}
\end{align}
While the second embedding follows from standard results,
the first one is a consequence of the elementary estimate
\begin{equation*}
 \begin{aligned}
   \norm{u(t)-u(s)}_{(L^2(\pa\Omega),H^2(\pa\Omega))_{\theta,2}}&\le
   \norm{u(t)-u(s)}
   ^{\theta}_{H^2(\pa\Omega)}\,\norm{u(t)-u(s)}_{L^2(\pa\Omega)}^{1-\theta}\\
   &\le c \norm{u}_{L^\infty(I,H^2(\pa\Omega))}^{\theta}\,
   \norm{u}_{W^{1,\infty}(I,L^2(\pa\Omega))}^{1-\theta}
   \,|t-s|^{1-\theta}.
 \end{aligned}
\end{equation*}

\begin{Proposition}\label{lemma:FortVonRandZeit}
  Let $\eta\in W^{1,\infty}(I,L^2(\pa\Omega))\cap
  L^\infty(I,H^2(\pa\Omega))$ be given with \linebreak
  ${\norm{\eta}_{L^\infty(I\times\pa\Omega)}<\kappa}$ and $\alpha$ a
  real number such that
  $\norm{\eta}_{L^\infty(I\times\pa\Omega)}<\alpha<\kappa$. The
  application of the extension operator from Proposition
  \ref{lemma:FortVonRand} at (almost) all times defines a bounded,
  linear extension operator $\F_\eta$ from
  \[
  \Big\{b\in H^{1}(I,L^2(\pa\Omega))\cap L^2(I,H^2(\pa\Omega))\ |\
  \int_{\pa\Omega}b(t,\cdot)\, \gamma(\eta(t,\cdot))\ dA=0\text{ for
    all }t\in I\Big\}
  \] 
  to
  \[
  \{\bphi\in H^{1}(I,L^2(B_\alpha))\cap C(\bar I,H^1(B_\alpha))\ |\
  \dv\bphi=0\}.
  \]
  The continuity constant depends only on $\Omega$, $\norm{\eta}_{
    W^{1,\infty}(I,L^2(\pa\Omega))\cap L^\infty(I,H^2(\Omega))}$, and
  $\alpha$; it stays bounded as long as $\norm{\eta}_{
    W^{1,\infty}(I,L^2(\pa\Omega))\cap L^\infty(I,H^2(\Omega))}$ and
  $\tau(\alpha)$ stay bounded.
\end{Proposition}
\beweis In view of
\begin{equation}\label{eqn:einbschnurz}
 \begin{aligned}
   H^1(I,L^2(\pa\Omega))\cap L^2(I,H^2_0(\pa\Omega))\embedding C(\bar
   I,H^1(\pa\Omega)),
 \end{aligned}
\end{equation}
the continuity to $\{\bphi\in L^\infty(I,H^1(B_\alpha))\ |\
\dv\bphi=0\}$ is a direct consequence of Proposition
\ref{lemma:FortVonRand}. Furthermore, from \eqref{eqn:einbschnurz},
\eqref{eqn:hoeldereinb}, and standard embedding results we deduce that
\begin{equation*}
 \begin{aligned}
   H^1(I,L^2(\pa\Omega))\cap L^2(I,H^2_0(\pa\Omega))
   \embedding C(\bar I,L^4(\pa\Omega)),\\
   W^{1,\infty}(I,L^2(\pa\Omega))\cap L^\infty(I,H^2(\pa\Omega))
   \embedding C(\bar I,W^{1,4}(\pa\Omega)).
 \end{aligned}
\end{equation*}
Thus, in view of \eqref{eqn:fort}, \eqref{eqn:fortnabla}, and the
mapping properties of the solution operator of the Stokes system we
deduce that $\F_\eta b\in C(\bar I,H^1(B_\alpha))$. Let us now
estimate $\F_\eta b$ in $H^1(I,L^2(B_\alpha))$. In $I\times S_\alpha$
we have
\begin{equation}\label{eqn:klopps}
 \begin{aligned}
   &\pa_t\, (\F_\eta b)(t,\cdot)\\
   &\hspace{1cm}= \Big[(\pa_tb)(t,q)\ \bnu\circ q -
   \beta(q+\eta(t,q)\, \bnu\circ q)\
   \pa_t\eta(t,q) \  b(t,q)\  \bnu\circ q\Big]\\
   &\hspace{1.6cm} e^{\int_{\eta\circ q}^s \beta(q+\tau\bnu\circ q) d\tau}.
\end{aligned}
\end{equation}
Since $\pa_tb\in L^2(I,L^2(\pa\Omega))$ the first summand can be
estimated in $L^2(I,L^2(S_\alpha))$. The estimate of the second
summand in $L^2(I,L^2(S_\alpha))$ follows from the fact that
$\pa_t\eta\in L^\infty(I,L^2(\pa\Omega))$ and $b\in
L^2(I,L^\infty(\pa\Omega))$. The time-derivative of the trace of
$\F_\eta b$ on $I\times\pa(\Omega\setminus\overline{S_\alpha})$ equals
\[
\Big[(\pa_tb)(t,q)\ \bnu\circ q - \beta(q+\eta(t,q)\,\bnu\circ q)\
\pa_t\eta(t,q) \ b(t,q)\ \bnu\circ q\Big]\,e^{\int_{\eta\circ
    q}^{-\alpha} \beta(q+\tau\bnu\circ q) d\tau},
\] 
and, thus, may be estimated in
$L^2(I,L^2(\pa(\Omega\setminus\overline{S_\alpha})))$.  Using the
mapping properties of the solution operator of the Stokes system we
conclude the proof. 
\qed

\begin{Remark}
  The proof above shows that, under the assumptions of Proposition
  \ref{lemma:FortVonRandZeit}, the application of the extension
  operator from Remark \ref{bem:nspur} at (almost) all times defines a
  bounded linear extension operator $\F_\eta$ from
  \[
  \Big\{b\in C(\bar I,L^2(\pa\Omega))\ |\ \int_{\pa\Omega}b(t,\cdot)\,
  \gamma(\eta(t,\cdot))\ dA=0\text{ for almost all }t\in I\Big\}
  \] 
  to
  \[
  \{\bphi\in C(\bar I,L^2(B_\alpha))\ |\ \dv\bphi=0\}.
  \]
  The continuity constant depends on the data as in Proposition
  \ref{lemma:FortVonRandZeit}.
\end{Remark}

\begin{Remark}\label{bem:tdelta2}
For $\eta\in C^2(I\times\pa\Omega)$ with $\norm{\eta}_{L^\infty(I\times\pa\Omega)}<\kappa$ an application of $\T_{\eta(t)}$
for each $t\in I$ defines isomorphisms between appropriate function spaces on $I\times\Omega$ respectively $\Omega_\eta^I$ as
long as the order of differentiability is not larger than $1$.
\end{Remark}

\section{Main result}
\label{sec:2}


We define 
\begin{equation*}
 Y^I:=W^{1,\infty}(I,L^2(M))\cap L^\infty(I,H_0^2(M)),
\end{equation*}
and for $\eta\in Y^I$ with $\norm{\eta}_{L^\infty(I\times M)}<\kappa$
we set
\begin{equation*}
X_\eta^I:=L^\infty(I,L^2(\Omega_{\eta(t)}))\cap L^2(I,H^1_{\dv}(\oet)).
\end{equation*}
Here and throughout the rest of the paper we tacitly extend functions
defined in $M$ by $0$ to $\pa\Omega$. Furthermore, we define the space of
test function $T_\eta^I$ to consists of  all couples
\[(b,\bphi)\in \big(H^1(I,L^2(M))\cap L^2(I,H^2_0(M))\big)\times
\big(H^1(\Omega_{\eta}^I)\cap
L^\infty(I,L^4(\oet))\big)\]
such that $b(T,\cdot)=0$, $\bphi(T,\cdot)=0$\footnote{We will see soon that it makes sense to evaluate $\bphi$ at a fixed
point in time.}, $\dv\bphi=0$
and $\bphi-\F_\eta b\in H_0$. Here $H_0$ denotes the closure in $H^1(\Omega_{\eta}^I)\cap
L^\infty(I,L^4(\oet))$ of the divergence-free elements of this space that vanish at $t=T$ and whose supports are contained in
$\Omega_{\eta}^{\bar I}$. From the last requirement we can infer that $\tr\bphi=\tr \F_\eta b=b\,\bnu$. The converse holds at least
for sufficiently regular $\eta$ and under the premise that we take the smaller space
\[H^1(I,H^2_0(M))\times H^1(I,H^1(\Omega_{\eta(t)}))\]
as a basis for the test functions.\footnote{Using the mapping $\T_{\eta}$ from the Remarks \ref{bem:tdelta} and \ref{bem:tdelta2} we see that in the case of
a sufficiently regular $\eta$ it suffices to consider the analogous situation in a spacetime cylinder $I\times\Omega$. There we
can construct suitable approximations by a standard argument, using the solution operator of the divergence equation, see for
example III.4.1 in \cite{b18}. For each $t\in\bar I$ these approximations lie in $C_0^\infty(\Omega)$ and converge in
$H^1(\Omega)$. Since this construction commutes with the time derivative we obtain the convergence in
$H^1(I,H^1(\Omega_{\eta(t)}))$.}
From $\tr \bphi = b\,\bnu$ we infer that $\bphi$ vanishes on $\Gamma$.

We say that the data $(\ff,g,\bu_0,\eta_0,\eta_1)$ is
\emph{admissible} if $\ff\in L_{\text{loc}}^{2}([0,\infty)\times
\setR^3)$,
$g\in L_{\text{loc}}^2([0,\infty)\times M)$, $\eta_0\in H^2_0(M)$ with
$\norm{\eta_0}_{L^\infty(M)}<\kappa$, $\eta_1\in L^2(M)$, and $\bu_0\in
L^2(\Omega_{\eta_0})$ with $\dv \bu_0=0$,
$\trnormaln\bu_0=\eta_1\,\gamma(\eta_0)$ .

\begin{Definition} A couple $(\eta,\bu)$ is a weak solution of
  \eqref{eqn:fluid}, \eqref{eqn:shell}, and \eqref{eqn:data} for the
  admissible data $(\ff,g,\bu_0,\eta_0,\eta_1)$ in the intervall $I$
  if $\eta\in Y^I$ with $\norm{\eta}_{L^\infty(I\times M)}<\kappa$,
  $\eta(0,\cdot)=\eta_0$, $\bu\in X_\eta^I$ with $\tr \bu =
  \pa_t\eta\,\bnu$, and
  \begin{align}
    &- \int_I\iot \bu\cdot\pa_t\bphi\ dxdt - \int_I\im(\pa_t\eta)^2\,
    b\, \gamma(\eta)\ dAdt+\int_I\iot(\bu\cdot\nabla)\bu\cdot\bphi\
    dxdt\notag
    \\
    &+ \int_I\iot \nabla\bu:\nabla\bphi\ dxdt -\int_I\im\pa_t\eta\,
    \pa_tb\ dAdt + 2\int_I K(\eta,b)\ dt\label{eqn:schwach}
    \\
    &\hspace{0.5cm}=\int_I\iot \ff\cdot\bphi\ dxdt + \int_I\im g\, b\
    dAdt+\int_{\Omega_{\eta_0}}\bu_0\cdot\bphi(0,\cdot)\ dx +
    \im\eta_1\, b(0,\cdot)\ dA \notag
  \end{align}
  for all test functions $(b,\bphi)\in T_\eta^I$.\footnote{Note that
    for this definition to make sense we need $\bphi(0,\cdot)\in
    L^2(\Omega_{\eta_0})$. But this is a consequence of the fact that
    the extension of $\bphi$ by $(b\,\bnu)\circ q$ on $I\times
    B_\alpha$, $\norm{\eta}_{L^\infty(I\times M)}<\alpha<\kappa$, lies
    in $H^1(I,L^2(B_\alpha))$. This fact is proved in Remark
    \ref{bem:ausw}.}
\end{Definition}

The weak formulation \eqref{eqn:schwach} arises formally by
multiplication of \eqref{eqn:fluid} with a test function $\bphi$,
integration over space and time, integration by parts, and taking into
account \eqref{eqn:shell}. More precisely, we integrate by parts the
stress tensor, use Remark \ref{bem:korn}, replace the resulting
boundary integral by means of \eqref{eqn:shell}$_1$, and exploit the
identity $(\grad_{L^2}K(\eta),b)_{L^2}=2K(\eta,b)$.  Furthermore, we
integrate by parts with respect to time the terms containing the first
time-derivative of $\bu$ and the second time-derivative of $\eta$.
With regard to the $\bu$-term we obtain a boundary integral that can
be calculated using Reynold's transport theorem \ref{theorem:reynolds}.
With $\bv=\bu=(\pa_t\eta\bnu)\circ \Phi_{\eta(t)}^{-1}$ and
$\xi=\bu\cdot\bphi$ we get
\begin{equation*}
 \begin{aligned}
   \frac{d}{dt}\iot\bu\cdot\bphi\ dx &= \iot \pa_t\bu\cdot\bphi\ dx +
   \iot \bu\cdot\pa_t\bphi\
   dx\\
   &\hspace{0.5cm}+\int_{\pa\Omega_{\eta(t)}}\bu\cdot\bphi\
   \big((\pa_t\eta\,\bnu)\circ\Phi_{\eta(t)}^{-1}\big)\cdot\bnu_{\eta(t)}\
   dA_{\eta(t)}.
 \end{aligned}
\end{equation*}
By Remark \ref{bem:trafo} the boundary integral equals
\[\im(\pa_t\eta)^2\, b\, \gamma(\eta)\ dA.\]
Hence, \eqref{eqn:schwach} follows. Now we want to make sure that the
third term in \eqref{eqn:schwach} is well-defined. It's not hard to
see that the extension of $\bphi$ by $(b\,\bnu)\circ q$ lies in
$H^1(I\times B_\alpha)$ for $\norm{\eta}_{L^\infty(I\times
  M)}<\alpha<\kappa$ so that by standard embedding results
we have $\bphi\in L^\infty(I,L^3(\oet))$. For $\bu\in
L^2(I,L^6(\oet))$ this would suffice for
$(\bu\cdot\nabla)\bu\cdot\bphi$ to be in $L^1(\Omega_\eta^I)$. But due
to low regularity of the boundary we merely have $\bu\in
L^2(I,L^r(\oet))$ for all $r<6$, see Corollary \ref{lemma:sobolev}.
For this reason the additional requirement $\bphi\in
L^\infty(I,L^4(\oet))$ is needed.

\begin{Remark}\label{bem:grenzen}
For each weak solution $(\eta,\bu)$, all test functions $(b,\bphi)\in T^I_\eta$, and almost all $t\in I$ we have
\begin{align}
  &- \int_{0}^t\ios \bu\cdot\pa_t\bphi\ dxds -
  \int_{0}^t\im(\pa_t\eta)^2\, b\, \gamma(\eta)\ 
  dAds + \int_{0}^t\ios(\bu\cdot\nabla)\bu\cdot\bphi\ dxds\notag \\
  & + \int_{0}^t\ios \nabla\bu:\nabla\bphi\ dxds 
  -\int_{0}^t\im\pa_t\eta\, \pa_tb\ dAds + 2\int_{0}^t K(\eta,b)\ ds \notag \\
  &\hspace{.7cm}=\int_{0}^t\ios \ff\cdot\bphi\ dxds + \int_{0}^t\im
  g\, b\ dAds +\int_{\Omega_{\eta_0}}\bu_0\cdot\bphi(0,\cdot)\
  dx\label{eqn:schwacht}\\
  &\hspace{1.2cm}+ \im\eta_1\, b(0,\cdot)\ dA
  -\int_{\Omega_{\eta(t)}}\bu(t,\cdot)\cdot\bphi(t,\cdot)\ dx -
  \im\pa_t\eta(t,\cdot)\, b(t,\cdot)\ dA.\notag 
\end{align}
In fact we can abandon the assumption $b(T,\cdot)=0$, $\bphi(T,\cdot)=0$, because in order to show \eqref{eqn:schwacht}
we use the test functions $(b\,\rho^t_\epsilon,\bphi\,\rho^t_\epsilon)$ in \eqref{eqn:schwach}. Here $\rho\in C^\infty(\setR)$,
$\rho(s)=1$ for $s\le0$, $\rho(s)=0$ for $s\ge 1$, and $\rho^t_\epsilon(s)=\rho(\epsilon^{-1}(s-t))$. Then
\begin{equation*}
 \begin{aligned}
   &-\int_I\im\pa_t\eta\, \pa_s(b\rho^t_\epsilon)\ dAds
   \\
   &= -\int_I\rho^t_\epsilon\im\pa_t\eta\, \pa_tb\
   dAds-\int_I\epsilon^{-1}\rho'(\epsilon^{-1}(s-t))\im\pa_t\eta\, b\
   dAds
   \\
   &\quad \rightarrow -\int_0^t\im\pa_t\eta\, \pa_tb\
   dAds+\im\pa_t\eta(t,\cdot)\, b(t,\cdot)\ dA
 \end{aligned}
\end{equation*}
for $\epsilon \rightarrow 0$ and almost all $t\in I$. The convergence of the first term follows by the dominated convergence
theorem, while for the second term we used Lebesgue's differentiation theorem and
the identity $\int_\setR \rho'\ ds=-1$. Analogously the convergence of the other terms in
\eqref{eqn:schwach} can be shown, so that we obtain \eqref{eqn:schwacht}.
\end{Remark}

Now we are can state the main result of this paper.
\begin{Theorem}\label{theorem:hs}
  For arbitrary admissible data $(\ff,g,\bu_0,\eta_0,\eta_1)$ there
  exist a time $T^*\in (0,\infty]$ and a couple $(\eta,\bu)$ such that
  for all $T<T^*$ $(\eta,\bu)$ is a weak solution of
  \eqref{eqn:fluid}, \eqref{eqn:shell}, and \eqref{eqn:data} in the
  intervall $I=(0,T)$. Furthermore we have\footnote{Here, in
    $X_\eta^I$ we use the equivalent norm
    $\norm{\,\cdot\,}_{L^\infty(I,L^2(\Omega_{\eta(t)}))} +
    \norm{\nabla\,\cdot\,}_{L^2(\Omega_{\eta}^I)}$.}
  \begin{equation}\label{ab:hs}
    \begin{aligned}
      \norm{\eta}_{Y^I}^2 + \norm{\bu}_{X_\eta^I}^2 &\le c\,e^T\,\Big(\norm{\bu_0}_{L^2(\Omega_{\eta_0})}^2 + \mnorm{\eta_1}^2 +
  \norm{\eta_0}_{H^2(M)}^2\\
&\hspace{2cm} + \int_0^T
  \norm{\ff(s,\cdot)}_{L^2(\Omega_{\eta(s)})}^2 + \mnorm{g(s,\cdot)}^2 ds\Big).\\
  \end{aligned}
  \end{equation}
Either $T^*=\infty$ or $\lim_{t\rightarrow T^*}\norm{\eta(t,\cdot)}_{L^\infty(M)}=\kappa$.
\end{Theorem}

In the following we will denote the right hand side of \eqref{ab:hs}
as a function of $T$, $\Omega_\eta^I$, and the data by
$c_0(T,\Omega_\eta^I,\ff,g,\bu_0,\eta_0,\eta_1)$.

\subsection{Compactness}
\label{subsec:21}

An important step in the proof of Theorem \ref{theorem:hs} consists in
showing that every sequence of approximate weak solutions that is
bounded in the energy norm is relatively compact in $L^2$. In case of
the Navier-Stokes system on a spacetime cylinder this is usually
achieved by applying the Aubin-Lions theorem. At first sight this
seems to be possible for our coupled system, too. In fact, assuming
that the test functions vanish at $t=0$ equation \eqref{eqn:schwach}
takes the form
\begin{equation*}
  \begin{aligned}
    - \int_I&\iot \bu\cdot\pa_t\bphi\ dxdt - \int_I\im(\pa_t\eta)^2\,
    b\, \gamma(\eta)\
    dAdt -\int_I\im\pa_t\eta\, \pa_tb\ dAdt\\
    & = -\int_I\iot(\bu\cdot\nabla)\bu\cdot\bphi\ dxdt - \int_I\iot
    \nabla\bu:\nabla\bphi\ dxdt
    - 2\int_I K(\eta,b)\ dt\\
    &\hspace{0.5cm} +\int_I\iot \ff\cdot\bphi\ dxdt + \int_I\im g\, b\
    dAdt.
  \end{aligned}
\end{equation*}
It seems natural to interpret the right hand side as the definition of
the time-derivative of the couple $(\pa_t\eta,\bu)$ and to conclude from the spatial regularity of the solution that this time-derivative is bound\-ed in some dual space. Following the Aubin-Lions theorem this
should suffice to obtain $L^2$-compactness. It is rather difficult to
make this idea precise, though. To begin with, it is not a trivial
task to construct substitutes for dual space-valued Bochner spaces on
non-cylindrical spacetime domains of low regularity. Moreover, it is
not clear in which strict sense one could speak of a time-derivative
of the couple $(\pa_t\eta,\bu)$. Finally, it is important to note that
the function spaces involved depend on the solution, and hence we have
to deal with sequences of (dual) spaces. In view of these difficulties
another approach is chosen in \cite{b27,b28}; the system is
tested with difference quotients with respect to time. However, the construction of
suitable test functions is subtle, due to the time-dependent domain. Whereas the construction is
possible for the simpler case of a plate, it seems to be very
difficult for a general shell. For this reason we turn back
towards the former, more natural approach. In fact, it is possible to transfer some
of the proofs of the Aubin-Lions theorem to the present context. We
will neither provide substitutes for dual space-valued Bochner spaces
nor give a precise meaning to the time-derivative. Instead, we will
avoid abstract notions and work directly with the weak formulation. In
order to make the argument more accessible we now give the proof of
the Aubin-Lions theorem upon which we build our approach. Later on, we
will refer to this proof.

\begin{Proposition}\label{theorem:aubinlions}
Let $I\subset\setR$ be an open, bounded intervall and $1\le p\le\infty$, $1<r\le\infty$.
Assume that for the Banach spaces $B_0$, $B$, and $B_1$ we have
\[B_0\compactembedding B\embedding B_1.\]
Then
\[W:=\big\{v\in L^p(I,B_0)\ \big|\ v'\in L^r(I,B_1)\big\}\compactembedding L^p(I,B).\]
\end{Proposition}
\proof Let $(v_n)\subset W$ be bounded. It suffices to show that a subsequence converges in $C(\bar
I,B_1)$, since an application of the Ehrling lemma shows that for each $\epsilon>0$ and all $n,m\in\setN$ we have
\[\norm{v_n-v_m}_{L^p(I,B)}\le \epsilon\,\norm{v_n-v_m}_{L^p(I,B_0)} + c(\epsilon)\,
\norm{v_n-v_m}_{L^p(I,B_1)}.\]
Obviously $(v_n)$ is bounded in $C^{0,1-1/r}(\bar I,B_1)$. Hence, we can infer the convergence of a subsequence in $C(\bar
I,B_1)$ from the Arzela-Ascoli theorem, provided that the sequence $(v_n(t))_n$, $t$ from a dense subset of $I$, is relatively
compact in $B_1$. But this follows from the embeddings
\begin{equation*}
 \begin{aligned}
W\embedding C(\bar I,(B_0,B_1)_{\theta,1/\theta})\ \text{ and }\
(B_0,B_1)_{\theta,1/\theta}\compactembedding B_1
 \end{aligned}
\end{equation*}
which hold for a suitable $0<\theta<1$; see Theorem $33$ in \cite{b59}.
\qed\medskip
The following
proposition will not be used in the proof of Theorem \ref{theorem:hs},
but we will transfer its proof almost literally to
situations occuring in the proof of Theorem \ref{theorem:hs}.

\begin{Proposition}\label{lemma:komp}
Let $(\ff,g,\bu_0^n,\eta_0^n,\eta_1^n)$ be a sequence of admissible data with
\begin{equation}\label{eqn:abeschraenkt}
 \begin{aligned}
\sup_n\big(\tau(\eta_0^n)+\norm{\eta_0^n}_{H^2_0(M)}+\norm{\eta_1^n}_{L^2(M)}+\norm{
\bu^n_0}_{
L^2(\Omega_ {\eta_0^n})}\big)<\infty,
 \end{aligned}
\end{equation} 
and let $(\eta_n,\bu_n)$ be a sequence of weak solutions of \eqref{eqn:fluid},
\eqref{eqn:shell}, and \eqref{eqn:data} for the above data in the intervall $I=(0,T)$ such that
\begin{equation}\label{ab:beschraenkt}
 \begin{aligned}
\sup_n\big(\tau(\eta_n) + \norm{\eta_n}_{Y^I} +
\norm{\bu_n}_{X_{\eta_n}^I}\big)<\infty.
 \end{aligned}
\end{equation}
Then the sequence $(\pa_t\eta_n,\bu_n)$ is relatively compact in $L^2(I\times M)\times L^2(I\times
\setR^3)$.\footnote{Here and throughout the rest of the paper we tacitly extend functions defined in a domain of $\setR^3$ by $0$ to the whole space.}
\end{Proposition}
\proof
We infer from \eqref{ab:beschraenkt} that for a subsequence\footnote{When passing over to a subsequence we will
tacitly always do so with respect to all involved sequences and use again the subscript $n$.} we have
\begin{equation}\label{eqn:schwkonv}
 \begin{aligned}
   \eta_n&\rightarrow\eta &&\text{ weak$^*$ in
   }L^\infty(I,H^2_0(M))\text{ and uniformly},\\ 
   \pa_t\eta_n&\weakastto\pa_t\eta &&\text{ weak$^*$ in }L^\infty(I,L^2(M)),\\
   \bu_n&\weakastto\bu &&\text{ weak$^*$ in } L^\infty(I,L^2(\setR^3)),\\
   \nabla\bu_n&\weakto\xi &&\text{ weak in } L^2(I\times \setR^3).
\end{aligned}
\end{equation}
Here, we extend the functions $\nabla\bu_n$, which a-priori are defined
only in $\Omega_{\eta_n}^I$, by $\bfzero$ to $I\times \setR^3$. It is easy
to verify that the limit $\xi$ is identical $\nabla\bu$, provided that
we extend this function by $\bfzero$, too. The uniform convergence of
$(\eta_n)$ is due to the embedding $(1/2<\theta<1)$
\[
Y^I\embedding C^{0,1-\theta}(\bar I, C^{0,2\theta -1}(\pa\Omega))
\compactembedding C(\bar I\times\pa\Omega).
\] 
Assuming that we can show that
\begin{equation}\label{eqn:finalkonv}
 \begin{aligned}
   \int_I&\int_{\Omega_{\eta_n(t)}}|\bu_n|^2\ dxdt +
   \int_I\im|\pa_t\eta_n|^2\ dAdt\\ 
   &\rightarrow\int_I\int_{\Omega_{\eta(t)}}|\bu|^2\ dxdt +
   \int_I\im|\pa_t\eta|^2\ dAdt,
 \end{aligned}
\end{equation}
the proposition follows by means of \eqref{eqn:schwkonv} and the
trivial identity
\begin{equation*}
  \begin{aligned}
    \int_I&\int_{\setR^3} |\bu_n-\bu|^2\ dxdt + \int_I\im
    |\pa_t\eta_n-\pa_t\eta|^2\ dAdt \\ 
    & = \int_I\int_{\Omega_{\eta_n(t)}} |\bu_n|^2\ dxdt + \int_I\im
    |\pa_t\eta_n|^2\ dAdt + \int_I\int_{\Omega_{\eta(t)}}
    |\bu|^2\ dxdt\\
    &\hspace{5mm} + \int_I\im |\pa_t\eta|^2\ dAdt -
    2\int_I\int_{\setR^3} \bu_n\cdot\bu\ dxdt - 2\int_I\im
    \pa_t\eta_n\ \pa_t\eta\ dAdt.
  \end{aligned}
\end{equation*}
The convergence in \eqref{eqn:finalkonv} in turn is a consequence of
\begin{equation}\label{eqn:l2konv1}
  \begin{aligned}
    \int_I\int_{\Omega_{\eta_n(t)}}\bu_n\cdot \F_{\eta_n}\pa_t\eta_n\
    dxdt
    +\int_I\im|\pa_t\eta_n|^2\ dAdt\\
    \rightarrow \int_I\int_{\Omega_{\eta(t)}}\bu\cdot
    \F_{\eta}\pa_t\eta\ dxdt + \int_I\im|\pa_t\eta|^2\ dAdt
  \end{aligned}
\end{equation}
and
\begin{equation}\label{eqn:l2konv2}
  \begin{aligned}
    \int_I\int_{\Omega_{\eta_n(t)}} \bu_n\cdot(\bu_n -
    \F_{\eta_n}\pa_t\eta_n)\ dxdt \rightarrow\int_I
    \int_{\Omega_{\eta(t)}} \bu\cdot(\bu - \F_{\eta}\pa_t\eta)\ dxdt.
  \end{aligned}
\end{equation}
Here, we assume that the number $\alpha$ in the defintion of $\F$, see
Proposition \ref{lemma:FortVonRandZeit}, satisfies the inequality
$\sup_n\norm{\eta_n}_{L^\infty(I\times M)}<\alpha<\kappa$.

The proofs of \eqref{eqn:l2konv1} and \eqref{eqn:l2konv2} proceed
quite similarly. Let us start with \eqref{eqn:l2konv1}. An arbitrary
function $b\in H^2_0(M)$ in general does not meet the mean value
condition \eqref{eqn:integralnull} with respect to $\eta_n(t,\cdot)$. Hence, it
is not extendible to a divergence-free function in
$\Omega_{\eta_n(t)}$. That's why we need the operators
$\M_{\eta_n}$ from Lemma \ref{lemma:mittelwert}. From this lemma,
Proposition \ref{lemma:FortVonRandZeit}, and \eqref{ab:beschraenkt} we
deduce that
\begin{equation}\label{ab:b}
 \begin{aligned}
   &\norm{\M_{\eta_n}b}_{H^1(I,L^2(M))\cap
     L^2(I,H^2_0(M))}+\norm{\F_{\eta_n}\M_{\eta_n}b}_{H^1(I,L^2(B_\alpha))\cap
     C(\bar I,H^1(B_\alpha))}\\
   &\hspace{7.6cm}\le c\,\norm{b}_{H^2_0(M)}.
 \end{aligned}
\end{equation}
Examining \eqref{eqn:schwacht} with the weak solutions
$(\eta_n,\bu_n)$ and $(\M_{\eta_n}b,\F_{\eta_n}\M_{\eta_n}b)\in
T_{\eta_n}^I$ as test functions we infer from \eqref{ab:beschraenkt}
and \eqref{ab:b} that the integrands of the involved integrals over
time are bounded in $L^{12/11}(I)$ indepedently of $b$ and $n$,
provided that $\norm{b}_{H^2_0(M)}\le 1$.\footnote{From now on will write
  ``independently of $\norm{b}_{H^2_0(M)}\le 1$'' instead of
  ``independently of $b$ provided that $\norm{b}_{H^2_0(M)}\le 1$''.} For
the convective term (for the other terms it is obvious) this claim
follows from
\begin{equation*}
 \begin{aligned}
&\Bignorm{\int_{\Omega_{\eta_n(t)}}(\bu_n\cdot\nabla)\bu_n\cdot\F_{\eta_n}\M_{\eta_n}b\ dx}_{L^{12/11}(I)}\\
&\quad\le \norm{\bu_n}_{L^{12/5}(I,L^4(\Omega_{\eta_n(t)}))}\,
\norm{\nabla\bu_n}_{L^{2}(\Omega_{\eta_n}^I)}\,
\norm{\F_{\eta_n}\M_{\eta_n}b}_{L^{\infty}(I,L^4(\Omega_{\eta_n(t)}))}  
 \end{aligned}
\end{equation*}
in view of the embedding
\[H^1(\Omega_{\eta_n(t)})\embedding L^5(\Omega_{\eta_n(t)}),\]
which holds uniformly in $n$ and $t$, see Corollary \ref{lemma:sobolev}, and the interpolation embedding ($\theta=5/6$)
\[L^\infty(I,L^2(\Omega_{\eta_n(t)}))\cap L^2(I,L^5(\Omega_{\eta_n(t)}))\embedding
L^{12/5}(I,L^4(\Omega_{\eta_n(t)})).\]
In fact, the boundedness in $L^{12/11}(I)$ is the analogue of the boundedness of $(v_n')$ in $L^r(I,B_1)$ in the Aubin-Lions
theorem. It follows that the first eight terms in \eqref{eqn:schwacht} are bounded in $C^{0,1/12}(\bar
I)$ independently of $\norm{b}_{H^2_0(M)}\le 1$ and $n$. As a consequence the same claim holds 
for the last two terms in this identity, i.e. for
\[c_{b,n}(t):=\int_{\Omega_{\eta_n(t)}}\bu_n(t,\cdot)\cdot(\F_{\eta_n}\M_{\eta_n}b)(t,\cdot)\ dx +
\im\pa_t\eta_n(t,\cdot)\ (\M_{\eta_n}b)(t,\cdot)\ dA,\]
since the ninth and the tenth term are real numbers that are bounded independently of
$\norm{b}_{H^2_0(M)}\le 1$ and $n$, due to \eqref{eqn:abeschraenkt}. In virtue of \eqref{eqn:schwkonv} and Lemma
\ref{lemma:konvergenzen} $(1.a)$, $(2.a)$ for fixed $b\in H^2_0(M)$ the sequence $(c_{b,n})_n$ converges to
\[c_{b}(t):= \int_{\Omega_{\eta(t)}}\bu(t,\cdot)\cdot(\F_\eta\M_\eta b)(t,\cdot)\ dx +
\im\pa_t\eta(t,\cdot)\, (\M_\eta b)(t,\cdot)\ dA.\]
in the distributional sense. An application of the Arzela-Ascoli theorem shows that the convergence is, in fact, uniform in $\bar
I$.

The next step is to show that this uniform convergence is independent of $\norm{b}_{H^2_0(M)}\le 1$, i.e. the functions
\begin{equation*}
h_n(t):=\sup_{\norm{b}_{H^2_0(M)}\le
1} \big(c_{b,n}(t)-c_{b}(t)\big)
\end{equation*} 
converge to $0$ uniformly in $\bar I$. This claim is the analogue of the convergence of $(v_n)$ in $C(\bar I,B_1)$ in the proof of
the Aubin-Lions theorem. In view of \eqref{ab:beschraenkt} we have for almost all $t\in I$
\[\sup_n \big(\norm{\bu_n(t)}_{L^2(\setR^3)}+\norm{\pa_t\eta_n(t)}_{L^2(M)}\big)<\infty.\] 
By a diagonal sequence argument we deduce that there exist a
countable, dense subset $I_0$ of $I$ and a subsequence
$(\eta_n,\bu_n)$ such that for all $t\in I_0$ the sequences
$(\bu_n(t,\cdot))_n$ and $(\pa_t\eta_n(t,\cdot))_n$ converge weakly in
$L^2(\setR^3)$ respectively $L^2(M)$. Now we show that for fixed $t\in
I_0$ the sequence $(c_{b,n}(t))_n$ converges independently of
$\norm{b}_{H^2_0(M)}\le 1$. Let $\eta^*$ denote the weak limit of
$(\pa_t\eta_n(t,\cdot))_n$ in $L^2(M)$. Then we deduce this claim for
the second term in $c_{b,n}(t)$ from the estimate
\begin{equation*}
 \begin{aligned}
   &\quad \big|\im \pa_t\eta_n(t,\cdot)\, (\M_{\eta_n}
   b)(t,\cdot)-\eta^*\,(\M_\eta b)(t,\cdot)\
   dA\big|\\
   &\hspace{2.0cm}\le \big|\im\big(\pa_t\eta_n(t,\cdot)-\eta^*\big)\,
   (\M_{\eta_n} b)(t,\cdot)\
   dA\big|\\
   &\hspace{2.5cm}+ \big|\im\eta^*\,
   \big((\M_{\eta_n}b)(t,\cdot)-(\M_{\eta}
   b)(t,\cdot)\big)\ dA\big|\\
   &\hspace{2.0cm}\le \norm{\pa_t\eta_n(t,\cdot)-\eta^*}_{(H^1(M))'}\,
   \norm{(\M_{\eta_n}b)(t,\cdot)}_{H^1(M)}\\
   &\hspace{2.5cm}+\norm{\eta^*}_{L^2(M)}\norm{(\M_{\eta_n}b)(t,\cdot)
     -(\M_{\eta}b)(t,\cdot)}_{L^2(M)},
 \end{aligned}
\end{equation*}
taking into account \eqref{ab:b}, \eqref{eqn:einbschnurz}, Lemma
\ref{lemma:konvergenzen} $(1.a)$, and
\[
L^2(M)\compactembedding (H^1(M))'.
\] 
Analogously, let $\bu^*$ denote the weak limit of $(\bu_n(t,\cdot))_n$
in $L^2(\setR^3)$. Then the claim for the first term in $c_{b,n}(t)$
follows from the estimate
\begin{equation*}
 \begin{aligned}
   &\quad\big|\int_{\Omega_{\eta_n(t)}}
   \bu_n(t,\cdot)\cdot(\F_{\eta_n}\M_{\eta_n} b)(t,\cdot)\
   dx-\int_{\Omega_{\eta(t)}} \bu^*\cdot(\F_{\eta}\M_{\eta}
   b)(t,\cdot)\ dx\big|\\
   &\hspace{2.0cm}\le
   \big|\int_{B_\alpha}\big(\bu_n(t,\cdot)-\bu^*\big)
   \cdot(\F_{\eta_n}\M_{\eta_n}
   b)(t,\cdot)\ dx\big|\\
   &\hspace{2.5cm}+
   \big|\int_{B_\alpha}\bu^*\cdot\big((\F_{\eta_n}\M_{\eta_n}b)
   (t,\cdot)-(\F_{\eta}\M_{\eta}
   b)(t,\cdot)\big)\ dx\big|\\
   &\hspace{2.0cm}\le \norm{\bu_n(t,\cdot)-\bu^*}_{(H^1(B_\alpha))'}\,
   \norm{(\F_{\eta_n}\M_{\eta_n}b)(t,\cdot)}_{H^1(B_\alpha)}\\
   &\hspace{2.5cm}+\norm{\bu^*}_{L^2(B_\alpha)}\norm{(\F_{\eta_n}\M_{\eta_n}b)
     (t,\cdot)-(\F_{\eta}\M_ {\eta}b)(t, \cdot)}_{L^2(B_\alpha)},
 \end{aligned}
\end{equation*}
taking into account \eqref{ab:b}, Lemma \ref{lemma:konvergenzen} $(2.a)$, and
\begin{equation}\label{eqn:keinb}
 L^2(B_\alpha)\compactembedding (H^1(B_\alpha))'.
\end{equation}
At last we show that the convergence of $(c_{b,n}(t))_n$, independent of $\norm{b}_{H^2_0(M)}\le
1$, does not merely hold for $t\in I_0$ but uniformly for all $t\in \bar I$. Due to the uniform boundedness of
$(c_{b,n})$ in $C^{0,1/12}(\bar I)$ we have for all $t,t'\in \bar I$ and $n,m\in\setN$
\begin{equation*}
 \begin{aligned}
  |c_{b,n}(t)-c_{b,m}(t)|&\le |c_{b,n}(t)-c_{b,n}(t')| + |c_{b,n}(t')-c_{b,m}(t')| +
|c_{b,m}(t')-c_{b,m}(t)|\\
&\le c\,|t-t'|^{1/12} + |c_{b,n}(t')-c_{b,m}(t')|. 
 \end{aligned}
\end{equation*}
For any given $\epsilon>0$ we can find a finite set $I_0^\epsilon\subset I_0$ such that for each $t\in \bar I$ there exists a
$t'\in I_0^\epsilon$ with $c\,|t-t'|^{1/12}<\epsilon/2.$ Furthermore we just showed that
$|c_{b,n}(t')-c_{b,m}(t')|<\epsilon/2$ provided that
$t'\in I_0^\epsilon$ and $n,m\ge N$, where $N$ depends on $\epsilon$ but not on $t'$ and
$\norm{b}_{H^2_0(M)}\le 1$. Hence, $(c_{b,n})_n$ converges to $c_b$ uniformly in $\bar I$ and independently
of $\norm{b}_{H^2_0(M)}\le 1$, i.e. $(h_n)$
converges to $0$ uniformly in $\bar I$.

Finally, we have to apply some argument of the type of the Ehrling lemma, namely Lemma \ref{lemma:ehrling}. Letting
\begin{equation*}
g_n(t):=\sup_{\norm{b}_{L^2(M)}\le
1}\big(c_{b,n}(t)-c_{b}(t)\big)
\end{equation*}
we infer from this lemma that for each $\epsilon>0$ there exists a constant $c(\epsilon)$ such that
\begin{equation*}
 \begin{aligned}
  \int_I g_n(t)\ dt \le \epsilon\,
c\,\big(\norm{\bu_n}_{L^2(I,H^1(\Omega_{\eta_n(t)}))}+\norm{\bu}_{L^2(I,H^1(\Omega_{\eta(t)}))}
\big) + c(\epsilon)\int_I h_n(t)\ dt.
 \end{aligned}
\end{equation*}
In view of the uniform convergence of $(h_n)$ to $0$ we conclude that
\begin{equation}\label{eqn:gn}
  \lim_n \int_I g_n(t)\ dt = 0.
\end{equation}
Adding a zero-sum we obtain the identity
\begin{align}
  &\hspace{0.5cm} \int_{\Omega_{\eta_n(t)}}\bu_n(t,\cdot)\cdot
  (\F_{\eta_n}\pa_t\eta_n)(t,\cdot)\ dx
  +\im|\pa_t\eta_n(t,\cdot)|^2\ dA\notag \\
  &\hspace{0.5cm}-\int_{\Omega_{\eta(t)}}\bu(t,\cdot)\cdot
  (\F_{\eta}\pa_t\eta)(t,\cdot)\ dx-
  \im|\pa_t\eta(t,\cdot)|^2\ dA\label{eqn:nulladd}\\
  &\hspace{0.5cm}= \int_{\Omega_{\eta_n(t)}}\bu_n(t,\cdot)\cdot
  (\F_{\eta_n}\pa_t\eta_n)(t,\cdot)\ dx +
  \im\pa_t\eta_n(t,\cdot)\, \pa_t\eta_n(t,\cdot)\ dA\notag \\
  &\hspace{1cm}-\int_{\Omega_{\eta(t)}}\bu(t,\cdot)\cdot
  (\F_{\eta}\M_{\eta}\pa_t\eta_n)(t,\cdot)\ dx
  -\im\pa_t\eta(t,\cdot)\,
  (\M_{\eta}\pa_t\eta_n)(t,\cdot)\ dA\notag \\
  &\hspace{1cm}+ \int_{\Omega_{\eta(t)}}\bu(t,\cdot)\cdot
  \big((\F_{\eta}\M_{\eta}\pa_t\eta_n)(t,\cdot)
  -\F_{\eta}\pa_t\eta)(t,\cdot)\big)\ dx\notag \\ 
  &\hspace{1cm} + \im\pa_t\eta(t,\cdot)\,
  \big((\M_\eta\pa_t\eta_n)(t,\cdot)-\pa_t\eta(t,\cdot)\big)\ dA.\notag \
\end{align}
Due to $\M_{\eta_n}\pa_t\eta_n=\pa_t\eta_n$ the first two lines of the
right hand side are equal to $c_{b,n}(t)-c_{b}(t)$ with
$b=\pa_t\eta_n(t,\cdot)$. In view of \eqref{ab:beschraenkt} their
absolute value is dominated by $c\,g_n(t)$ for almost all $t$.
Integrating \eqref{eqn:nulladd} over $I$ and using \eqref{eqn:gn} as
well as the weak convergences of $(\M_\eta \pa_t\eta_n)$ and
$(\F_\eta\M_\eta\pa_t\eta_n)$ in $L^2(I\times M)$ respectively
$L^2(I\times B_\alpha)$, which follow from the weak convergence of
$(\pa_t\eta_n)$ in $L^2(I\times M)$, we obtain \eqref{eqn:l2konv1};
note that $\M_{\eta}\pa_t\eta=\pa_t\eta$.

Let us proceed with the proof of \eqref{eqn:l2konv2}. Let
$\sigma>0$ be sufficiently small and $\delta_\sigma\in C^4(\bar
I\times \pa\Omega)$ with $\norm{\delta_\sigma-\eta}_{L^\infty(I\times
  \pa\Omega)}<\sigma$ and $\delta_\sigma<\eta$ in $\bar
I\times\pa\Omega$. For $\bphi\in H(\Omega)$ we define
\[c_{\bphi,n}^\sigma(t):=\int_{\Omega_{\eta_n(t)}}\bu_n(t,\cdot)\cdot \T_{\delta_\sigma(t)}\bphi\
dx,\quad
c_{\bphi}^\sigma(t):=\int_{\Omega_{\eta_n(t)}}\bu(t,\cdot)\cdot \T_{\delta_\sigma(t)}\bphi\ dx;\]
see Remarks \ref{bem:tdelta} and \ref{bem:tdelta2}. Now, we can proceed as before to show that the functions
\[h_n^\sigma(t):=\sup_{\norm{\bphi}_{H^1_{0,\dv}(\Omega)}\le
1} \big(c_{\bphi,n}^\sigma(t)-c_{\bphi}^\sigma(t)\big)\]
converge to $0$ uniformly in $\bar I$. For $\bphi\in H^1_{0,\dv}(\Omega)$ we extend
the function $\T_{\delta_\sigma}\bphi$, which a-priori is defined only in $\Omega_{\delta_\sigma}^I$, by $\boldsymbol{0}$ to
$I\times B_\alpha$. Then
\begin{equation}\label{eqn:schnarch}
 \begin{aligned}
\norm{\T_{\delta_\sigma}\bphi}_{H^1(I,L^2(B_\alpha))\cap C(\bar I,H^1(B_\alpha))}\le
c\,\norm{\bphi}_{H^1_{0,\dv}(\Omega)};
 \end{aligned}
\end{equation}
see Remarks \ref{bem:tdelta} and \ref{bem:tdelta2}. The identity \eqref{eqn:schwacht} with the weak solutions $(\eta_n,\bu_n)$ and the test functions
$(0,\T_{\delta_\sigma}\bphi)\in T^I_{\eta_n}$ ($n$ sufficiently large) shows in view of \eqref{ab:beschraenkt} and
\eqref{eqn:schnarch} that the functions $c_{\bphi,n}^\sigma$ are bounded in $C^{0,1/12}(\bar I)$ independently of
$\norm{\bphi}_{H^1_{0,\dv}(\Omega)}\le 1$
and $n$. From \eqref{eqn:schwkonv} and the Arzela-Ascoli theorem we obtain as before the convergence of $(c_{\bphi,n}^\sigma)_n$
to $c_{\bphi}^\sigma$ uniformly $\bar I$. Futhermore, from \eqref{eqn:keinb} we deduce that for $t\in I_0$ the
sequence $(c_{\bphi,n}^\sigma(t))_n$ converges independently of $\norm{\bphi}_{H^1_{0,\dv}(\Omega)}\le 1$. Again, we
conclude that $(h_n^\sigma)_n$ converges to $0$ uniformly.

An application of Lemma \ref{lemma:ehrling} shows that
\begin{equation}\label{eqn:gsigma}
\lim_n\int_I g^\sigma_n(t)\ dt=0,
\end{equation}
where
\[g_n^\sigma(t):=\sup_{\norm{\bphi}_{H(\Omega)}\le
1} \big(c_{\bphi,n}^\sigma(t)-c_{\bphi}^\sigma(t)\big).\]
Due to \eqref{ab:beschraenkt} the $L^2(\Omega_{\eta_n(t)})$ norms of
$\bu_n(t,\cdot)-(\F_{\eta_n}\pa_t\eta_n)(t,\cdot)\in
H(\Omega_{\eta_n(t)})$ are bounded for almost all $t$, independently of $t$ and
$n$. By Lemma \ref{lemma:divdichtglm}, for each $\epsilon > 0$ there exists a
$\sigma>0$ such that for almost all $t$ and all sufficiently large $n$ there are functions
$\bpsi_{t,n}\in H(\Omega_{\eta_n(t)})$ having $L^2(\Omega_{\eta_n(t)})$ norms bounded independently of $t$ and $n$, $\supp
\bpsi_{t,n}\subset
\Omega_{\delta_\sigma(t)}$ and
\begin{equation*}
  \norm{\bu_n(t,\cdot)-(\F_{\eta_n}\pa_t\eta_n)(t,\cdot) -
\bpsi_{t,n}}_{(H^{1/4}(\setR^3))'} <
\epsilon,
\end{equation*}
in particular $\bpsi_{t,n}\in H(\Omega_{\delta_\sigma(t)})$ with norms bounded independently of $t$ and $n$. By adding a zero-sum
we obtain the identity
\begin{equation}\label{eqn:nulladd2}
 \begin{aligned}
   & \int_{\Omega_{\eta_n(t)}} \bu_n(t,\cdot)\cdot\big(\bu_n(t,\cdot)
   -
   (\F_{\eta_n}\pa_t\eta_n)(t,\cdot)\big)\ dx \\
   &-\int_{\Omega_{\eta(t)}} \bu(t,\cdot)\cdot\big(\bu(t,\cdot) -
   (\F_{\eta}\pa_t\eta)(t,\cdot)\big)\
   dx\\
   &\hspace{.5cm}= \int_{\Omega_{\eta_n(t)}}
   \bu(t,\cdot)\cdot\big(\bu_n(t,\cdot) -
   (\F_{\eta_n}\pa_t\eta_n)(t,\cdot)\big)\, dx \\
   &\hspace{1cm}-\int_{\Omega_{\eta(t)}}
   \bu(t,\cdot)\cdot\big(\bu(t,\cdot) -
   (\F_{\eta}\pa_t\eta)(t,\cdot)\big)\,
   dx\\
   &\hspace{1cm} + \int_{\Omega_{\eta_n(t)}}
   \bu_n(t,\cdot)\cdot\bpsi_{t,n}\ dx -
   \int_{\Omega_{\eta(t)}} \bu(t,\cdot)\cdot\bpsi_{t,n}\ dx\\
   &\hspace{1cm}+ \int_{\Omega_{\eta_n(t)}}
   \bu_n(t,\cdot)\cdot\big(\bu_n(t,\cdot) -
   (\F_{\eta_n}\pa_t\eta_n)(t,\cdot)-
   \bpsi_{t,n}\big)\ dx\\
   &\hspace{1cm} - \int_{\Omega_{\eta(t)}}
   \bu(t,\cdot)\cdot\big(\bu_n(t,\cdot) -
   (\F_{\eta_n}\pa_t\eta_n)(t,\cdot)- \bpsi_{t,n}\big)\ dx.
 \end{aligned}
\end{equation}
The absolute value of the second line of the right-hand side is
dominated by $c\, g_n^\sigma(t)$, while the absolute values of the
last two lines can be estimated by
\[
\epsilon\, ( \norm{\bu_n(t,\cdot)}_{H^{1/4}(\setR^3)} +
\norm{\bu(t,\cdot)}_{H^{1/4}(\setR^3)}),
\] 
see Proposition \ref{lemma:nullfort}. Integrating \eqref{eqn:nulladd2}
over $I$ and using \eqref{eqn:gsigma} as well as the weak convergence
of $(\chi_{\Omega^I_{\eta_n}}(\bu_n - \F_{\eta_n}\pa_t\eta_n))$ to
$\chi_{\Omega^I_{\eta}}(\bu - F_{\eta}\pa_t\eta)$ in
$L^2(I\times\setR^3)$, which follows from \eqref{eqn:schwkonv} and
Lemma \ref{lemma:konvergenzen} $(2.c)$, we obtain \eqref{eqn:l2konv2}. This concludes the proof.
\qed

\subsection{The decoupled system}
Let us now construct approximate solutions. Since the domain of
definition of the solution depends on the solution
itself it's not possible to use a Galerkin approach directly. For this reason we shall decouple the
problem and use a fixed point argument. There seem to be several possibilities to do so. We choose one which preserves the coupling of the fluid with the structure. More precisely, for a given motion $\delta$ of the boundary with $\delta(0,\cdot)=\eta_0$ we seek functions
$\eta\in Y^I$ with $\norm{\eta}_{L^\infty(I\times M)}<\kappa$ and
$\eta(0,\cdot)=\eta_0$ as well as $\bu\in X_\delta^I$ with $\trd \bu =
\pa_t\eta\,\bnu$ satisfying
\begin{align}
  &- \int_I\int_{\Omega_{\delta(t)}} \bu\cdot\pa_t\bphi\ dxdt -
  \int_I\im\pa_t\eta\,\pa_t\delta\, b\, \gamma(\delta)\
  dAdt+\int_I\int_{\Omega_{\delta(t)}}(\bu\cdot\nabla)\bu\cdot\bphi\
  dxdt \notag \\
  &+\int_I\int_{\Omega_{\delta(t)}} \nabla\bu:\nabla\bphi\ dxdt
  -\int_I\im\pa_t\eta\, \pa_tb\ dAdt + 2\int_I K(\eta,b)\ dt \label{eqn:zwischen}\\
  &\hspace{.3cm}=\int_I\int_{\Omega_{\delta(t)}} \ff\cdot\bphi\ dxdt +
  \int_I\im g\, b\
  dAdt+\int_{\Omega_{\eta_0}}\bu_0\cdot\bphi(0,\cdot)\ dx +
  \im\eta_1\, b(0,\cdot)\ dA\notag
\end{align}
for all test functions $(b,\bphi)\in T_\delta^I$. In order to preserve the energy estimate we have to modify the convective term. Otherwise, an application of Reynold's transport theorem
\ref{theorem:reynolds} does not lead to the desired identity since the velocity of the boundary and the velocity of the fluid at the boundary do not coincide anymore. To solve this
problem we rewrite the convective term in the form
\begin{equation*}
  \begin{aligned}
    &\int_I\iot(\bu\cdot\nabla)\bu\cdot\bphi\ dxdt \\
    &\quad = \frac{1}{2}\int_I\iot(\bu\cdot\nabla)\bu\cdot\bphi\ dxdt -
    \frac{1}{2}\int_I\iot(\bu\cdot\nabla)\bphi\cdot\bu\ dxdt\\
    &\quad\quad +\frac{1}{2}\int_I\im(\pa_t\eta)^2\,b\,\gamma(\eta)\ dAdt
  \end{aligned}
\end{equation*}
and replace $\eta$ by the given function $\delta$ in the
appropriate places. Thus, instead of \eqref{eqn:zwischen} we obtain
\begin{align}
  &- \int_I\int_{\Omega_{\delta(t)}} \bu\cdot\pa_t\bphi\ dxdt -
  \frac{1}{2}\int_I\im\pa_t\eta\,\pa_t\delta\, b\,\gamma(\delta)\ dAdt
  \notag\\
  &+\frac{1}{2}\int_I\int_{\Omega_{\delta(t)}}(\bu\cdot\nabla)\bu\cdot\bphi\
  dxdt -
  \frac{1}{2}\int_I\int_{\Omega_{\delta(t)}}(\bu\cdot\nabla)\bphi\cdot\bu\
  dxdt \label{eqn:zwischen2}\\
  &+\int_I\int_{\Omega_{\delta(t)}} \nabla\bu:\nabla\bphi\ dxdt
  -\int_I\im\pa_t\eta\,\pa_t b\ dAdt +
  2\int_I K(\eta,b)\ dt \notag\\
  &\hspace{0.3cm}=\int_I\int_{\Omega_{\delta(t)}} \ff\cdot\bphi\ dxdt
  + \int_I\im g\ b\
  dAdt\notag +\int_{\Omega_{\eta_0}}\bu_0\cdot\bphi(0,\cdot)\ dx
  + \im\eta_1\, b(0,\cdot)\ dA.\notag
\end{align}
Testing this equation formally with $(\pa_t\eta,\bu)$, we see that
the two ``convective terms'' cancel each other while the first two
terms yield 
\[
  \int_I\frac{1}{2}\frac{d}{dt}\int_{\Omega_{\delta(t)}}|\bu|^2\ dxdt.
\]
For the fixed-point argument we need to regularize the motion of the domain as well as the convective term in \eqref{eqn:zwischen2}. For a detailed explanation see \cite{b62}. However, even with these regularisations the uniqueness of weak solutions of the decoupled system is not obvious. For this reason we shall employ a variant of Schauder's fixed-point theorem for multi-valued mappings, Theorem \ref{theorem:kakutani}. Hence, we need to linearise the ``convective
terms'' to ensure that the image sets of the multi are convex. The resulting weak formulation is essentially the one used in \cite{b27}.

Let us now construct appropriate regularization operators
$\R_\epsilon$, $\epsilon>0$. For a function $\omega\in
C_0^\infty(\setR^3)$ with $\int_{\setR^3}\omega\ dxdt=1$ and
$\supp\omega\subset\{(t,x)\in\setR^3\ |\ t>0\}$ we set
$\omega_\epsilon:=\epsilon^{-3}\omega(\epsilon^{-1}\cdot)$.  Let
$(\phi_k,U_k)_k$ be a finite atlas of $\pa\Omega$ with subordinate
partition of unity $(\psi_k)_k$, see for instance \cite{b1}.
Furthermore, we extend functions $\delta\in C(\bar I\times\pa\Omega)$
with $\delta(0,\cdot)=\eta_0$ by $\eta_0$ to
$(-\infty,T]\times\pa\Omega$. Let
\begin{equation}\label{eqn:reg}
  \R_\epsilon\delta:=\sum_k(\omega_\epsilon\ast((\psi_k\,
  \delta)\circ\phi^{-1}_k))\circ\phi_k +\epsilon^{1/2},
\end{equation}
where the summand with index $k$ is extended by $0$ to the complement
of $U_{k}$. Taking into account the support of $\omega$ we see that
$\R_\epsilon\eta_0:=(\R_\epsilon\delta)(0,\cdot)$ depends only on
$\eta_0$ (and $\epsilon$). Due to the H{\"o}lder continuity of $\eta_0$
we get $\R_\epsilon\eta_0\ge\eta_0$ for all $0<\epsilon\le c$ with a
constant $c=c(\eta_0)$. In fact, this follows from the estimate
\begin{align*}
  |\omega_\epsilon&\ast((\psi_k\, \delta)\circ\phi^{-1}_k)-(\psi_k\,
  \delta)\circ\phi^{-1}_k|(0,x)\\
  &=|\int_{\setR^4}\omega_\epsilon(-s,x-y)\Big(((\psi_k\,
  \eta_0)\circ\phi^{-1}_k)(y)-((\psi_k\,
  \eta_0)\circ\phi^{-1}_k)(x)\Big)\ dyds|\\
  &\le c\int_{\setR^4}|\omega_\epsilon(-s,x-y)||x-y|^{3/4}\ dyds=c\
  \epsilon^{3/4},
\end{align*}
where we used the $3/4$-H{\"o}lder continuity of $\eta_0$ and the
regularity of $\phi_k$, $\psi_k$. Note that the term
\[(\omega_\epsilon\ast((\psi_k\,
\delta)\circ\phi^{-1}_k))\circ\phi_k\] belongs to $C^4(\bar
I\times\pa\Omega)$ for sufficiently small $\epsilon$. By Euclidean
theory we deduce the uniform convergence of $(\R_\epsilon\delta)$ to
$\delta$. If $\delta$ possesses a time derivative $\pa_t\delta$ in \linebreak 
$L^2(I\times\pa\Omega)$, then we analogously deduce the convergence of
$(\pa_t\R_\epsilon\delta)$ to $\pa_t\delta$ in $L^2(I\times
\pa\Omega)$ using the identity
$\pa_t\R_\epsilon\delta=\R_\epsilon\pa_t\delta$. Finally, for
$\widetilde\omega\in C_0^\infty(\setR^4)$ with
$\int_{\setR^4}\widetilde \omega\ dx=1$ we set
$\widetilde\omega_\epsilon:=\epsilon^{-4}\widetilde\omega(\epsilon^{-1}\cdot)$.
For functions $v\in L^2(I\times\setR^3)$ we set $\R_\epsilon
v:=\widetilde\omega_\epsilon\ast v$ where again $v$ is extended by $0$
to $\setR^4$. $(\R_\epsilon v)$ converges to $v$ in
$L^2(I\times\setR^3)$.

In general, the domain of definition of the initial value $\bu_0$ is different from $\Omega_{\R_\epsilon\eta_0}$. Thus, we must
modify our initial values. Since $(\R_\epsilon\eta_0)$ approximates
the function $\eta_0$ from above it suffices to extend $\bu_0$ appropriately.
To this end, using Remark \ref{bem:nspur} we extend $\eta_1$ to a divergence-free $L^2$-vector field $\bpsi$ in $B_\alpha$, where $\alpha$
satisfies $\norm{\eta_0}_{L^\infty(M)}<\alpha<\kappa$, and define
\begin{align*}
  \bu_0^\epsilon:=\left\{\begin{array}{cl} \bu_0 \qquad & \text{in
      }\Omega_{\eta_0}, \\
      \bpsi \qquad & \text{in }
      \Omega_{\R_\epsilon\eta_0}\setminus\Omega_{\eta_0}\end{array}\right..
\end{align*}
Clearly, we have $\bu_0^\epsilon\in L^2(\Omega_{\R_\epsilon\eta_0})$,
and from \eqref{eqn:nspur}, $\trnormaln\bu_0=\eta_1\gamma(\eta_0)$, and Proposition
\ref{lemma:nspur} we deduce that $\bu_0^\epsilon$ is divergence-free.
Defining 
\[
  \eta_1^\epsilon:= \exp\Big(\int_{\eta_0}^{\R_\epsilon\eta_0}
  \beta(\cdot+\tau\bnu)\ d\tau\Big)\ \eta_1
\]
we also have $\trnormale
\bu_0^\epsilon=\eta_1^\epsilon\,\gamma(\R_\epsilon\eta_0)$. This can
be easily seen from the definition of $\bpsi$ and an approximation
argument as in the proof of Proposition \ref{lemma:FortVonRand}.
Extending $\bu_0$ and $\bu_0^\epsilon$ by $\boldsymbol{0}$ to
$\setR^3$ we obtain, using the uniform convergence of
$\R_\epsilon\eta_0\rightarrow\eta_0$ in $\pa\Omega$, the convergence
\begin{equation}\label{eqn:konv7}
  \begin{aligned}
    \bu_0^\epsilon\rightarrow\bu_0 \qquad &\text{ in } L^2(\setR^3),\\
    \eta_1^\epsilon\rightarrow\eta_1 \qquad&\text{ in } L^2(M).
  \end{aligned}
\end{equation}

In the following, let $I=(0,T)$, $T>0$, be a fixed time interval and $\bv\in
L^2(I\times \setR^3)$ and $\delta\in C(\bar I\times \pa\Omega)$ with $\norm{\delta}_{L^\infty(I\times\pa\Omega)}<\kappa$ and $\delta(0,\cdot)=\eta_0$ be arbitrary but fixed functions. Based on the above discussion we now define
our approximate problem. For the sake of a better readability we will
suppress for the  moment the parameter $\epsilon$ in the notation. In
particular stand $\bu_0$ and $\eta_1$ for the modified initial
conditions $\bu_0^\epsilon$ and $\eta_1^\epsilon$, respectively.

\begin{Definition}\label{def:ent}
  A couple $(\eta,\bu)$ is called weak solution of the decoupled and
  regularized problem with the data $(\delta,\bv)$ in the interval
  $I$ if $\eta\in Y^I$ with $\eta(0,\cdot)=\eta_0$, $\bu\in
  X_{\R\delta}^I$ with $\trrd \bu =\pa_t\eta\,\bnu$ satisfy 
  \begin{align}
    &- \int_I\iortd \bu\cdot\pa_t\bphi\ dxdt -
    \frac{1}{2}\int_I\im\pa_t\eta\,
    \pa_t\mathcal{R}\delta\, b\, \gamma(\mathcal{R}\delta)\ dAdt
    \notag \\
    &+\frac{1}{2}\int_I\iortd(\mathcal{R}\bv\cdot\nabla)\bu\cdot\bphi\
    dxdt -
    \frac{1}{2}\int_I\iortd(\mathcal{R}\bv\cdot\nabla)\bphi\cdot\bu\
    dxdt\label{eqn:ent} \\
    & + \int_I\iortd \nabla\bu:\nabla\bphi\ dxdt -\int_I\im\pa_t\eta\,
    \pa_tb\ dAdt + 2\int_I K(\eta,b)\
    dt\notag  \\
    &\hspace{0.2cm}=\int_I\iortd \ff\cdot\bphi\ dxdt + \int_I\im g\,
    b\ dAdt
    +\int_{\Omega_{\mathcal{R}\eta_0}}\bu_0\cdot\bphi(0,\cdot)\ dx +
    \im \eta_1\, b(0,\cdot)\ dA\notag 
  \end{align}
  for all test functions $(b,\bphi)\in T_{\R\delta}^I$.
\end{Definition}
\medskip

The existence of a solution of the decoupled and regularized problem
is shown in \cite{b27} in the case of a flat geometry by transforming
the problem to a spacetime cylinder and a subsequent 
Galerkin approximation.  We  show the existence of a solution without
the transformation of the problem. Nevertheless, we need a
diffeomorphism of the time varying  domain to a spacetime cylinder
for the construction of suitable basis functions. 

\begin{Proposition}\label{lemma:ent}
There exists a weak solution $(\eta,\bu)$ of the decoupled and
regularized problem with the data  $(\delta,\bv)$ in the interval $I$,
which satisfies the estimate 
\begin{equation}\label{ab:ent}
 \begin{aligned}
 &\norm{\eta}_{Y^I}^2 + \norm{\bu}_{X_{\R\delta}^I}^2\le
c_0(T,\Omega_{\R\delta}^I,\ff,g,\bu_0,\eta_0,\eta_1).
  \end{aligned}
\end{equation}
In particular, the left-hand side is bounded independently of the
parameter $\epsilon$ and the data $(\delta,\bv)$. 
\end{Proposition}

Before we can give a proof of this proposition we need to construct appropriate basis
functions for our Galerkin approach. To this end, we choose a basis 
$(\widehat\bX_k)_{k\in\setN}$ of $H^1_{0,\dv}(\Omega)$ and a basis
$(\widehat Y_k)_{k\in\setN}$ of the space
\begin{align*}
  \Big\{Y\in H^2_0(M)\big|\int_M Y\ dA = 0\Big\}.
\end{align*}
Theorem \ref{theorem:stokes} ensures the existence of divergence-free
vector fields $\widehat\bY_k$ solving the Stokes system in $\Omega$ with
boundary data $\widehat Y_k\,\bnu$ (as usual extended by $0$ to
$\pa\Omega$).  For $t\in\bar I$  we set (see Remarks \ref{bem:tdelta} and \ref{bem:tdelta2})
\begin{align*}
  \bX_k(t,\cdot):=\T_{\R\delta(t)}\widehat\bX_k,\quad
  \bY_k(t,\cdot):=\T_{\R\delta(t)}\widehat\bY_k.
\end{align*}
The vector fields $\bX_k(t,\cdot)$ obviously form a basis of
$H^1_{0,\dv}(\Omega_{\R\delta(t)})$. Note that for $q\in\pa\Omega$
the differential $d\Psi_{\R\delta(t)}(q)$ merely scales the normal
$\bnu(q)$. Thus, the definition 
\begin{align*}
  Y_k(t,\cdot)\,\bnu:=\trrdt \bY_k(t,\cdot)=d\Psi_{\R\delta(t)}\,
  (\det d\Psi_{\R\delta(t)})^{-1}\,\widehat Y_k\,\bnu
\end{align*}
makes sense. The identity $\trnormald \bY_k=Y_k\,\gamma(\R\delta)$ and
Proposition \ref{lemma:partInt} imply that
\begin{align*}
  \int_M Y_k(t,\cdot)\,\gamma(\R\delta(t,\cdot))\ dA = 0.
\end{align*}
Since the fields $\widehat Y_k$ are a basis we see that the functions
$Y_k(t,\cdot)$ form a basis of the space 
\begin{equation}\label{eqn:raum}
  \Big\{Y\in H^2_0(M)\ \big|\ \int_M
  Y\, \gamma(\R\delta(t,\cdot))\ dA = 0\Big\}. 
\end{equation}
To simplify the notation we choose an enumeration\footnote{For instance, we can
choose $\bW_{2k}:=\bX_k$ and $\bW_{2k-1}:=\bY_k$.}
$(\bW_k)_{k\in\setN}$ of  $\bX_k,\bY_k$ and set $W_k\,\bnu:=\trrd
\bW_k$. 

Now  we show that 
\begin{align*}
  \spann\{(\phi\, W_k,\phi\,\bW_k)\ |\ \phi\in C_0^1([0,T)),\,
  k\in\setN\}
\end{align*}
is dense in the space of all couples 
\begin{align*}
  (b,\bphi)\in \big(H^1(I,L^2(M))\cap L^2(I,H^2_0(M))\big)\times
  H^1(\Omega_{\R\delta}^I)
\end{align*}
with $b(T,\cdot)=0$, $\bphi(T,\cdot)=0$, $\dv\bphi=0$, and
$\trrd\bphi=b\,\bnu$. Clearly, $T_{\R\delta}^I$ embeds into this
space. Since $\T_{\R\delta}$ induces isomorphisms between corresponding function spaces the assertion is equivalent to the density of 
\begin{align*}
  \spann\{(\phi\,\widehat W_k,\phi\,\widehat\bW_k)\ |\ \phi\in
  C_0^1([0,T)),\, k\in\setN\}
\end{align*}
in the space $T$ of all couples 
\begin{align*}
  (b,\bphi)\in H^1(I,L^2(M))\cap L^2(I,H^2_0(M))\times
  H^1(I\times\Omega)
\end{align*}
with $b(T,\cdot)=0$, $\bphi(T,\cdot)=\boldsymbol{0}$, $\dv\bphi=0$, and
$\bphi|_{I\times\pa\Omega} = b\,\bnu$. Let $(b,\bphi)\in T$. We
approximate $b$ by functions $\tilde b\in C^\infty_0([0,T),H^2_0(M))$
in $H^1(I,L^2(M))\cap L^2(I,H^2_0(M))$ with
\begin{align*}
  \int_M \tilde b(t,\cdot)\ dA=0
\end{align*}
for $t\in \bar I$.\footnote{We can construct $\tilde b$ in the
  following way. We first extend $b$ by $0$ to $[0,\infty)$ and then
  convolute $b(\cdot+h)$, $h>0$, with a standard mollification kernel.
  Note that the mean value is preserved.}  Now, we can approximate
$\pa_t\tilde b$ by a sequence $(\phi^k_n\, \widehat Y_k)_n$ (summation
with respect to $k$ from $1$ to $n$) with $\phi^k_n\in C_0^1([0,T))$,
in $L^2(I,H^2_0(M))$. Indeed, for $f\in L^2(I,H^2_0(M))$ with
\begin{align*}
  \int_M f(t,\cdot)\ dA=0
\end{align*}
for almost  all $t\in I$ and
\begin{align*}
  \int_I\phi(t)\, (\widehat Y_k,f(t,\cdot))_{H^2(M)}\ dt=0
\end{align*}
for all $\phi\in C_0^1([0,T))$ and all $k\in\setN$, we see that the 
scalar product in the integrand vanishes almost everywhere. Thus we
get $f=0$, since the functions $\widehat Y_k(t,\cdot)$ form a basis. From 
\begin{equation*}
  \Bignorm{\tilde b(t,\cdot) + \int_t^T\phi^k_n(s)\, ds\,
    \widehat Y_k}_{H^2(M)}\le \int_0^T \norm{\pa_s\tilde b(s,\cdot) 
    - \phi^k_n(s)    \,\widehat Y_k}_{H^2(M)}\ ds
\end{equation*}
we infer that the sequence $(-\int_t^T\phi^k_n(s)\, ds\ \widehat
Y_k)_{n}$ converges to $\tilde b$ in $C(\bar I, H^2_0(M))$. Thus, these linear combinations of the fields $\widehat Y_k$ converge to
$b$ in $H^1(I,L^2(M))\cap L^2(I,H^2_0(M))$. In view of the continuity
properties of the solution operator of the Stokes system we deduce that
the corresponding linear combinations of the vector fields
$\widehat\bY_k$ converge to some function $\bY$ in $H^1(I\times\Omega)$.  Moreover,
we have $(\bY-\bphi)|_{I\times\pa\Omega}=\bfzero$, so that it remains
to show that we can approximate each $\bphi$ with $(0,\bphi)\in
T$ by a sequence  $(\phi^k_n\, \widehat \bX_k)_{n}$, $\phi^k_n\in
C_0^1([0,T))$, in $H^1(I\times \Omega)$. This, however, can be done analogously to the  approximation of $b$.\\

\noindent {\bf Proof} (of Proposition \ref{lemma:ent}){\bf :} We use
the Galerkin method. We seek
functions $\alpha_n^k:[0,T]\rightarrow\setR$, $k,n\in\setN$, such that
$\bu_n:=\alpha_n^k\,\bW_k$ and $\eta_n(t,\cdot):=\int_0^t\alpha^k_n\,
W_k\ ds+\eta_0$ (summation with respect to $k$ from $1$ to $n$) solve
the equation
\begin{equation}\label{eqn:gal}
  \begin{aligned}
    &\iortd \pa_t\bu_n\cdot\bW_j\ dx +\frac{1}{2} \im\pa_t\eta_n
    \, \pa_t\mathcal{R}\delta\, W_j\, \gamma(\mathcal{R}\delta)\ dA\\
    &+\frac{1}{2}\iortd(\mathcal{R}\bv\cdot\nabla)\bu_n\cdot\bW_j\ dx
    -
    \frac{1}{2}\iortd(\mathcal{R}\bv\cdot\nabla)\bW_j\cdot\bu_n\ dx\\
    &+\iortd \nabla\bu_n:\nabla\bW_j\ dx +\im\pa^2_t\eta_n\, W_j\ dA +
    2K(\eta_n, W_j) \\ 
    &\hspace{4cm}=\iortd \ff_n\cdot\bW_j\ dx + \im g_n\, W_j\ dA
 \end{aligned}
\end{equation}
for all $1\le j \le n$. Here $\ff_n$ and $g_n$ a smooth functions
which converge\footnote{We can take, e.g., $\ff_n:=\R_{1/n}\ff$ and
  $g_n:=\sum_\alpha(\omega_{1/n}\ast((\psi_\alpha
  g)\circ\phi^{-1}_\alpha))\circ\phi_\alpha$ (cf. Definition
  \eqref{eqn:reg}).} to $\ff$ and $g$ in $L^2_\loc([0,\infty)\times
\setR^3)$ and $L^2_\loc([0,\infty)\times M)$, respectively. 
We choose initial conditions $\alpha^k_n(0)$ such that
\begin{align*}
  \pa_t\eta_n(0,\cdot)&\rightarrow\eta_1\qquad \text{ in } L^2(M),
  \\
  \bu_n(0,\cdot)&\rightarrow \bu_0\qquad \text{ in }L^2(\Omega_{\R\eta_0}).
\end{align*}
To this end, we choose the coefficients of $\bY_k$ at $t=0$ such that the
first convergence holds. This is possible since the functions $Y_k(0,\cdot)$ form a
basis of the space \eqref{eqn:raum} with $t=0$ and since Proposition
\ref{lemma:nspur} yields the identity 
\begin{align*}
  \int_M \eta_1\, \gamma(\R\eta_0)\ dA=\int_{\Omega_{\R\eta_0}}\dv\bu_0\ dx=0.
\end{align*}
By Theorem \ref{theorem:stokes} the solution operator of
the Stokes system maps 
\begin{align*}
  \Big\{Y\,\bnu\in L^2(M)\ |\ \int_M Y\,\gamma(\R\eta_0)\ dA=0\Big\}
\end{align*}
continuously into $L^2(\Omega_{\R\eta_0})$. Thus, from the convergence of the linear
combinations of the functions $Y_k(0,\cdot)$ to $\eta_1$ in $L^2(M)$ we deduce the convergence of the corresponding linear combinations of the fields $\bY_k(0,\cdot)$ to some $\bY$ in $L^2(\Omega_{\R\eta_0})$. We also have $\trnormaleo(\bu_0-\bY) = 0$. Since the fields
$\bX_k(0,\cdot)$ form a basis of 
\begin{align*}
  \big\{\bX\in L^2(\Omega_{\R\eta_0})\ |\ \dv \bX=0,\, \trnormaleo
  \bX=0\big\}, 
\end{align*}
we can choose their coefficients at $t=0$ such that the sequence
$(\alpha_n^k\,\bW_k)_{n}$ converges to $\bu_0$ in $L^2(\Omega_{\R\eta_0})$.

Thus, \eqref{eqn:gal} with the above initial conditions is a
Cauchy problem for a linear system of ordinary
integro-differential equations of the form ($1\le j\le n$, summation
with respect to $k$ from $1$ to $n$)
\begin{equation*}
  A_{jk}(t)\,\dot\alpha^k(t) = B_{jk}\,\alpha^k(t) + \int_0^t
  C_{jk}(t,s)\,\alpha^k(s)\ ds +   D_j(t).
\end{equation*}
The coefficients, given through
\begin{equation*}
  \begin{aligned}
    A_{jk}(t)&=\iortd \bW_k\cdot\bW_j\ dx + \im W_k\, W_j\ dA,\\
    B_{jk}(t)&=\iortd \pa_t\bW_k\cdot\bW_j\ dx + \frac{1}{2}\im W_k\,
    W_j\
    \pa_t\R\delta\, \gamma(\R\delta)\ dA\\
    &\hspace{0.3cm} +\frac{1}{2} \iortd
    (\R\bv\cdot\nabla)\bW_k\cdot\bW_j\ dx - \frac{1}{2} \iortd
    (\R\bv\cdot\nabla)\bW_j\cdot\bW_k\ dx \\
    &\hspace{0.3cm} + \iortd \nabla\bW_k:\nabla\bW_j\ dx + \im \pa_t W_k\, W_j\ dA,\\
    C_{jk}(t,s)&=2 K(W_k(s),W_j(t)),\\
    D_j(t)&=\iortd \ff_n\cdot \bW_j\ dx + \im g_n\, W_j\ dA,
  \end{aligned}
\end{equation*}
are continuous. The tensor $A(t)$ is symmetric and positive
definite since the basis functions are linearly independent and we have
($\beta^j\in\setR$, $j=1,\ldots,n$)
\begin{align*}
  \beta^j\beta^k A_{jk}(t)=\iortd |\beta^j\bW_j|^2\ dx+\im |\beta^j W_j|^2\ dA.
\end{align*}
In particular, $A(t)$ is invertible and we infer from Proposition
\ref{pro:ode} that \eqref{eqn:gal} is solvable on the interval $[0,T]$
for all $n$. In order to derive energy estimates we test
\eqref{eqn:gal} with $(\pa_t\eta_n,\bu_n)$ and obtain
\begin{align*}
  &\iortd \pa_t\bu_n\cdot\bu_n\ dx +\frac{1}{2} \im(\pa_t\eta_n)^2
  \, \pa_t\mathcal{R}\delta\, \gamma(\mathcal{R}\delta)\ dA\\
  &+\iortd |\nabla\bu_n|^2\ dx
  +\im\pa^2_t\eta_n\, \pa_t\eta_n\ dA + 2K(\eta_n,\pa_t\eta_n)\\
  &\hspace{4cm}=\iortd \ff_n\cdot\bu_n\ dx + \im g_n\, \pa_t\eta_n\
  dA.
\end{align*}
By Reynold's transport theorem we obtain
\begin{align*}
  \frac{1}{2}\frac{d}{dt}&\iortd |\bu_n|^2\ dx + \iortd |\nabla\bu_n|^2\ dx
  +\frac{1}{2}\frac{d}{dt}\im|\pa_t\eta_n|^2\ dA + \frac{d}{dt}
  K(\eta_n)
  \\
  &\hspace{5cm}=\iortd \ff_n\cdot\bu_n\ dx + \im g_n\, \pa_t\eta_n\ dA.
\end{align*}
Now, we can proceed as in the end of Subsection \ref{subsec:apriori} to
obtain
\begin{align*}
  \norm{\eta_n}_{Y^I}^2 + \norm{\bu_n}_{X_{\R\delta}^I}^2\le
  c_0(T,\Omega_{\R\delta}^I,\ff_n,g_n,\bu_n(0,\cdot),\eta_0,\pa_t\eta_n(0,\cdot)).
\end{align*}
This in turn implies\footnote{Note that 
  \begin{equation*}
    \begin{aligned}
      X_{\R\delta}^I&\simeq L^\infty(I,L^2(\Omega))\cap
      L^2(I,H^1_{\dv}(\Omega))\\
      &\simeq (L^1(I,L^2(\Omega))+ L^2(I,H^1_{\dv}(\Omega)))'\\
      &\simeq (L^1(I,L^2(\Omega_{\R\delta(t)}))+
      L^2(I,H^1_{\dv}(\Omega_{\R\delta(t)})))',
    \end{aligned}
  \end{equation*}
  where the first and the third identification are induced by the
  mapping $\T_{\R\delta}$. }
\begin{equation*}
  \begin{aligned}
    \eta_n&\overset{*}{\weakto} \eta \hspace{0.3cm}&&\text{
      weakly$^*$ in } L^\infty(I,H^2_0(M)),\\
    \pa_t\eta_n&\overset{*}{\weakto} \pa_t\eta &&\text{ weakly$^*$ in
    } L^\infty(I,L^2(M)),\\ 
    \bu_n&\overset{*}{\weakto} \bu &&\text{ weakly$^*$ in }
    X_{\R\delta}^I.
  \end{aligned}
\end{equation*}
The lower semi-continuity of the norm with respect to the weak*
convergence shows \eqref{ab:ent}. Moreover, the above convergences and
$\trrd \bu_n=\pa_t\eta_n\,\bnu$ imply the identity $\trrd
\bu=\pa_t\eta\,\bnu$.

In order to show that $(\eta, \bu)$ satisfies \eqref{eqn:ent} we multiply
\eqref{eqn:gal} by $\phi(t)$, where $\phi\in C_0^1([0,T))$,
integrate over $I$ and subsequently use integration by parts with
respect to time. This yields for $1\le j\le n$
\begin{equation*}
  \begin{aligned}
    &-\int_I\iortd \bu_n\cdot\pa_t(\phi\,\bW_j)\ dxdt
    -\frac{1}{2}\int_I\im\pa_t\eta_n \, \pa_t\mathcal{R}\delta\,
    \phi\, W_j\, \gamma(\mathcal{R}\delta)\
    dAdt\\
    & +\frac{1}{2}\int_I\iortd(\mathcal{R}\bv\cdot\nabla)\bu_n\cdot(\phi\,\bW_j)\ dxdt\\
    &-\frac{1}{2}\int_I\iortd(\mathcal{R}\bv\cdot\nabla)(\phi\,\bW_j)\cdot\bu_n\
    dxdt
    +\int_I\iortd \nabla\bu_n:\nabla(\phi\,\bW_j)\ dxdt\\
    &+\int_I\im\pa_t\eta_n\, \pa_t (\phi\, W_j)\ dAdt + 2\int_I
    K(\eta_n, \phi\, W_j)\
    dAdt \\
    &\hspace{1cm}=\int_I\iortd \ff_n\cdot(\phi\,\bW_j)\ dxdt +
    \int_I\im g_n\, \phi\, W_j\
    dAdt\\
    &\hspace{1.5cm}
    +\int_{\Omega_{\mathcal{R}\eta_0}}\bu_n(0)\cdot(\phi(0)\bW_j(0,\cdot))\
    dx + \im\pa_t\eta_n(0)\, \phi(0)W_j(0,\cdot)\ dA.
  \end{aligned}
\end{equation*}
The limiting process for $n\to \infty$ in this identity shows that
$\eta$ and $\bu$ satisfy \eqref{eqn:ent} for all test functions from
\begin{align*}
  \spann\{(\phi\, W_j,\phi\,\bW_j)\ |\ \phi\in C_0^1([0,T)),\,
  j\in\setN\}.
\end{align*}
By density we get the validity of \eqref{eqn:ent}  for all test
functions from $T^I_{\R\delta}$. Note that due to the regularisation
and linearisation of the convective terms the convergence of the test
functions in $H^1(\Omega_{\R\delta}^I)$ is sufficient.
\qed
\smallskip
\begin{Remark}
  Note that the uniqueness of  weak solutions we have just constructed
  is not obvious, due to the mixed character of the system and the
  non-cylindrical spacetime domain. We circumvent this difficulty by
  using a multi-valued version of Schauder's fixed point theorem,
  namely the Kakutani-Glicksberg-Fan Theorem \ref{theorem:kakutani}.
\end{Remark}

\subsection{Fixed-point argument}
 
Let us now define solutions of our regularized problem. 

\begin{Definition} \label{def:epsilon} A couple $(\eta,\bu)$ is a
  weak solution of the regularized problem with the parameter
  $\epsilon$ in the interval $I$ if $\eta=\eta_\epsilon\in Y^I$ with
  $\norm{\eta}_{L^\infty(I\times M)}<\kappa$, $\eta(0,\cdot)=\eta_0$,
  and $\bu=\bu_\epsilon\in X_{\mathcal{R}_\epsilon\eta}^I$ with $\trre
  \bu = \pa_t\eta\,\bnu$ satisfy
  \begin{align}
    &- \int_I\iorte \bu\cdot\pa_t\bphi\ dxdt -
    \frac{1}{2}\int_I\im\pa_t\eta \, \pa_t\mathcal{R}_\epsilon\eta\
    b\,
    \gamma(\mathcal{R}_\epsilon\eta)\ dAdt\notag\\
    &+\frac{1}{2}\int_I\iorte(\mathcal{R}_\epsilon\bu\cdot\nabla)
    \bu\cdot\bphi\ dxdt -
    \frac{1}{2}\int_I\iorte(\mathcal{R}_\epsilon\bu\cdot\nabla)
    \bphi\cdot\bu\ dxdt\label{eq:ent-a}\\
    & + \int_I\iorte \nabla\bu:\nabla\bphi\ dxdt-\int_I\im\pa_t\eta\,
    \pa_tb\ dAdt + 2\int_I K(\eta,b)\
    dt\notag \\
    &\hspace{3.5cm}=\int_I\iorte \ff\cdot\bphi\ dxdt + \int_I\im g\,
    b\ dAd\notag\\
    &\hspace{4cm}
    +\int_{\Omega_{\R_\epsilon\eta_0}}\bu_0^\epsilon\cdot\bphi(0,\cdot)\
    dx + \im\eta_1^\epsilon\, b(0,\cdot)\ dA \notag 
  \end{align}
  for all test functions $(b,\bphi)\in  T_{\R_\epsilon\eta}^I$.
\end{Definition}
\smallskip

\begin{Proposition}\label{lemma:approx}
  There exists a $T>0$ such that for all sufficiently small $\epsilon>0$
  there exists a weak solution $(\eta_\epsilon,\bu_\epsilon)$ of the regularized system with the
  parameter $\epsilon$ in the interval $I=(0,T)$. Furthermore, we have
  \begin{equation}\label{ab:approx}
    \begin{aligned}
      &\norm{\eta_\epsilon}_{Y^I}^2 + \norm{\bu_\epsilon}_{X_{\R_\epsilon\eta}^I}^2\le
      c_0(T,\Omega_{\R_\epsilon\eta}^I,\ff,g,\bu^\epsilon_0,\eta_0, 
      \eta^\epsilon_1) 
    \end{aligned}
  \end{equation}
  and $\sup_{\epsilon}\tau(\eta_\epsilon)<\infty$. The time $T$ can be
  chosen to depend only on $\tau(\eta_0)$ and the bound in
  \eqref{ab:approx} for the $Y^I$ norm of $\eta_\epsilon$.
\end{Proposition}
\beweis We set $\alpha:=(\norm{\eta_0}_{L^\infty(M)}+\kappa)/2$ and
fix an arbitrary but sufficiently small  $\epsilon>0$. For a better
readability we will not write this parameter $\epsilon >0$ explicitly.
We want to use Theorem \ref{theorem:kakutani}. To this end, we define 
the space 
\begin{align*}
  Z:=C(\bar I\times\pa\Omega)\times L^2(I\times \setR^3)
\end{align*}
and the convex set 
\begin{align*}
  D:=\big\{(\delta,\bv)\in Z\ |\ \delta(0,\cdot)=\eta_0,\,
  \norm{\delta}_{L^\infty(I\times\pa\Omega)}\le\alpha,\,
  \norm{\bv}_{L^2(I\times \setR^3)}\le c_1\big\},
\end{align*}
where $c_1>0$ is sufficiently large. Let 
\begin{align*}
  F:D\subset Z\rightarrow 2^Z
\end{align*}
be the mapping which assigns to each couple $(\delta,\bv)$ the set of
all weak solutions $(\eta,\bu)$ of the decoupled and regularized
problem with  data $(\delta,\bv)$ that satisfy the estimate 
\begin{equation}\label{ab:schranke}
  \norm{\eta}_{Y^I} + \norm{\bu}_{X_{\R\delta}^I} \le
  c_0(T,\Omega_{\R\delta}^I,\ff,g,\bu_0,\eta_0,\eta_1). 
\end{equation}
Proposition \ref{lemma:ent} implies that
$F(\delta,\bv)$ is non-empty. Let $(\eta,\bu)\in F(\delta,\bv)$. From
\eqref{ab:schranke} we deduce that $\norm{\bu}_{L^2(I\times \setR^3)} \le
c_1$ and that the norm of $\eta$ in
\begin{align}\label{eq:emb}
  Y^I\embedding C^{0,1-\theta}(\bar I, C^{0,2\theta
    -1}(\pa\Omega))\quad (1/2<\theta<1)
\end{align}
is bounded. Since $\eta(0,\cdot)=\eta_0$ we can choose the time
interval $I=(0,T)$ so small that $\norm{\eta}_{L^\infty(I\times M)}\le\alpha$, independently of the parameter $\epsilon$; in particular, $\tau(\eta)\le c(\alpha)$. Thus,
$F$ maps the set $D$ into its power set, i.e.~$F(D)\subset 2^D$. Since
the problem is linear the set $F(\delta,\bv)$ is convex. Moreover, the
set $F(\delta,\bv)$ is closed in $Z$. Indeed, let
$(\eta_n,\bu_n)\subset F(\delta,\bv)$ be a sequence which converges to
some  $(\eta,\bu)$ in $Z$. The estimate \eqref{ab:schranke} implies
that for a subsequence we have 
\begin{equation*}
 \begin{aligned}
   \eta_n&\weakastto \eta \hspace{0.3cm}&&\text{ weakly$^*$ in }
   L^\infty(I,H^2_0(M)),\\ 
   \pa_t\eta_n&\weakastto \pa_t\eta &&\text{ weakly$^*$ in }
   L^\infty(I,L^2(M)),\\ 
   \bu_n&\weakastto\bu &&\text{ weakly$^*$ in } X_{\R\delta}^I. 
 \end{aligned}
\end{equation*}
Thus, we can pass to the limit in \eqref{eqn:ent} and we obtain that $(\eta,\bu)\in F(\delta,\bv)$.

Next, we want to show that $F(D)$ is relatively compact in $Z$. For this purpose, let $(\delta_n,\bv_n)\subset D$ be a sequence and $(\eta_n,\bu_n)\in F(\delta_n,\bv_n)$. We have to show that there exists a subsequence of $(\eta_n,\bu_n)$ which converges in $Z$. In
view of the construction of $F$ and of \eqref{eq:emb} we immediately
get the uniform convergence of a subsequence of $(\eta_n)$. The
mollification operator 
\begin{align*}
  \mathcal{R}:\{\delta\in C(\bar I\times\pa\Omega)|\
  \delta(0,\cdot)=\eta_0 \}\rightarrow C^3(\bar
  I\times\pa\Omega)\compactembedding C^2(\bar I\times\pa\Omega)
\end{align*}
gives compactness and thus there exists a subsequence of
$(\mathcal{R}\delta_n)$ converging to $\mathcal{R}\delta$ in $C^2(\bar
I\times\pa\Omega)$. The proof of the relative compactness of $(\bu_n)$
in $L^2(I\times\setR^3)$ can be taken almost literally from the proof of Proposition \ref{lemma:komp}. The only differences are the slightly changed form \eqref{eqn:ent} of the system
and that in some places one has to replace the sequence $(\eta_n)$ with
limit $\eta$ by the sequence $(\R\delta_n)$ with limit $\R\delta$.
Moreover, due to the regularization one can simplify some of the
arguments. This reasoning also yields the relative compactness of
$(\pa_t\eta_n)$ in ${L^2(I\times M)}$.

It remains to show that the mapping $F$ has a closed graph.\footnote{We
  prove this in such a way that the arguments also hold under weaker
  assumptions on the regularity of the boundary, since we will need
  these arguments in the proof of \ref{theorem:hs} again. } 
Let $(\delta_n,\bv_n)\subset D$ and $(\eta_n,\bu_n)\in F(\delta_n,\bv_n)$ be sequences with $(\delta_n,\bv_n)\rightarrow(\delta,\bv)$ and $(\eta_n,\bu_n)\rightarrow(\eta,\bu)$ in $Z$. We have to show that $(\eta,\bu)\in F(\delta,\bv)$. From \eqref{ab:schranke} and the relative compactness we just showed we can deduce that for a subsequence
we have 
\begin{equation}\label{eqn:wichtigekonv}
  \begin{aligned}
    \eta_n&\rightarrow\eta &&\text{ uniformly and weakly$^*$ in
    }L^\infty(I,H^2_0(M)),\\ 
    \pa_t\eta_n&\rightarrow\pa_t\eta &&\text{ strongly in
    }L^2(I\times M)\text{ and weakly$^*$ in }L^\infty(I,L^2(M)),\\
    \bu_n&\rightarrow\bu &&\text{ strongly in
    }L^2(I\times \setR^3)\text{ and weakly$^*$ in }
    L^\infty(I,L^2(\setR^3)),\\ 
    \nabla\bu_n&\weakto\nabla\bu &&\text{ weakly in } L^2(I\times
    \setR^3).
  \end{aligned}
\end{equation}
We extend $\nabla\bu_n$ and $\nabla\bu$, which are a-priori defined on
$\Omega_{\R\delta_n}^I$ and $\Omega_{\R\delta}^I$, resp., by $\bfzero$
to the whole of $I\times \setR^3$. The lower semi-continuity of the
norms implies that also $\eta$ and $\bu$ satisfy the estimate
\eqref{ab:schranke}. The property $\eta(0,\cdot)=\eta_0$ follows
immediately from the uniform convergence of $(\eta_n)$. Moreover, we
have 
\begin{equation*}
  \begin{aligned}
    \pa_t\eta_n\,\bnu=\trrdn \bu_n = \bw_n|_{\pa\Omega},
 \end{aligned}
\end{equation*}
where $\bw_n:=\bu_n\circ\Psi_{\R\delta_n}$. By Lemma
\ref{lemma:psi} there exists a subsequence of $(\bw_n)$ which converges to some $\bw$ weakly in $L^2(I,W^{1,3/2}(\Omega))$. This implies
$\pa_t\eta\,\bnu=\bw|_{\pa\Omega}$. 
In the estimate 
\begin{align*}
  &\norm{\bw_n-\bu\circ\Psi_{\R\delta}}_{L^1(I\times \Omega)}
  \\
  &\hspace{0.5cm}\le
  \norm{(\bu_n-\bu)\circ\Psi_{\R\delta_n}}_{L^1(I\times \Omega)}+
  \norm{\bu\circ\Psi_{\R\delta_n}-\bu\circ\Psi_{\R\delta}}_{L^1(I\times
    \Omega)},
\end{align*}
the first term tends to zero due to Lemma \ref{lemma:psi} and
\eqref{eqn:wichtigekonv}$_3$ while the second term converges to zero
due to Remark \ref{bem:konv}. Thus we have $\bw=\bu\circ\Psi_{\R\delta}$ and, hence, 
$\pa_t\eta\,\bnu=\trrd\bu$. It remains to show that
\eqref{eqn:ent} is satisfied.  We have for all $n$ and all test
functions $(b_n,\bphi_n)\in T_{\R\delta_n}^I$
\begin{align}
  &- \int_I\iortdn \bu_n\cdot\pa_t\bphi_n\ dxdt -
  \frac{1}{2}\int_I\im\pa_t\eta_n
  \, \pa_t\mathcal{R}\delta_n\, b_n\, \gamma(\mathcal{R}\delta_n)\
  dAdt\label{eqn:graphab} \\
  &+\frac{1}{2}\int_I\iortdn(\mathcal{R}\bv_n\cdot\nabla)\bu_n\cdot\bphi_n\
  dxdt - \frac{1}{2}\int_I\iortdn(\mathcal{R}\bv_n\cdot\nabla)
  \bphi_n\cdot\bu_n\ dxdt \notag \\
  & + \int_I\iortdn \nabla\bu_n:\nabla\bphi_n\
  dxdt-\int_I\im\pa_t\eta_n\, \pa_tb_n\ dAdt + 2\int_I K(\eta_n,b_n)\
  dt \notag\\
  &\hspace{0.5cm}=\int_I\iortdn \ff\cdot\bphi_n\ dxdt + \int_I\im g\,
  b_n\    dAdt\notag\\
  &\hspace{1cm}
  +\int_{\Omega_{\mathcal{R}\eta_0}}\mathcal{R}\bu_0\cdot\bphi_n(0,\cdot)\
  dx + \im \R\eta_1\, b_n(0,\cdot)\ dA.\notag
  \end{align}
  The limiting procedure is not immediate since the test functions
  depend on $n$. Hence, for arbitrary $(b,\bphi)\in
  T_{\R\delta}^I$ we choose the special test functions
  $(b_n,\bphi_n):=(\M_{\R\delta_n}b,\F_{\R\delta_n}\M_{\R\delta_n}b)\in
  T_{\R\delta_n}^I$ which already have been used in the proof of
  Proposition \ref{lemma:komp}. We choose the number $\alpha$ in Propositon
  \ref{lemma:FortVonRandZeit} to be $(\alpha+\kappa)/2$. In view of the assertions $\rm (1.b)$, $\rm (2.b)$ in Lemma
  \ref{lemma:konvergenzen} and the convergences in
  \eqref{eqn:wichtigekonv} we can pass to the limit in
  \eqref{eqn:graphab} for the above special test functions. This
  yields the validity of \eqref{eqn:ent} for
  $(b,\bphi)=(b,\F_{\R\delta}b)\in T_{\R\delta}^I$. Due to the
  definition of $T_{\R\delta}^I$ it remains to show the validity of
  \eqref{eqn:ent} for test functions $(0,\bphi)\in T_{\R\delta}^I$
  with $\bphi(T,\cdot)=0$ and
  $\supp\bphi\subset\Omega_{\R\delta}^{\bar I}$. From the uniform
  convergence of $(\R\delta_n)$ it follows for sufficient large $n$ that 
  $(0,\bphi)\in T_{\R\delta_n}^I$. This in turn enables the limiting
  process in  \eqref{eqn:graphab}. 

  Now, Theorem \ref{theorem:kakutani} guarantees the existence of a 
  fixed point of $F$, i.e.~there exists a couple $(\eta,\bu)\in D$
  with $(\eta,\bu)\in F(\eta,\bu)$. This concludes the proof.\qed

\smallskip

\subsection{Limiting process}
Now, we can prove our main result by letting the
regularizing parameter $\epsilon$ in Definition \ref{def:epsilon} tend
to zero.\\

\beweis (of Theorem \ref{theorem:hs}) We have shown that there exists
a $T>0$ such that for all $\epsilon=1/n$, $n\in\setN$ sufficiently
large, there exists a weak solution $(\eta_\epsilon,\bu_\epsilon)$ of
the regularized problem with the parameter $\epsilon$ in the interval
$I=(0,T)$. The estimate \eqref{ab:approx} and the compact embedding
$Y^I\compactembedding C(\bar I\times\pa\Omega)$ yields the following
convergences for a subsequence 
\begin{equation}\label{eqn:konv3}
  \begin{aligned}
    \eta_\epsilon,\,\R_\epsilon\eta_\epsilon&\rightarrow\eta &&\text{
      uniformly and weakly$^*$ in }  L^\infty(I,H^2_0(M)),\\
    \pa_t\eta_\epsilon,\,
    \pa_t\R_\epsilon\eta_\epsilon&\weakastto\pa_t\eta &&\text{
      weakly$^*$ in }
    L^\infty(I,L^2(M)),\\
    \bu_\epsilon&\weakastto\bu &&\text{ weakly$^*$ in }
    L^\infty(I,L^2(\setR^3)),\\
    \nabla\bu_\epsilon&\weakto\nabla\bu &&\text{ weakly in
    }L^2(I\times \setR^3).
\end{aligned}
\end{equation}
We extend $\nabla\bu_n$ and $\nabla\bu$, which are a-priori defined on
$\Omega_{\R\delta_n}^I$ and $\Omega_{\R\delta}^I$, resp., by $\bfzero$
to the whole of $I\times \setR^3$.  The uniform convergence of
$(\R_\epsilon\eta_\epsilon)$ follows from the estimate 
\begin{align*}
  \norm{ \R_\epsilon\eta_\epsilon-\eta}_{L^\infty(I\times \pa\Omega)}
  \le\norm{\R_\epsilon(\eta_\epsilon-\eta)}_{L^\infty(I\times
    \pa\Omega)}+\norm{\R_\epsilon\eta-\eta}_{L^\infty(I\times
    \pa\Omega)}.
\end{align*}
Now, we can repeat the
proof of  Proposition \ref{lemma:komp} almost literally to show that
\begin{equation}\label{eqn:konv31}
  \begin{aligned}
    \pa_t\eta_\epsilon&\rightarrow\pa_t\eta &&\text{ strongly in }
    L^2(I\times M),\\
    \bu_\epsilon&\rightarrow\bu &&\text{ strongly in }
    L^2(I\times\setR^3).
  \end{aligned}
\end{equation}
As in the proof of Proposition \ref{lemma:approx} we obtain the
identity $\tr\bu=\pa_t\eta\,\bnu$. Moreover, using \eqref{eqn:konv31},
the interpolation inequality 
\begin{align*}
  \norm{\bu_\epsilon-\bu}_{L^2(I,L^4(\setR^3))}\le
  \norm{\bu_\epsilon-\bu}_{L^2(I\times
    \setR^3)}^{1/6}\,\norm{\bu_\epsilon-\bu}_{L^2(I,L^5(\setR^3))}^{5/6}
\end{align*}
and Corollary \ref{lemma:sobolev} we get 
\begin{equation}\label{eqn:konv20}
  \begin{aligned}
    &\bu_\epsilon\rightarrow\bu&&\text{ in }L^2(I,L^4(\setR^3)).
  \end{aligned}
\end{equation}
Consequently, we have
\begin{equation}\label{eqn:konv8}
  \begin{aligned}
    \pa_t\R_\epsilon\eta_\epsilon&\rightarrow\pa_t\eta &&\text{ in
    }L^2(I\times \pa\Omega),\\ 
    \R_\epsilon\bu_\epsilon&\rightarrow\bu &&\text{ in
    }L^2(I,L^4(\setR^3)).
  \end{aligned}
\end{equation}
The lower semi-continuity of the norms yields the estimate
\eqref{ab:hs}, while the uniform convergence of $(\eta_\epsilon)$
gives $\eta(0,\cdot)=\eta_0$. For all $\epsilon$ and all
$(b_\epsilon,\bphi_\epsilon)\in T_{\R_\epsilon\eta_\epsilon}^I$ we
have
\begin{align}
  &- \int_I\int_{\Omega_{R_\epsilon\eta_\epsilon(t)}}
  \bu_\epsilon\cdot\pa_t\bphi_\epsilon\ dxdt -
  \frac{1}{2}\int_I\im\pa_t\eta_\epsilon\,
  \pa_t\mathcal{R}_\epsilon\eta_\epsilon\, b_\epsilon\,
  \gamma(\mathcal{R}_\epsilon\eta_\epsilon)\
  dAdt\notag \\
  &+\frac{1}{2}\int_I\int_{\Omega_{R_\epsilon\eta_\epsilon(t)}}(\mathcal{R}
  _\epsilon\bu_\epsilon\cdot\nabla)\bu_\epsilon\cdot\bphi_\epsilon\
  dxdt -
  \frac{1}{2}\int_I\int_{\Omega_{R_\epsilon\eta_\epsilon(t)}}(\mathcal{R}
  _\epsilon\bu_\epsilon\cdot\nabla)\bphi_\epsilon\cdot\bu_\epsilon\
  dxdt\notag \\ 
  & + \int_I\int_{\Omega_{R_\epsilon\eta_\epsilon(t)}}
  \nabla\bu_\epsilon:\nabla\bphi_\epsilon\
  dxdt-\int_I\im\pa_t\eta_\epsilon\, \pa_tb_\epsilon\ dAdt + 2\int_I
  K(\eta_\epsilon,b_\epsilon)\
  dt \notag \\
  &\hspace{3.5cm}=\int_I\int_{\Omega_{R_\epsilon\eta_\epsilon(t)}}
  \ff\cdot\bphi_\epsilon\ dxdt + \int_I\im g\, b_\epsilon\
  dAdt\label{eqn:sepp}\\
  &\hspace{4cm}
  +\int_{\Omega_{\eta_0^\epsilon}}\bu_0^\epsilon\cdot\bphi_\epsilon(0,\cdot)\
  dx + \im\eta_1^\epsilon\, b_\epsilon(0,\cdot)\ dA.\notag 
  \end{align}
  Just like in the proof of Proposition \ref{lemma:approx}, we make use of the special test functions
  $(b_\epsilon,\bphi_\epsilon):=(\M_{\R_\epsilon\eta_\epsilon}b,
  \F_{\R_\epsilon\eta_\epsilon}\M_{ \R_\epsilon\eta_\epsilon}b)\in
  T_{\R_\epsilon\eta_\epsilon}^I$ for $(b,\bphi)\in T_\eta^I$. Here, we assume that $\alpha$ in
  Propositon \ref{lemma:FortVonRandZeit} satisfies $\sup_{\epsilon}\norm{\eta_\epsilon}_{L^\infty(I\times
    M)}<\alpha<\kappa$. Since the sequence
  $(\nabla\,\F_{\R_\epsilon\eta_\epsilon}\M_{\R_\epsilon\eta_\epsilon}b)$
  is bounded in $L^\infty(I,L^{2}(B_\alpha))$, we get
  \begin{equation}\label{eqn:konv5}
    \begin{aligned}
      \nabla\,\F_{\R_\epsilon\eta_\epsilon}\M_{\R_\epsilon\eta_\epsilon}
      b\weakto\nabla\,\F_{\eta}\M_{\eta}b\quad\text{ weakly$^*$
        in }L^\infty(I,L^2(B_\alpha)).
    \end{aligned}
  \end{equation}
  Using \eqref{eqn:konv3}, \eqref{eqn:konv31}, \eqref{eqn:konv20},
  \eqref{eqn:konv8} and \eqref{eqn:konv7}, as well as the assertions $\rm
  (1.b)$, $\rm (2.b)$ in Lemma \ref{lemma:konvergenzen}, and
  \eqref{eqn:konv5} we can pass to the limit in \eqref{eqn:sepp}.  The
  convergences \eqref{eqn:konv31}$_1$ and \eqref{eqn:konv8}$_2$ are
  used for the second term, while the convergences \eqref{eqn:konv20},
  \eqref{eqn:konv8}$_3$ and \eqref{eqn:konv5} are needed in the third
  and the fourth term. This implies the validity of \eqref{eqn:schwach}
  for $(b,\bphi)=(b,\F_{\eta}b)\in T_{\eta}^I$. The limiting process
  for test functions $(0,\bphi)\in T_\eta^I$ with $\bphi(T,\cdot)=0$
  and $\supp\bphi\subset\Omega^{\bar I}_{\eta}$ can be done like in the
  proof of Proposition \ref{lemma:approx}. 

  The length of the time interval depends only on the $L^\infty$ norm of
  $\eta$ at $t=0$ and on the bound for the H{\"o}lder norm
  of $\eta$. We have $\norm{\eta}_{L^\infty(I\times M)}<\kappa$.  Note that the
  quantities $\norm{\bu(t)}_{L^2(\oet)}$, $\norm{\eta(t)}_{H^2(M)}$,
  and $\norm{\pa_t\eta(t)}_{L^2(M)}$ are uniformly bounded for almost
  all $t\in I$. Constructing solutions with this initial data, we get from
  \eqref{ab:hs} that the H{\"o}lder norms of the displacements are bounded
  from above, independently of the initial time $t$. This, in turn, implies
  that the length of the existence interval is uniformly bounded from
  below. Choosing $t\in I=(0,T)$ sufficiently close to $T$, we obtain a weak solution on the interval $(t,\widetilde T)$, $\widetilde T>T$. Extending the original solution by the new one, by Remark \ref{bem:grenzen} we obtain a solution $(\widetilde\eta,\widetilde\bu)$ on the interval $(0,\widetilde T)$. Moreover, the solution $(\widetilde\eta,\widetilde\bu)$ satisfies
  the estimate \eqref{ab:hs} on the interval $(0,\widetilde T)$, since
  for $t_0<t<\widetilde T$ we have
  \begin{align*}
    &\onorm{\widetilde\bu(t,\cdot)}^2 +
    \int_0^t\norm{\nabla\widetilde\bu(s,\cdot) 
    }^p_{L^p(\Omega_{\widetilde\eta(s)})}\ ds +
    \mnorm{\pa_t\widetilde\eta(t,\cdot)}^2 +
    \norm{\widetilde\eta(t,\cdot)}_{H^2(M)}^2 \\ 
    &\hspace{.0cm}\le
    \big(\norm{\bu(t_0)}_{L^2(\Omega_{\widetilde\eta(t_0)})}^2 +
    \mnorm{\pa_t\eta(t_0)}^2 +
    \norm{\eta(t_0)}_{H^2(M)}^2\big)e^{c(t-t_0)}\\
    &\hspace{0.5cm} + \int_{t_0}^t
    \big(\norm{\ff(s,\cdot)}_{L^2(\Omega_{\eta(s)})}^2 +
    \mnorm{g(s,\cdot)}^2 \big)e^{c(t-s)}\
    ds + \int_0^{t_0}\norm{\nabla\bu(s,\cdot)
    }^p_{L^p(\Omega_{\widetilde\eta(s)})}\ ds\\
    &\hspace{.0cm}\le \big(\norm{\bu_0}_{L^2(\Omega_{\eta_0})}^2 +
    \mnorm{\eta_1}^2 + \norm{\eta_0}_{H^2(M)}^2\big)e^{ct}\\
    &\hspace{.5cm} + \int_0^t
    \big(\norm{\ff(s,\cdot)}_{L^2(\Omega_{\eta(s)})}^2 +
    \mnorm{g(s,\cdot)}^2\big)e^{c(t-s)}\ ds.
  \end{align*}
  Repeating this procedure, we obtain a maximal time $T^*\in
  (0,\infty]$ and a couple $(\eta,\bu)$ that solves our problem on
  each interval $(0,T)$, $T<T^*$, and that satisfies the estimate
  \eqref{ab:hs}.  If $T^*$ is finite, then this estimate implies that
  the H{\"o}lder norm of $\eta$ in $[0,T^*]\times M$ is bounded.
  Consequently, we have
  $\norm{\eta(t,\cdot)}_{L^\infty(M)}\nearrow\kappa$ for $t\nearrow
  T^*$.  \qed

\appendix
\section{Appendix}\label{sec:app} \label{sec:ode}

\begin{Proposition}\label{theorem:dicht}
  Let $\Omega\subset\setR^d$, $d\in\setN$, be a bounded domain with
  $C^0$ boundary and $1\le p<\infty$. Then $C^\infty_0(\setR^d)$ is
  dense in $W^{1,p}(\Omega)$ and in $E^p(\Omega)$.
\end{Proposition}
\proof The proof for $W^{1,p}(\Omega)$ can be found in
\cite[Theorem 5.5.9]{b5}, and it can be carried over to prove the assertion for $E^p(\Omega)$. For the convenience of the reader we briefly sketch the proof. Let $B\subset\setR^{d-1}$ be an open
ball, $g:\overline B\rightarrow\setR$ a continuous function and let
$\widetilde\Omega:=\{(x',x_d)\in\setR^d \ |\ x'\in B,$ 
$g(x')< x_d\}$. Let $v\in W^{1,p}(\widetilde\Omega)$ be a function
with bounded support in $\{(x',x_d)\in\setR^d\ |\ x'\in B,$ $\,
g(x')\le x_d\}$. We extend $v$ by $0$ to $\setR^d$. Let
$v_t(x):=v(x',x_d+t)$ with $t>0$. Thus, for sufficiently small
$\epsilon>0$, we have $\nabla(\omega_\epsilon\ast
v_t)=\omega_\epsilon\ast\nabla v_t$ in $\widetilde\Omega$. Here,
$\omega_\epsilon$ is a standard mollification kernel.  Clearly, we obtain
that $\omega_\epsilon\ast v_t\in C^\infty_0(\setR^d)$ converges to
$v_t$ in $W^{1,p}(\widetilde\Omega)$ for $\epsilon\rightarrow 0$. The 
convergence of $v_t$ to $v$ in $W^{1,p}(\widetilde\Omega)$ follows
directly from the continuity of translation in $L^p(\setR^d)$.
Now, the denseness in $W^{1,p}(\Omega)$ follows by localization. The denseness in $E^p(\Omega)$ can be shown in the same way by simply replacing $\nabla$ by $\dv$.  \qed
\medskip

The following classical result can be found in  \cite{b13}.
\begin{Proposition}{\bf (Reynolds
    transport theorem)}\label{theorem:reynolds}
  Let $\Omega\subset\setR^3$ be a bounded domain with $C^1$-boundary,
  let $I\subset\setR$ be an interval, and let $\Psi\in
  C^1(I\times\overline\Omega,\setR^3)$ such that
  \[
    \Psi_t:=\Psi(t,\cdot):\overline\Omega\rightarrow
    \Psi_t(\overline\Omega)
  \] 
  is a diffeomorphism for all $t\in I$. We set
  $\Omega_t:=\Psi_t(\Omega)$ and $\bv:=(\pa_t\Psi)\circ\Psi_t^{-1}$.
  Then we have  for  all $\xi\in C^1(\bigcup_{t\in I}\{t\}\times
  \overline\Omega_t)$ and $t\in I$ 
  \[
    \frac{d}{dt}\int_{\Omega_t} \xi(t,x)\ dx= \int_{\Omega_t} \pa_t
    \xi(t,x)\ dx + \int_{\pa\Omega_t} \bv\cdot \bnu_t\ \xi(t,\cdot)\
    dA_t.
  \]
  Here, $dA_t$ denotes the surface measure and $\bnu_t$ denotes the
  outer unit normal of $\pa\Omega_t$.
\end{Proposition}

Let us now state an existence result for ordinary integro-differential
equations.  
\begin{Proposition}\label{pro:ode}
  Let $d\in\setN$, $\bA\in C([0,\infty)\times\setR^d,\setR^d)$ and
  $\bB\in C([0,\infty)^2\times\setR^d,\setR^d)$. Then, for
  all $\balpha_0\in\setR^d$ there exists a $T^*\in (0,\infty]$ and a solution 
  $\balpha\in C^1([0,T^*),\setR^d)$ of 
  \begin{equation}\label{eqn:gidgl}
    \begin{aligned}
      \dot \balpha(t)&=\bA(t,\balpha(t)) + \int_0^t
      \bB(t,s,\balpha(s))\ ds\,,
      \\
      \balpha(0)&=\balpha_0\,.
    \end{aligned}
  \end{equation}
  for all $t\in [0,T^*)$. If $T^*<\infty$ then $\lim_{t\nearrow
    T^*}|\balpha(t)|=\infty$. If $\bA$ and $\bB$ are
  affine linear in $\balpha$, then  $T^*=\infty$.
\end{Proposition}
\proof The proof is an easy adaptation of the proof
of the existence theorem of Peano. Details can be found in
\cite{b62}. \qed
\medskip

Let us now consider the inhomogeneous Stokes system
\begin{equation}
  \begin{aligned}
    -\Delta\bu+\nabla\pi&=\bfzero &&\text{ in }\Omega,\\
    \dv \bu&=0&&\text{ in }\Omega,\\
    \bu&=\bg&&\text{ auf }\pa\Omega.
  \end{aligned}\label{eq:stokes}
\end{equation}

\begin{Theorem}\label{theorem:stokes}
  Let $\Omega\subset\setR^3$ be a bounded domain with a boundary of class
  $C^{\max(2,k)}$, $k\in\setN$, and let ${1<p<\infty}$.
  Moreover, let $\bg\in W^{k-1/p,p}(\pa\Omega)$ satisfy 
  \[
     \int_{\pa\Omega}\bg\cdot\bnu\ dA=0,
  \] 
  where $\bnu$ and $dA$ denote the unit outer normal and the surface
  measure of $\pa\Omega$, respectively. Then there exists exactly
  one (very weak) solution $\bu\in W^{k,p}(\Omega)$ of
  \eqref{eq:stokes}, i.e., we have $\dv\bu=0$, and
  \[
    -\int_\Omega\bu\cdot\Delta\bphi\
    dx+\int_{\pa\Omega}\bg\cdot\pa_{\bnu}\bphi\ dA=0 
  \]
  is satisfied for all $\bphi\in C^2(\overline\Omega)$ with $\dv\bphi=0$ and
  $\bphi=0$ on $\pa\Omega$. The mapping $\bg\mapsto\bu$ defines
  a continuous, linear operator from $W^{k-1/p,p}(\pa\Omega)$ to 
  $W^{k,p}(\Omega)$.

  For $\pa\Omega$ belonging to $C^{2,1}$, the statement holds also if
  we replace $W^{k-1/p,p}(\pa\Omega)$ by $L^p(\pa\Omega)$ and
  $W^{k,p}(\Omega)$ by $L^p(\Omega)$, respectively. 
\end{Theorem}
\proof The first assertion is proved in \cite{b18},
while the second one is a direct consequence of Theorem 3 in
\cite{b19}.\qed 
\medskip

A proof of the following variant of Schauder's  fixed point theorem
for set-valued mappings can be found in \cite{b29}.  
\begin{Theorem}{\bf (Kakutani-Glicksberg-Fan)}\label{theorem:kakutani}
  Let $C$ be a convex subset of a normed vector space $Z$ and let 
  $F: C\rightarrow 2^C$ be an upper-semicontinuous set-valued mapping,
  i.e., for every open set $W\subset C$ the set $\{c\in
  C\ |\ F(c)\subset W\}\subset C$ is open. Moreover, let $F(C)$ be
  contained in a compact subset of $C$, and let  $F(z)$ be 
  non-empty, convex, and compact for all $z\in C$. Then $F$ possesses 
  a fixed point, i.e., there exists a $c_0\in C$ with $c_0\in
  F(c_0)$.
\end{Theorem}
It is not hard to see that the requirement of
upper-semicontinuity is equivalent to the requirement that the graph
of $F$ is closed, i.e., the convergences $c_n\rightarrow c$ in $C$ and
$z_n\rightarrow z$ in $Z$, where $z_n \in F(c_n) $, imply that 
$z\in F(c)$.

\begin{Lemma}\label{bem:korn}
In the context of subsection \ref{statement} and assuming all involved functions to be sufficiently regular we have
\[\int_{\Omega_{\eta(t)}} (\nabla\bu)^T:\nabla\bphi\ dx =0\]
for all functions $\bphi$ with $\tr \bphi=b\,\bnu$ for some scalar function $b$, in particular for
$\bphi=\bu$.
\end{Lemma}
\proof In view of $\dv\bu=0$ integration by parts shows that
\[\int_{\Omega_{\eta(t)}} (\nabla\bu)^T:\nabla\bphi\ dx = \int_{\pa\Omega_{\eta(t)}}
(\bphi\cdot\nabla)\bu\cdot\bnu_{\eta(t)}\ dA.\]
Thus, it suffices to prove that on $\pa\Omega_{\eta(t)}$ we have
\[((\bnu\circ q)\cdot\nabla)\bu\cdot\bnu_{\eta(t)}=0.\]
To this end, we set $\be_1:=\bnu\circ q$. Moreover, on
$\pa\Omega_{\eta(t)}$ we choose two linearly independent, tangential vector fields\footnote{These exist locally.} and extend these
constantly along $\be_1$. We denote the resulting vector fields by $\be_2$ und $\be_3$. Thus, we have
$\Gamma_{1i}^k\,\be_k:=\nabla_{\be_1}\be_i=0$ and hence $\Gamma_{1i}^k=0$ for all $i,k$.
Letting $\bu=u^i\,\be_i$ we infer that
\[((\bnu\circ
q)\cdot\nabla)\bu=\nabla_{\be_1}\bu=du^i\,\be_1\ \be_i+u^i\,\nabla_{\be_1}\be_i=du^i\,\be_1\
\be_i,\]
and thus on $\pa\Omega_{\eta(t)}$
\begin{equation}\label{eqn:zw}
 \begin{aligned}
((\bnu\circ
q)\cdot\nabla)\bu\cdot\bnu_{\eta(t)}=du^1\,\be_1\ \be_1\cdot\bnu_{\eta(t)}.  
 \end{aligned}
\end{equation}
On the other hand we have
\[0=\dv\bu=du^i\,\be_i+ u^k\,\Gamma_{ik}^i.\]
The components $u^2$ and $u^3$ as well as their tangential derivatives $du^2\,\be_2$ respectively $du^3\,\be_3$ vanish on
$\pa\Omega_{\eta(t)}$. Hence, on $\pa\Omega_{\eta(t)}$
\[0=du^1\,\be_1+u^1\,\Gamma_{1i}^i=du^1\,\be_1.\]
This identity and \eqref{eqn:zw} prove the claim.\qed

\begin{Lemma}\label{lemma:mittelwert}
Let $\eta\in Y^I$ with $\norm{\eta}_{L^\infty(I\times M)}<\kappa$. There exists a linear operator
$\M_\eta$ such that
\begin{equation*}
 \begin{aligned}
  \norm{\M_\eta b}_{L^r(I\times M)}&\le c\, \norm{b}_{L^r(I\times M)},\\
  \norm{\M_\eta b}_{C(\bar I,L^r(M)}&\le c\, \norm{b}_{C(\bar I,L^r(M)},\\
  \norm{\M_\eta b}_{L^r(I,H^2_0(M))}&\le c\, \norm{b}_{L^r(I,H^2_0(M))},\\
  \norm{\M_\eta b}_{H^1(I,L^2(M))}&\le c\, \norm{b}_{H^1(I,L^2(M))}\\
 \end{aligned}
\end{equation*}
for all $1\le r\le\infty$ and
\[\int_M (\M_\eta b)(t,\cdot)\,\gamma(\eta(t,\cdot))\ dA = 0\]
for almost all $t\in I$. The constant $c$ depends only on $\Omega$, $\norm{\eta}_{Y^I}$ and
$\tau(\eta)$; it stays bounded as long as $\norm{\eta}_{Y^I}$ and $\tau(\eta)$ stay bounded.
\end{Lemma}
\proof We fix an arbitrary function $\psi\in C_0^\infty(\inn M)$ with $\psi\ge
0$, $\psi\not\equiv 0$ and define
\begin{equation*}
 \begin{aligned}
  (M_\eta b)(t,\cdot):= b- \psi\, \frac{a(b(t,\cdot),\eta(t,\cdot))}{a(\psi,\eta(t,\cdot))},
 \end{aligned}
\end{equation*}
where
\[a(b(t,\cdot),\eta(t,\cdot)):=\int_M b\,\gamma(\eta(t,\cdot))\ dA.\]
Now it's easy to prove the claims provided we note that $a(\psi,\eta)\ge c$ with a constant $c>0$ depending only
$\tau(\eta)$. But this follows from Remark \ref{bem:groessernull} since $\gamma$ is a continuous function of $\eta$. 
\qed

\pagebreak 
\begin{Lemma}\label{lemma:konvergenzen}
Let the sequence $(\eta_n)\subset Y^I$ satisfiy $\sup_n\norm{\eta_n}_{L^\infty(I\times
M)}<\alpha<\kappa$ and \eqref{eqn:schwkonv}$_{(1,2)}$.
\begin{itemize}
 \item[\rm (1.a)] Provided that $b\in
C(\bar I,L^2(M))$ the sequence $(\M_{\eta_n}b)$ converges to $\M_{\eta}b$ in $C(\bar I,L^2(M))$ independently of $\norm{b}_{C(\bar
I,L^2(M))}\le 1$.
\item[\rm (1.b)] Provided that $b\in H^1(I,L^2(M))\cap
  L^2(I,H^2_0(M))$ and, additionally, $(\pa_t\eta_n)$ converges in \linebreak
  $L^2(I\times M)$ the sequence $(\M_{\eta_n}b)$ converges to $\M_\eta
  b$ in $H^1(I,L^2(M))\cap L^2(I,H^2_0(M))$.
\item[\rm (2.a)] Provided that $b\in C(\bar I,L^2(M))$ the sequence
  $(\F_{\eta_n}\M_{\eta_n}b)$ converges to $\F_{\eta}\M_{\eta}b$ in \linebreak
  $C(\bar I,L^2(B_\alpha))$ independently of $\norm{b}_{C(\bar
    I,L^2(M))}\le 1$.
 \item[\rm (2.b)] On the conditions of $(1.b)$ the sequence $(\F_{\eta_n}\M_{\eta_n}b)$
converges to $\F_{\eta}\M_\eta b$ in
$H^1(I\times B_\alpha)\cap L^\infty(I,L^4(B_\alpha))$.
\item[\rm (2.c)] Provided that $(b_n)$ converges to $b$ weakly in $L^2(I\times M)$ the sequence $(\F_{\eta_n}b_n)$ converges to
$\F_\eta b$ weakly in $L^2(I\times B_\alpha)$.
\end{itemize}
\end{Lemma}
\proof $\rm (1.a)$ is a consequence of
\begin{equation*}
 \begin{aligned}
  \norm{a(b,\eta_n)-a(b,\eta)}_{L^\infty(I)}&=\Bignorm{\int_M b\,
(\gamma(\eta_n)-\gamma(\eta))\ dA}_{L^\infty(I)}\\
&\le c\,\norm{b}_{C(\bar I,L^2(M))}\,\norm{\gamma(\eta_n)-\gamma(\eta)}_{L^\infty(I\times M)},
 \end{aligned}
\end{equation*}
the analogous estimate with $b$ replaced by $\psi$, the inequality $a(\psi,\eta_n)\ge c>0$, and the fact that $(\gamma(\eta_n))$
converges to $\gamma(\eta)$ uniformly. Assertion $(1.b)$ follows from the convergence of $(a(b,\eta_n))$, $(a(\psi,\eta_n))$
to $a(b,\eta)$ respectively $a(\psi,\eta)$ in $H^1(I)$. The proof of these convergences is very easy.

Now let us prove $(2.a)$. We infer from
\eqref{eqn:fort} and $(1.a)$ that $(\F_{\eta_n}\M_{\eta_n}b)$ converges in $C(\bar I,L^2(S_\alpha))$
independently of $\norm{b}_{C(\bar I,L^2(M))}\le 1$. Similarly, the formal trace
\begin{equation*}
\exp\Big(\int_{\eta_n(t,q)}^{-\alpha}
\beta(q+\tau\,\bnu\circ q))\ d\tau\Big)\,
(\M_{\eta_n}b)(t,q)\,\bnu\circ q 
\end{equation*}
of $\F_{\eta_n}\M_{\eta_n}b$ on
$I\times\pa(\Omega\setminus\overline{S_\alpha})$ converges in $C(\bar
I,L^2(\Omega\setminus\overline{S_\alpha}))$ independently of
${\norm{b}_{C(\bar I,L^2(M))}\le 1}$. Thus, the claim follows from the
continuity properties of the solution operator of the Stokes system,
cf. Remark \ref{bem:nspur}. The proof of $\rm (2.c)$ proceeds
analogously.

Now we show $\rm (2.b)$. From \eqref{eqn:fort}, $(1.b)$, and
\[H^1(I,L^2(M))\cap L^2(I,H^2_0(M))\embedding C(\bar I,H^1(M))\embedding L^\infty(I,L^4(M))\]
we deduce the convergence of $(\F_{\eta_n}\M_{\eta_n}b)$ in $L^\infty(I,L^4(S_\alpha))$. In order to prove the convergence in
$H^1(I\times S_\alpha)$ we infer by interpolation $(\theta=2/3)$ and
Sobolev embedding that
\begin{equation}\label{eqn:lurchi}
  L^\infty(I,L^2(M))\cap L^2(I,H^2_0(M)) \embedding L^{3}(I,H^{4/3}(M))\embedding
  L^3(I,L^\infty(M)).
\end{equation}
Now, the convergence in $L^2(I,H^1(S_\alpha))$ is consequence of
\eqref{eqn:fortnabla}, $\rm (1.b)$, \eqref{eqn:lurchi}, and
\[Y^I\compactembedding L^6(I,H^1_0(M)),\]
while the convergence in $H^1(I,L^2(S_\alpha))$ can be inferred from \eqref{eqn:klopps}, $\rm (1.b)$, \eqref{eqn:lurchi}, and the
convergence of $(\pa_t\eta_n)$ in $L^6(I,L^2(M))$. The latter follows by interpolation.
Similarly we can show the convergence of the trace in $H^1(I\times\pa(\Omega\setminus\overline{S_\alpha}))\cap
L^\infty(I,L^4(\pa(\Omega\setminus\overline{S_\alpha})))$. As before, in order to prove the claim it suffices now to invoke the
continuity properties of the solution operator of the Stokes system.
\qed

\begin{Lemma}\label{lemma:ehrling} For  all $N\in\setN$, $3/2<p\le\infty$ and $\epsilon>0$ 
there exists a constant $c$ such that for all $\eta,\,\tilde\eta\in H^2_0(M)$ with
$\norm{\eta}_{H^2_0(M)}+\norm{\tilde\eta}_{H^2_0(M)}+\tau(\eta)+\tau(\tilde\eta)\le N$
and all $\bv\in W^{1,p}(\Omega_\eta)$, ${\tilde\bv\in W^{1,p}(\Omega_{\tilde\eta})}$ we have
\begin{equation*}
 \begin{aligned}
&\sup_{\norm{b}_{L^2(M)}\le
1}\bigg(\int_{\Omega_{\eta}} \bv\cdot\F_{\eta}\M_{\eta}b\ dx - \int_{\Omega_{\tilde\eta}}
\tilde\bv\cdot\F_{\tilde\eta}\M_{\tilde\eta}b\ dx\\
&\hspace{2cm} + \int_M \tr\bv\cdot\bnu\ \M_{\eta}b - \trt\tilde\bv\cdot\bnu\ \M_{\tilde\eta}b\
dA\bigg)\\
&\hspace{0.5cm}\le c\sup_{\norm{b}_{H^2_0(M)}\le
1}\bigg(\int_{\Omega_{\eta}} \bv\cdot\F_{\eta}\M_{\eta}b\ dx - \int_{\Omega_{\tilde\eta}}
\tilde\bv\cdot\F_{\tilde\eta}\M_{\tilde\eta}b\ dx\\
&\hspace{1cm}  + \int_M \tr\bv\cdot\bnu\ \M_{\eta}b -
\trt\tilde\bv\cdot\bnu\ \M_{\tilde\eta}b\ dA\bigg) + \epsilon\, \big(\norm{\bv}_{W^{1,p}(\Omega_{\eta})} +
\norm{\tilde\bv}_{W^{1,p}(\Omega_{\tilde\eta})}\big).
 \end{aligned}
\end{equation*}
Similarly, for all $N\in\setN$, $6/5<p,r\le\infty$ and $\epsilon>0$ there exists a constant $c$
such that for all $\eta,\,\tilde\eta\in
H^2_0(M)$ and $\delta\in C^4(\pa\Omega)$
with
\[\norm{\eta}_{H^2_0(M)}+\norm{\tilde\eta}_{H^2_0(M)}+
\norm{\delta}_{C^4(\pa\Omega)}+\tau(\eta)+\tau(\tilde\eta)+\tau(\delta)\le N\] 
and all $\bv\in W^{1,p}(\Omega_\eta)$, $\tilde\bv\in W^{1,p}(\Omega_{\tilde\eta})$ we have
\begin{equation*}
 \begin{aligned}
&\sup_{\norm{\bphi}_{H(\Omega)}\le
1}\bigg(\int_{\Omega_{\eta}} \bv\cdot\T_\delta\bphi\ dx - \int_{\Omega_{\tilde\eta}}
\tilde\bv\cdot\T_\delta\bphi\ dx\bigg)\\
&\hspace{1cm}\le c\sup_{\norm{\bphi}_{W^{1,r}_{0,\dv}(\Omega)}\le
1}\bigg(\int_{\Omega_{\eta}} \bv\cdot\T_\delta\bphi\ dx - \int_{\Omega_{\tilde\eta}}
\tilde\bv\cdot\T_\delta\bphi\ dx\bigg)\\
&\hspace{4cm}+\epsilon\, \big(\norm{\bv}_{W^{1,p}(\Omega_{\eta})} +
\norm{\tilde\bv}_{W^{1,p}(\Omega_{\tilde\eta})}\big) .
 \end{aligned}
\end{equation*}
\end{Lemma}
\proof We prove these assertions of Ehrling lemma-type by the usual contradiction argument. Let us start with the first claim.
Assuming that it is wrong there exists a $3/2<p\le\infty$, an $\epsilon>0$, bounded sequences
$(\eta_n),\,(\tilde\eta_n)\subset H^2_0(M)$ with
$\sup_n\big(\tau(\eta_n)+\tau(\tilde\eta_n)\big)<\infty$ as well as sequences
$(\bv_n),\,(\tilde\bv_n)$ with
\[\norm{\bv_n}_{W^{1,p}(\Omega_{\eta_n})}
+ \norm{\tilde\bv_n}_{W^{1,p}(\Omega_{\tilde\eta_n})}=1\]
and 
\begin{equation}\label{eqn:wid}
 \begin{aligned}
&\sup_{\norm{b}_{L^2(M)}\le
1}\bigg(\int_{\Omega_{\eta_n}} \bv_n\cdot\F_{\eta_n}\M_{\eta_n}b\ dx - \int_{\Omega_{\tilde\eta_n}}
\tilde\bv_n\cdot\F_{\tilde\eta_n}\M_{\tilde\eta_n}b\ dx\\
&\hspace{2cm} + \int_M \tren\bv_n\cdot\bnu\ \M_{\eta_n}b - \trent\tilde\bv_n\cdot\bnu\
\M_{\tilde\eta_n}b\
dA\bigg)\\
&\hspace{1cm}> \epsilon + n\sup_{\norm{b}_{H^2_0(M)}\le
1}\bigg(\int_{\Omega_{\eta_n}} \bv_n\cdot\F_{\eta_n}\M_{\eta_n}b\ dx - \int_{\Omega_{\tilde\eta_n}}
\tilde\bv_n\cdot\F_{\tilde\eta_n}\M_{\tilde\eta_n}b\ dx\\
&\hspace{4.5cm}  + \int_M \tren\bv_n\cdot\bnu\ \M_{\eta_n}b -
\trent\tilde\bv_n\cdot\bnu\ \M_{\tilde\eta_n}b\ dA\bigg).
\end{aligned}
\end{equation}
Due to Corollary \ref{lemma:spur} the sequences $(\tren\bv_n)$, $(\trent\tilde\bv_n)$ are bounded in
$W^{1-1/r,r}(M)$ for some $r>3/2$ and by Sobolev embedding in $H^s(M)$ for
some $s>0$. Thus, there exist subsequences with
\begin{equation*}
 \begin{aligned}
\tren\bv_n\cdot\bnu\rightarrow d,\ \trent\tilde\bv_n\cdot\bnu\rightarrow \tilde d&\text{\quad in
}L^2(M),\\
\eta_n\rightarrow\eta,\ \tilde\eta_n\rightarrow\tilde\eta &\text{\quad weakly in }
H^2_0(M) \text{, in particular uniformly}. 
 \end{aligned}
\end{equation*}
Lemma \ref{lemma:psi} and the common Sobolev embeddings show that a subsequence of
$(\bw_n:=\bv_n\circ\Psi_{\eta_n})$ converges to some $\bw$ in $L^3(\Omega)$. If we extend all involved functions by
$\boldsymbol{0}$ to $\setR^3$, then from the estimate
\begin{equation*}
 \begin{aligned}
  \norm{\bv_n-\bw\circ \Psi^{-1}_{\eta}}_{L^2(\setR^3)} \le \norm{(\bw_n-\bw)\circ
\Psi^{-1}_{\eta_n}}_{L^2(\setR^3)} + \norm{\bw\circ\Psi^{-1}_{\eta_n} -\bw\circ
\Psi^{-1}_{\eta}}_{L^2(\setR^3)},
 \end{aligned}
\end{equation*}
Lemma \ref{lemma:psi} and Remark \ref{bem:konv} we infer that $(\bv_n)$ converges to $\bv:=\bw\circ \Psi^{-1}_{\tilde\eta}$ in
$L^2(\setR^3)$. Similarly, the sequence $(\tilde\bv_n)$ converges to $\tilde \bv$ in $L^2(\setR^3)$. Moreover, Lemma
\ref{lemma:konvergenzen} $(1.a)$, $(2.a)$ show that the sequences $(\M_{\eta_n}b)$ and $(\M_{\tilde\eta_n}b)$ converge in $L^2(M)$
while $(\F_{\eta_n}\M_{\eta_n}b)$, $(\F_{\tilde\eta_n}\M_{\tilde\eta_n}b)$ converge in $L^2(\setR^3)$, each independently of
$\norm{b}_{L^2(M)}\le 1$. Hence, the supremum on the right hand side of \eqref{eqn:wid} tends to
\begin{equation*}
 \begin{aligned}
\sup_{\norm{b}_{H^2_0(M)}\le
1}\bigg(\int_{\Omega_{\eta}} \bv\cdot\F_{\eta}\M_{\eta}b\ dx - \int_{\Omega_{\tilde\eta}}
\tilde\bv\cdot\F_{\tilde\eta}\M_{\tilde\eta}b\ dx + \int_M d\ \M_{\eta}b - \tilde d\
\M_{\tilde\eta}b\ dA\bigg).
 \end{aligned}
\end{equation*}
Since the left hand side of \eqref{eqn:wid} is bounded this limit must vanish. In view of the denseness of $H^2_0(M)$ in $L^2(M)$
the limit
\begin{equation*}
 \begin{aligned}
\sup_{\norm{b}_{L^2(M)}\le
1}\bigg(\int_{\Omega_{\eta}} \bv\cdot\F_{\eta}\M_{\eta}b\ dx - \int_{\Omega_{\tilde\eta}}
\tilde\bv\cdot\F_{\tilde\eta}\M_{\tilde\eta}b\ dx + \int_M d\ \M_{\eta}b - \tilde d\
\M_{\tilde\eta}b\ dA\bigg) 
\end{aligned}
\end{equation*}
of the left hand side of \eqref{eqn:wid} must vanish, too, contradicting $\epsilon>0$.

The proof of the second claim proceeds analogously. Therefore, we
merely show that for bounded sequences $(\eta_n)\subset H^2_0(M)$,
$(\delta_n)\subset C^4(\pa\Omega)$ and a sequence $(\bv_n)$ with
\[\sup_n\big(\norm{\bv_n}_{W^{1,p}(\Omega_{\eta_n})} + \tau(\eta_n)+\tau(\delta_n)\big)<\infty\]
there exist subsequences such that
\begin{equation*}
 \begin{aligned}
   \int_{\Omega_{\eta_n}} \bv_n\cdot\T_{\delta_n}\bphi\
   dx\quad\rightarrow\quad \int_{\Omega_{\eta}}
   \bv\cdot\T_{\delta}\bphi\ dx
 \end{aligned}
\end{equation*}
for $n\rightarrow\infty$ independently of $\norm{\bphi}_{H(\Omega)}\le 1$. As before we can find
subsequences such that
\begin{equation*}
 \begin{aligned}
\bv_n\rightarrow \bv&\text{\quad in }L^2(\setR^3),\\
\eta_n\rightarrow\eta &\text{\quad weakly in }
H^2_0(M) \text{, in particular uniformly},\\
\delta_n\rightarrow\delta &\text{\quad in } C^3(\pa\Omega).
 \end{aligned}
\end{equation*}
Adding a zero-sum we obtain the estimate
\begin{equation*}
 \begin{aligned}
  \big|\int_{\Omega_{\eta_n}}\bv_n\cdot\T_{\delta_n}\bphi\
dx-\int_{\Omega_{\eta}}\bv\cdot\T_{\delta}\bphi\ dx\big|&\le
\norm{\bv_n-\bv}_{L^2(\setR^3)}\,\norm{\T_{\delta_n}\bphi}_{L^2(\setR^3)}\\
&\hspace{0.5cm}+\norm{\bv}_{W^{1,p}(\Omega_\eta)}\,\norm{\T_{\delta_n}\bphi-\T_\delta\bphi}_{
(W^{1,p}(\Omega_\eta))'}.
 \end{aligned}
\end{equation*}
Hence, we merely have to show that $(\T_{\delta_n}\bphi)$ converges to $\T_{\delta}\bphi$
in $(W^{1,p}(\Omega_\eta))'$ independently of $\norm{\bphi}_{H(\Omega)}\le 1$. Assuming that this is false there exists an
$\epsilon>0$ and a sequence $(\bphi_n)\subset H(\Omega)$
converging to some $\bphi$ weakly in $H(\Omega)$ such that
\[\norm{\T_{\delta_n}\bphi_n-\T_{\delta}\bphi_n}_{(W^{1,p}(\Omega_\eta))'}>\epsilon.\]
But, in fact, this contradicts
\begin{equation}\label{ab:fin}
 \begin{aligned}
\norm{\T_{\delta_n}\bphi_n-\T_{\delta}\bphi_n}_{(W^{1,p}(\Omega_\eta))'}\le
&\norm{\T_{\delta_n}(\bphi_n-\bphi)}_{(W^{1,p}(\Omega_\eta))'}\\
& +\norm{\T_{\delta_n}\bphi-\T_{\delta}\bphi_n}_{(W^{1,p}(\Omega_\eta))'}.
 \end{aligned}
\end{equation}
In this respect note that the sequences $(\T_{\delta}\bphi_n)$ and $(\T_{\delta_n}\bphi)$ converge to $\T_\delta\bphi$
weakly respectively strongly in $L^2(\setR^3)$. The strong convergence of $(\T_{\delta_n}\bphi)$ can be easily shown by
approximating $\bphi$ by smooth functions; cf. Remark \ref{bem:konv}. Moreover, the identity
\[\int_{\setR^3} \T_{\delta_n}(\bphi_n-\bphi)\ h\ dx = \int_\Omega
d\Psi_{\delta_n}(\bphi_n-\bphi)\ h\circ\Psi_{\delta_n}\ dx,\]
$h\in L^2(\setR^3)$, shows that $(\T_{\delta_n}(\bphi_n-\bphi))$ converges to $\boldsymbol{0}$ weakly in $L^2(\setR^3)$. Now, from
\[L^2(\Omega_\eta)\compactembedding (W^{1,p}(\Omega_\eta))'\]
we infer that the right hand side of \eqref{ab:fin} gets small for large $n$.
\qed
\medskip

\begin{Lemma}\label{lemma:divdichtglm}
For all $N\in\setN$, $s>0$ and $\epsilon>0$ there exists a small $\sigma>0$ such that
for all $\eta\in H^2_0(M)$ with $\norm{\eta}_{H^2_0(M)}+\tau(\eta)\le N$ and all $\bphi\in
H(\Omega_\eta)$ with $\norm{\bphi}_{L^2(\Omega_\eta)}\le
1$ there exists a $\bpsi\in H(\Omega_\eta)$ with
$\norm{\bpsi}_{L^2(\Omega_\eta)}\le
2$, $\norm{\bphi-\bpsi}_{(H^{s}(\setR^3))'}<\epsilon$ and $\supp \bpsi\subset\Omega_{\eta-\sigma}$.\footnote{When constructing
$\Omega_{\eta-\sigma}$ we \emph{first} extend $\eta$ by $0$ to $\pa\Omega$ and \emph{then} subtract $\sigma$.
}
\end{Lemma}
\proof Let us assume that the claim is wrong. Then there exist $s,\,\epsilon>0$, a sequence $(\sigma_n)_{n\in\setN}$ of positive
numbers converging to $0$, a bounded sequence $(\eta_n)\subset H^2_0(M)$ with $\sup_n\tau(\eta_n)<\infty$ and a sequence
$(\bphi_n)\subset H(\Omega_{\eta_n})$ with
$\norm{\bphi_n}_{L^2(\Omega_{\eta_n})}\le 1$, $\norm{\bphi_n-\bpsi}_{(H^{s}(\setR^3))'}\ge
\epsilon$ for all $\bpsi\in H(\Omega_{\eta_n})$ with $\norm{\bpsi}_{L^2(\Omega_{\eta_n})}\le
2$ and $\supp \bpsi\subset\Omega_{\eta_n-\sigma_n}$ as well as
\begin{equation*}
 \begin{aligned}
  \eta_n&\rightarrow\eta &&\text{ weakly in }H^2_0(M), \text{ in particular uniformly},\\
  \bphi_n&\rightarrow \bphi &&\text{ weakly in }L^2(\setR^3).
\end{aligned}
\end{equation*}
From the compact embedding
\[L^2(B)\compactembedding  (H^{s}(B))'\]
for a suitable ball $B\subset\setR^3$ we deduce the strong convergence of $(\bphi_n)$ in $(H^{s}(\setR^3))'$ using some
extension operator for fractional Sobolev spaces, see, e.g., \cite{b4}. Let $\psi$ be a $C^1$-function defined in a
neighbourhood of $\Omega_\eta$. Then Proposition \ref{lemma:nspur} shows that
\begin{equation*}
  \int_{\Omega_\eta} \bphi\cdot\nabla\psi\ dx = \lim_n \int_{\Omega_{\eta_n}}
\bphi_n\cdot\nabla\psi\ dx =0,
\end{equation*}
because each term of the sequence vanishes. Hence, $\trnormal \bphi = 0$, and by Proposition
\ref{lemma:divdicht} there exists a $\bpsi\in H(\Omega_\eta)$ with $\supp \bpsi\subset\Omega_\eta$ and
$\norm{\bphi-\bpsi}_{L^2(\Omega_\eta)}<\epsilon/2$. It follows that for sufficiently large $n$
\begin{equation*}
  \norm{\bphi_n-\bpsi}_{(H^{s}(\setR^3))'}\le
\norm{\bphi_n-\bphi}_{(H^{s}(\setR^3))'} +
\norm{\bphi-\bpsi}_{L^2(\Omega_{\eta})}< \epsilon.
\end{equation*}
This is a contradiction provided that $\norm{\bpsi}_{L^2(\Omega_{\eta_n})}\le
2$. But this estimate is a consequence of
$\norm{\bphi-\bpsi}_{L^2(\Omega_\eta)}<\epsilon/2$
for sufficiently small $\epsilon$. The latter can be assumed without loss of generality.
\qed

\begin{Remark}\label{bem:ausw}
Let us show that for $\eta\in Y^I$, $\norm{\eta}_{L^\infty(I\times
M)}<\alpha<\kappa$ and $(b,\bphi)\in T^I_\eta$ the extension of
$\bphi$ by
$(b\,\bnu)\circ q$ to $I\times B_\alpha$ is in $H^1(I,L^2(B_\alpha))$. We approximate
$\bphi$ by functions $(\bphi_k)\subset C_0^\infty(\setR^4)$ in $H^1(\Omega_\eta^I)$ and $\eta$ by $(\eta_n)\subset
C^4(\bar I\times\pa\Omega)$ such that $(\eta_n)$ converges to $\eta$ in
$L^\infty(I\times\pa\Omega)$ and $(\pa_t\eta_n)$ converges to $\pa_t\eta$ in $L^2(I\times\pa\Omega)$; cf. Definition
\eqref{eqn:reg}. Using Reynolds' transport theorem with $\bv=(\pa_t\eta_n\bnu)\circ
\Phi_{\eta_n(t)}^{-1}$ and
$\boldsymbol{\xi}=\bphi_k\,\psi$, $\psi\in C_0^\infty(I\times\setR^3)$ we obtain
\begin{equation*}
 \begin{aligned}
0 = \int_I\frac{d}{dt}\int_{\Omega_{\eta_n(t)}}
\bphi_k\,\psi\ dxdt&=\int_I\int_{\Omega_{\eta_n(t)}} \pa_t\bphi_k\, \psi +
\bphi_k\, \pa_t\psi\ dxdt\\
&\hspace{0.5cm}+ \int_I\int_{\pa\Omega} \tren(\bphi_k\,\psi)\, \pa_t\eta_n\, \gamma(\eta_n)\ dAdt.
 \end{aligned}
\end{equation*}
Letting first $n$ and then $k$ tend to infinity results in
\begin{equation*}
 \begin{aligned}
\int_I\int_{\Omega_{\eta(t)}} \bphi\, \pa_t\psi\ dxdt=-
\int_I\int_{\Omega_{\eta(t)}} \pa_t\bphi\, \psi\ dxdt-\int_I\int_{\pa\Omega}b\,\bnu\,
\tr\psi\ \pa_t\eta\, \gamma(\eta)\ dAdt.
 \end{aligned}
\end{equation*}
Similarly we show that
\begin{equation*}
 \begin{aligned}
\int_I\int_{B_\alpha\setminus\overline{\Omega_{\eta(t)}}} (b\,\bnu)\circ q\ \pa_t\psi\ dxdt&=-
\int_I\int_{B_\alpha\setminus\overline{\Omega_{\eta(t)}}} (\pa_tb\ \bnu)\circ q\ \psi\ dxdt\\
&\hspace{0.5cm}+\int_I\int_{\pa\Omega}b\,\bnu\, \tr\psi\ \pa_t\eta\, \gamma(\eta)\ dAdt.
 \end{aligned}
\end{equation*}
Adding the last two equalities proves the claim.
\end{Remark}

\section*{Acknowledgements}
We would like to thank Lars Diening for his valuable and helpful
comments and the fruitful discussions on the topic. The research of the
authors was part of the project C1 of the SFB/TR~71 ``Geometric Partial Differential
Equations''.


\ifx\undefined\bysame
\newcommand{\bysame}{\leavevmode\hbox to3em{\hrulefill}\,}
\fi


\end{document}

%% file: Bild1b.pstex_t
\begin{picture}(0,0)%
\includegraphics{Bild1b.pstex}%
\end{picture}%
\setlength{\unitlength}{3108sp}%
\begingroup\makeatletter\ifx\SetFigFont\undefined%
\gdef\SetFigFont#1#2#3#4#5{%
  \reset@font\fontsize{#1}{#2pt}%
  \fontfamily{#3}\fontseries{#4}\fontshape{#5}%
  \selectfont}%
\fi\endgroup%
\begin{picture}(3111,2515)(-42,-1304)
\put(2161,-601){\makebox(0,0)[lb]{\smash{{\SetFigFont{9}{10.8}{\rmdefault}{\mddefault}{\updefault}{\color[rgb]{0,0,0}$x$}%
}}}}
\put(2341,-1231){\makebox(0,0)[lb]{\smash{{\SetFigFont{9}{10.8}{\rmdefault}{\mddefault}{\updefault}{\color[rgb]{0,0,0}$q(x)$}%
}}}}
\put(1621,884){\makebox(0,0)[lb]{\smash{{\SetFigFont{9}{10.8}{\rmdefault}{\mddefault}{\updefault}{\color[rgb]{0,0,0}$\eta(q)$}%
}}}}
\put(1531,614){\makebox(0,0)[lb]{\smash{{\SetFigFont{9}{10.8}{\rmdefault}{\mddefault}{\updefault}{\color[rgb]{0,0,0}$q$}%
}}}}
\put(1486,-286){\makebox(0,0)[lb]{\smash{{\SetFigFont{9}{10.8}{\rmdefault}{\mddefault}{\updefault}{\color[rgb]{0,0,0}$\Omega$}%
}}}}
\put(1756,-961){\makebox(0,0)[lb]{\smash{{\SetFigFont{9}{10.8}{\rmdefault}{\mddefault}{\updefault}{\color[rgb]{0,0,0}$s(x)$}%
}}}}
\end{picture}%